\newif\ifHAL
\HALfalse 
\HALtrue 

\newif\ifRR
\RRfalse 
\RRtrue

\ifRR
\documentclass[11pt]{article}

\usepackage[centertags]{amsmath}
\usepackage{amsthm}
\usepackage{amssymb}
\usepackage{graphicx}
\usepackage[english]{babel}
\usepackage{minitoc}
\usepackage{url}
\usepackage{color}
\usepackage[dvipdfm,bookmarks=true,colorlinks=true,linkcolor=blue,urlcolor=green,citecolor=red,pdfstartview=FitH,pagebackref]{hyperref}

\theoremstyle{plain}
\newtheorem{theorem}{Theorem}[section]
\newtheorem{corollary}[theorem]{Corollary}
\newtheorem{lemma}[theorem]{Lemma}
\theoremstyle{definition}
\newtheorem{definition}{Definition}[section]
\theoremstyle{remark}
\newtheorem{remark}{Remark}[section]
\numberwithin{equation}{section}

\else
\documentclass{llncs}
\usepackage[centertags]{amsmath}
\usepackage{amssymb}
\usepackage[francais,english]{babel}
\usepackage{graphicx}
\usepackage{url}
\usepackage{color}
\usepackage[dvipdfm,bookmarks=true,colorlinks=true,linkcolor=blue,urlcolor=green,citecolor=red,pdfstartview=FitH,pagebackref]{hyperref}
\fi

\title{Fast learning rates in statistical inference through aggregation}

\author{J.-Y. Audibert$^{1,2}$}

\def\lhs{\text{l.h.s.}}
\def\rhs{\text{r.h.s.}}
\def\wrt{\text{w.r.t.}}

\def\ds1{{\mathbf 1}} 
\def\dsV{{\mathbb V}} 
\def\Var{{\dsV\text{ar}}\,}

\def\E{\mathbb{E}}
\def\N{\mathbb{N}}
\def\P{\mathbb{P}}
\def\Q{\mathbb{Q}}
\def\R{\mathbb{R}}
\def\Rp{\mathbb{R}_+}

\def\begar{$$\begin{array}{lll}}
\def\endar{\end{array}$$}
\def\begarc{$}
\def\endarc{$ }
\def\bigbegar{\begin{eqnarray*}} 
\def\bigendar{\end{eqnarray*}}
\def\lbegar{$$\left\{ \begin{array}{lll}}
\def\rendar{\end{array} \right.$$}
\def\rendarp{\end{array} \right..$$}
\def\begarlab{\begin{equation} \begin{array}{lll} \label}
\def\endarlab{\end{array} \end{equation}}
\def\bigbegarlab{\begin{eqnarray} \label}
\def\bigendarlab{\end{eqnarray}}
\def\lbegarlab{\begin{equation} \left\{ \begin{array}{lll} \label}
\def\rendarlab{\end{array} \right. \end{equation}}
\def\rendarplab{\end{array} \right.. \end{equation}}

\def\A{\mathcal{A}}
\def\B{\mathcal{B}}
\def\C{\mathcal{C}}
\def\D{\mathcal{D}}
\def\calE{\mathcal{E}}

\def\G{\mathcal{G}}
\def\H{\mathcal{H}}
\def\L{\mathcal{L}}
\def\M{\mathcal{M}}

\def\calP{\mathcal{P}}
\def\calQ{\mathcal{Q}}

\def\S{\mathcal{S}}

\def\W{\mathcal{W}}
\def\X{\mathcal{X}}
\def\Y{\mathcal{Y}}
\def\Z{\mathcal{Z}}

\newcommand\und[2]{\underset{#2}{#1}\;}
\newcommand\undc[2]{{#1}_{#2}\;}
\newcommand\jyproof[1]{\begin{proof} See Section \ref{sec:proof#1}. \end{proof}}
\newcommand\secproofth[1]{\subsection{Proof of Theorem \ref{th:#1}} \label{sec:proofth:#1}}

\newcommand\subsecprooflem[1]{\subsubsection{Proof of Lemma \ref{lem:#1}} \label{sec:prooflem:#1}}

\newcommand\fracl[2]{{(#1)}/{#2}}
\newcommand\fracc[2]{{#1}/{#2}}

\newcommand\expe[2]{\undc{\E}{#1\sim#2}}
\newcommand\expec[2]{\undc{\E}{#1\sim#2}}
\newcommand\expecc[2]{\E_{#1}}
\newcommand\expecd[2]{\E_{#2}}

\newcommand\refp[1]{\ref{#1} [p.\pageref{#1}]}

\def\eqdef{\triangleq}
\def\eps{\epsilon}
\def\logeps{\log(\eps^{-1})}
\def\lam{\lambda}

\newcommand\integ[1]{\left\lfloor{#1}\rfloor\right.}

\def\pto{\text{o}}
\def\gdo{\text{O}}

\def\hg{\hat{g}}
\def\hh{\hat{h}}

\def\zun{Z_1^n}
\def\Ezun{\E_{\zun}}

\def\hpi{\hat{\pi}}
\def\hrho{\hat{\rho}}
\def\hmu{\hat{\mu}}
\def\hmua{\hmu_{\text{a}}}

\def\barG{{\bar{\cal G}}}

\def\pirho{\hpi({\rho})}

\def\tildf{\tilde{f}}
\def\tildg{\tilde{g}}

\def\tildR{\tilde{R}}

\def\sign{\text{sign}}

\def\nuv{{\bar \nu}}
\def\sigmav{{\bar \sigma}}

\def\PN{\und{\E}{Z_1^n\sim P^{\otimes n}}}
\def\demi{\frac{1}{2}}
\def\demic{1/2}

\def\bbb{B}
\def\aaa{A}

\def\pa{p_{+}}
\def\pb{p_{-}}
\def\qa{q_{+}}
\def\qb{q_{-}}
\def\ha{h_1}
\def\hb{h_2}
\def\ya{y_1}
\def\yb{y_2}
\def\da{d_{\textnormal{I}}}
\def\db{d_{\textnormal{II}}}
\def\Prestz{P_{\X_1,0}}
\def\Prestm{P_{\X_1,-}}
\def\Prestp{P_{\X_1,+}}
\def\Prestr{P_{\X_1,r}}
\newcommand\Prm[1]{P_{-,#1}}
\newcommand\Prp[1]{P_{+,#1}}
\newcommand\Prr[1]{P_{r,#1}}
\newcommand\Prz[1]{P_{0,#1}}

\def\tpsi{\tilde{\psi}}
\def\begitem{\begin{itemize}}
\def\enditem{\end{itemize}} 
\def\point{\bullet}

\def\tm{\tilde{m}}
\def\tw{\tilde{w}}
\def\tdb{\tilde{\db}}

\def\jyem{\em}

\newcommand\subsectionjy[1]{\ifRR \subsection{#1} \else \section{#1} \fi}
\newcommand\subsubsectionjy[1]{\ifRR \subsubsection{#1} \else \subsection{#1} \fi}

\allowdisplaybreaks[1]

\begin{document}

\date{}
\maketitle

\begin{center}
{ \it
\vspace{-0.5cm}
$^1$Certis - Ecole des Ponts - Paris Est\\
$^2$Willow - ENS/INRIA\\
}
\end{center}
\ifRR

\abstract{
We develop minimax optimal risk bounds for the 
general learning task consisting in predicting as well as the
best function in a reference set $\G$ up to the smallest possible additive term,
called the convergence rate. 
When the reference set is finite and when $n$ denotes the size of 
the training data, we provide minimax convergence rates of the form $C \big(\frac{\log |\G|}{n}\big)^{v}$
with tight evaluation of the positive constant $C$ and with exact $0<v\le 1$, the latter value
depending on the convexity of the loss function and on the level of noise in the output distribution.

The risk upper bounds are based 
on a sequential randomized algorithm, which at each step
concentrates on functions having both low risk and low variance
with respect to the previous step prediction function.
Our analysis puts forward the links between
the probabilistic and worst-case viewpoints, and allows to obtain
risk bounds unachievable with the standard statistical learning approach.
One of the key idea of this work is to use probabilistic inequalities 
with respect to appropriate (Gibbs) distributions on the prediction function space
instead of using them with respect to the distribution generating the data.

The risk lower bounds are based on refinements of the Assouad lemma taking particularly into account 
the properties of the loss function. Our key example to illustrate 
the upper and lower bounds is to consider the $L_q$-regression setting
for which an exhaustive analysis of the convergence rates is given
while $q$ ranges in $[1;+\infty[$.

}


\section{Introduction}

We are given a family $\G$ of functions 
and we want to learn from data a function that predicts as well as the
best function in $\G$ up to some additive term called the convergence rate.
Even when the set $\G$ is finite, this learning task is crucial since
\begin{itemize}
\item any continuous set of
prediction functions can be viewed through its covering nets with respect to ($\wrt$) appropriate (pseudo-)distances
and these nets are generally finite.
\item one way of doing model selection 
among a finite family of submodels
is to cut the training set into two parts, use the first part to learn the best prediction function of
each submodel and use the second part to learn a prediction function which performs as well as the best of 
the prediction functions learned on the first part of the training set. 
\end{itemize}

From this last item, our learning task for finite $\G$ is often referred to as model selection aggregation.
It has two well-known variants. Instead of looking 
for a function predicting as well as the best in $\G$,
these variants want to perform as well as the best convex combination of functions in $\G$ or 
as well as the best linear combination of functions in $\G$. These three aggregation tasks 
are linked in several ways (see \cite{Tsy03} and references within).

Nevertheless, among these learning tasks, model selection aggregation has rare properties.
First, in general an algorithm picking functions in the set $\G$ is not optimal 
(see e.g. \cite[Theorem 2]{Aud07b}, \cite[Theorem 3]{Lec07}, \cite[p.14]{Cat99}). 

This means that the estimator
has to look at an enlarged set of prediction functions.
Secondly, in the statistical community,
the only known optimal algorithms are all based on a Cesaro mean of 
Bayesian estimators (also referred to as progressive mixture rule). 
Thirdly, the proof of their optimality is not 
achieved by the most prominent tool in statistical learning theory:
bounds on the supremum of empirical processes (see \cite{Vap82}, 
and refined works as \cite{Bar02,Kol06,Mas00,bou05} and references within).

The idea of the proof, which comes back to Barron \cite{Bar87}, is based 
on a chain rule and appeared to be successful for least square and entropy losses
\cite{Cat97,Cat99,Bar99,Yan00,Bun05} and for general loss in \cite{Jud06}.

In the online prediction with expert advice setting, without any probabilistic assumption
on the generation of the data, appropriate weighting methods have been
showed to behave as well as the best expert up to a minimax-optimal additive remainder term
(see \cite{Mer98,CesLug06} and references within).
In this worst-case context, amazingly sharp constants have been found (see in particular \cite{Hau98,Ces97,Ces99,Yar04}).
These results are expressed in cumulative loss and can be transposed to model selection aggregation
to the extent that the expected risk of the randomized procedure based on sequential 
predictions is proportional to the expectation of the cumulative loss of the sequential procedure
(see Lemma \ref{lem:rand} for precise statement).

This work presents a sequential algorithm, which iteratively updates a prior
distribution put on the set of prediction functions. 
Contrarily to previously mentioned works, these updates take into account the variance of
the task. As a consequence, posterior distributions concentrate on 
simultaneously low risk functions and functions close to the previously drawn prediction function.
This conservative law is not surprising in view of previous works on
high dimensional statistical tasks, such as wavelet thresholding, shrinkage procedures,
iterative compression schemes (\cite{Aud03b}), iterative feature selection (\cite{Alq05}).

The paper is organized as follows. 
Section \ref{sec:notation} introduces the notation and the existing algorithms.
Section \ref{sec:algo} proposes a unifying setting 
to combine worst-case analysis tight results
and probabilistic tools. It details our sequentially randomized estimator and 
gives a sharp expected risk bound. In Sections \ref{sec:online} and \ref{sec:jud},
we show how to apply our main result under assumptions coming 
respectively from sequential prediction and model selection aggregation.
While all this work concentrates on stating results when the data are 
independent and identically distributed, Section \ref{sec:seqpred} collects 
new results for sequential predictions, i.e. when no probabilistic assumption is made 
and when the data points come one by one (i.e. not in a batch manner).
Section \ref{sec:var} contains algorithms that satisfy sharp standard-style 
generalization error bounds. To the author's knowledge, these bounds are 
not achievable with classical statistical learning approach based on supremum of empirical processes.
Here the main trick is to use probabilistic inequalities 
$\wrt$ appropriate distributions on the prediction function space
instead of using them $\wrt$ the distribution generating the data.
Section \ref{sec:pow} presents an improved bound for $L_q$-regression ($q>1$) when the noise has just a bounded moment 
of order $s~\ge~q$. This last assumption is much weaker than the traditional exponential moment assumption.
Section \ref{sec:assouad} refines Assouad's lemma in order to obtain sharp constants and
to take into account the properties of the loss function of the learning task.
We illustrate our results by providing lower bounds matching the upper bounds obtained in 
the previous sections and by improving significantly the constants in lower bounds 
concerning Vapnik-Cervonenkis classes in classification.
Section \ref{sec:sum} summarizes the contributions of this work and lists some related open problems.

\section{Notation and existing algorithms} \label{sec:notation}

We assume that we observe $n$ pairs $Z_1=(X_1,Y_1),\dots,Z_n=(X_n,Y_n)$ of input-output 
and that each pair has been independently drawn from the same unknown distribution denoted $P$. The input and output space
are denoted respectively $\X$ and $\Y$, so that $P$ is a probability distribution on the product space 
$\Z \eqdef \X\times\Y$. 
The target of a learning algorithm is to predict the output $Y$ associated with an input $X$
for pairs $(X,Y)$ drawn from the distribution $P$. 
In this work, $Z_{n+1}$ will denote a random variable 
independent of the training set $Z_1^n\eqdef(Z_1,\dots,Z_n)$ and with the same distribution $P$.
The quality of a prediction function 
$g:\X\rightarrow\Y$ is measured by the {\jyem risk} (also called expected loss or regret):
        $$R(g) \eqdef \undc{\E}{Z\sim P} L(Z,g),$$
where $L(Z,g)$ assesses the loss of considering the prediction function $g$ on the data $Z\in\Z$.
The symbol $\eqdef$ is used to underline that the equality is a definition.
When there is no ambiguity on the distribution that a random variable has, the expectation
$\wrt$ this distribution will simply be written by indexing the expectation sign $\E$
by the random variable.
For instance, we can write $R(g) \eqdef \undc{\E}{Z} L(Z,g).$ More generally, 
when they are multiple sources of randomness, $\undc{\E}{Z}$ means that we take 
the expectation with respect to the conditional distribution of 
$Z$ knowing all other sources of randomness.

We use $L(Z,g)$ rather than $L[Y,g(X)]$ to underline that our results are not restricted to non-regularized
losses, where we call non-regularized loss a loss that can be written as $\ell[Y,g(X)]$ 
for some function $\ell:\Y\times\Y\rightarrow\R$. 

For any $i\in\{0,\dots,n\}$, the {\jyem cumulative loss} suffered by the prediction function
 $g$ on the first $i$ pairs of input-output, denoted $Z_1^i$ for short, is 
        $$\Sigma_i(g) \eqdef \sum_{j=1}^i L(Z_j,g),$$
where by convention we take $\Sigma_0$ identically equal to zero.
The symbol $\equiv$ is used to underline when a function is identical to a constant (e.g. $\Sigma_0\equiv 0$).
With slight abuse, a symbol denoting a constant function may be used to 
denote the value of this function.

We assume that the set, denoted $\barG$, of all prediction functions
has been equipped with a $\sigma$-algebra. Let $\D$ be the
set of all probability distributions on $\barG$.
By definition, a randomized algorithm produces a prediction
function drawn according to a probability in $\D$.
Let $\calP$ be a set of probability distributions on $\Z$
in which we assume that the true unknown distribution generating the data is.
The learning task is essentially described by the 3-tuple
$(\G,L,P)$ since we look for a possibly randomized estimator (or algorithm) $\hg$ such that 
        \bigbegar
        \und{\sup}{P\in\calP} \Big\{ \Ezun R(\hg_{Z_1^n}) - \und{\min}{g\in\G} \, R(g) \Big\}
        \bigendar
is minimized, where we recall that $R(g) \eqdef \undc{\E}{Z\sim P} L(Z,g).$
To shorten notation, when no confusion can arise, the dependence of $\hg_{Z_1^n}$ $\wrt$ the training sample ${Z_1^n}$
will be dropped and we will simply write $\hg$. 
This means that we use
the same symbol for both the algorithm and the prediction function produced by the algorithm on
a training sample. 

We implicitly assume that the quantities we manipulate are measurable: in particular,
we assume that a prediction function is a measurable function from $\X$ to $\Y$,
the mapping $(x,y,g) \mapsto L[(x,y),g]$ is measurable, 
the estimators considered in our lower bounds are measurable, \dots

The $n$-fold product of a distribution $\mu$, which is the
distribution of a vector consisting in $n$ 
i.i.d. realizations of $\mu$, is denoted $\mu^{\otimes n}$.
For instance the distribution of $(Z_1,\dots,Z_n)$ is $P^{\otimes n}$.

The symbol $C$ will denote some positive constant whose value may differ
from line to line. The set of non-negative real numbers is denoted $\Rp=[0;+\infty[$.
We define $\integ{x}$ as the largest integer $k$ such that $k \le x$.
To shorten notation, any finite sequence $a_1,\dots,a_n$ will 
occasionally be denoted $a_1^n$. For instance, 
the training set is $Z_1^n$.

To handle possibly continuous set $\G$, we consider that $\G$ is a measurable space and that
we have some {\jyem prior distribution $\pi$} on it. 
The set of probability distributions on $\G$ will be denoted $\M$.
The {\jyem Kullback-Leibler divergence} between a distribution $\rho\in\M$ and the prior distribution $\pi$ is
        \begar
        K(\rho,\pi) \eqdef \left\{ \begin{array}{ll}
        \undc{\E}{g\sim\rho} \log\big( \frac{\rho}{\pi}(g) \big) \quad & \text{if } \rho \ll \pi,\\
                + \infty & \text{otherwise}
        \end{array} \right.
        \endar
where $\frac{\rho}{\pi}$ denotes the density of $\rho$ $\wrt$ $\pi$ when it exists (i.e.  $\rho \ll \pi$).
For any $\rho\in\M$, we have $K(\rho,\pi) \ge 0$ and when $\pi$ is the 
uniform distribution on a finite set $\G$, we also have $K(\rho,\pi) \le \log |\G|$.
The Kullback-Leibler divergence satisfies the duality formula (see e.g. \cite[p.10]{Cat03b}): for any 
real-valued measurable function $h$ defined on $\G$,
        \begarlab{eq:legendre}
        \und{\inf}{\rho\in\M} \big\{ \undc{\E}{g\sim\rho} h(g) +
        K(\rho,\pi) \big\} = -\log \expec{g}{\pi} e^{-h(g)}.
        \endarlab
and that the infimum is reached for the {\jyem Gibbs distribution} 
        \begarlab{eq:gibbs}
        \pi_{-h}(dg) \eqdef \frac{e^{-h(g)}}{\undc{\E}{g'\sim\pi} e^{-h(g')}} \cdot \pi(dg).
        \endarlab
Intuitively, the Gibbs distribution $\pi_{-h}$ concentrates on prediction functions $g$
that are close to minimizing the function $h:\G\rightarrow \R$.

For any $\rho\in\M$, $\undc{\E}{g\sim\rho} g: x \mapsto \undc{\E}{g\sim\rho} g(x) = \int g(x) \rho(dg)$
is called a mixture of prediction functions. 
When $\G$ is finite, a mixture is simply a convex combination.
Throughout this work, whenever we consider mixtures of prediction functions,
we implicitly assume that $\undc{\E}{g\sim\rho} g(x)$ belongs to $\Y$ for any $x$ so that the mixture 
is a prediction function. This is typically the case when $\Y$ is an interval of $\R$.

We will say that the loss function is convex when the function $g\mapsto L(z,g)$ is convex 
for any $z\in\Z$, equivalently $L(z,\undc{\E}{g\sim\rho} g) \le \undc{\E}{g\sim\rho} L(z,g)$ for any $\rho\in\M$ and $z\in\Z$.
In this work, we do not assume the loss function to be convex except when it is explicitly mentioned.


The algorithm used to prove optimal convergence rates for several different losses 
(see e.g. \cite{Cat97,Cat99,Bar99,Bla99,Yan00,Bun05,Jud06}) is the following:

\vspace{0.2cm}
\noindent {\bf Algorithm $\aaa$:}
Let $\lam>0$.
Predict according to
         $\frac{1}{n+1} \sum_{i=0}^n \expec{g}{\pi_{-\lam \Sigma_i}} g$, 
where we recall that $\Sigma_i$ maps a function $g\in\G$ to its cumulative loss up to time $i$. 
\vspace{0.2cm}

In other words, for a new input $x$, the prediction of the output given by Algorithm $\aaa$ is
        $
        \frac{1}{n+1} \sum_{i=0}^n \frac{\int g(x) e^{-\lam \Sigma_i(g)} \pi(dg)}{\int e^{-\lam \Sigma_i(g)} \pi(dg)}.
        $
Algorithm $\aaa$ has also been used with the classification loss. For this non-convex loss, it has the same
properties as the empirical risk minimizer on $\G$ (\cite{Lec07a,Lec07b}).
To give the optimal convergence rate, the parameter $\lam$ and the distribution $\pi$ should be appropriately chosen.
When $\G$ is finite, the estimator belongs to the convex hull of the set $\G$. 

From Vovk, Haussler, Kivinen and Warmuth works (\cite{Vov90,Hau98,Vov98})
and the link between cumulative loss in online setting and expected risk in the
batch setting (see later Lemma \ref{lem:rand}), an ``optimal'' algorithm is:

\vspace{0.2cm}
\noindent {\bf Algorithm $\bbb$:}
Let $\lam>0$. For any $i\in\{0,\dots,n\}$, let $\hh_i$ be a prediction function such that 
        \bigbegar
        \forall \, z \in\Z \qquad
                L(z,\hh_i) \le -\frac{1}{\lam} \log \expe{g}{\pi_{-\lam \Sigma_i}} e^{-\lam L(z,g)}.
        \bigendar
If one of the $\hh_i$ does not exist, the algorithm is said to fail. Otherwise it
predicts according to $\frac{1}{n+1} \sum_{i=0}^n \hh_i$.
\vspace{0.2cm}

In particular, for appropriate $\lam>0$, this algorithm does not fail when
the loss function is the square loss (i.e. $L(z,g)=[y-g(x)]^2$) and when the output space is bounded.
Algorithm $\bbb$ is based on the same Gibbs distribution $\pi_{-\lam \Sigma_i}$ as Algorithm $\aaa$.
Besides, in \cite[Example 3.13]{Hau98}, it is shown that Algorithm $\aaa$ is not in general
a particular case of Algorithm $\bbb$, and that Algorithm $\bbb$ will not generally produce
a prediction function in the convex hull of $\G$ unlike Algorithm $\aaa$.
In Sections \ref{sec:online} and \ref{sec:jud}, we will see how both algorithms
are connected to the SeqRand algorithm presented in the next section.

\section{The algorithm and its generalization error bound} \label{sec:algo}

The aim of this section is to build an algorithm with the best 
possible minimax convergence rate. The algorithm relies on the following central 
condition for which we recall that $\G$ is a subset of the set $\barG$
of all prediction functions and that $\M$ and $\D$ are the sets of all probability distributions on respectively 
$\G$ and~$\barG$.

For any $\lam>0$, let $\delta_\lam$ be a real-valued function defined on $\Z\times\G\times\barG$
that satisfies the following inequality, which will be referred to as the \emph{variance inequality} \label{varcond}
        \begin{multline*}
        \forall \, \rho\in\M \quad \exists \, \pirho \in \D \quad\\
        \qquad
                \und{\sup}{P\in\calP} \bigg\{ 
                \expec{Z}{P} \expec{g'}{\pirho} \log \expec{g}{\rho}
                e^{\lam \big[L(Z,g') - L(Z,g) - \delta_\lam(Z,g,g')\big]} \bigg\} \le 0.
        \end{multline*}

The variance inequality is our probabilistic version of the generic algorithm condition in
the online prediction setting (see \cite[proof of Theorem 1]{Vov90} or more explicitly in \cite[p.11]{Hau98}),
in which we added the variance function $\delta_\lam$. 
Our results will be all the sharper as this variance function is small. 
To make the variance inequality more readable, let us say for the moment that
\begitem
\item
without any assumption on $\calP$, for several usual ``strongly'' convex loss functions, 
we may take $\delta_\lam \equiv 0$ provided that $\lam$ is a small enough constant (see Section \ref{sec:online}).
\item the variance inequality can be seen as a ``small expectation'' inequality. 
The usual viewpoint is to control the quantity $L(Z,g)$ by its expectation $\wrt$ $Z$
and a variance term. Here, roughly, $L(Z,g)$ is mainly controlled by $L(Z,g')$ 
where $g'$ is appropriately chosen through the choice of $\pirho$,
plus the additive term $\delta_\lam$. By definition this additive term does not depend on 
the particular probability distribution generating the data and leads to empirical compensation.
\item 
in the examples we will be interested in throughout this work, 
$\pirho$ will be either equal to $\rho$ or
to a Dirac distribution on some function, which is \emph{not necessarily in $\G$}.
\item for any loss function $L$, any set $\calP$ and any $\lam>0$, 
one may choose $\delta_\lam(Z,g,g')=\frac{\lam}{2} \big[L(Z,g)-L(Z,g')\big]^2$
(see Section \ref{sec:var}).
\enditem


Our results concern the sequentially randomized algorithm described in Figure~\ref{fig:seqr}, which for sake of shortness
we will call the SeqRand algorithm.

\begin{figure} 
\hspace*{0.5cm}\hbox{\raisebox{0.4em}{\vrule depth 0pt height 0.4pt width 12cm}}
\begin{enumerate}
\item[] Input: $\lam>0$ and $\pi$ a distribution on the set $\G$. 
\item Define $\hrho_0 \eqdef \hpi(\pi)$ in the sense of the variance inequality (p.\pageref{varcond})
and draw a function $\hg_{0}$ according to this distribution.
Let $S_0(g)=0$ for any $g\in\G$. 
\item For any $i\in\{1,\dots,n\}$, iteratively define 
        \begarlab{eq:defSi}
        S_{i}(g) \eqdef S_{i-1}(g) + L(Z_{i},g) + \delta_\lam(Z_{i},g,\hg_{i-1}) \quad \text{for any } g\in\G.
        \endarlab
and
        \begar
        \hrho_i \eqdef \hpi(\pi_{-\lam S_i}) \qquad\text{in the sense of the variance inequality (p.\pageref{varcond})}
        \endar
and draw a function $\hg_i$ according to the distribution $\hrho_i$.
\item Predict with a function drawn according to the uniform distribution on the finite set $\{\hg_0,\dots,\hg_n\}$.\\
Conditionally to the training set, the distribution of the output prediction function will 
be denoted $\hmu$.
\end{enumerate}
\hspace*{0.5cm}\hbox{\raisebox{0.4em}{\vrule depth 0pt height 0.4pt width 12cm}}
\caption{The SeqRand algorithm} \label{fig:seqr}
\end{figure}

\begin{remark}
When $\delta_\lam(Z,g,g')$ does not depend on $g$, we recover a more standard-style algorithm
to the extent that we then have $\pi_{-\lam S_i}= \pi_{-\lam \Sigma_i}$.
Precisely our algorithm becomes the randomized version of Algorithm~$\aaa$.
When $\delta_\lam(Z,g,g')$ depends on $g$, the posterior distributions
tend to concentrate on functions having small risk and
small variance term. 
In Section \ref{sec:var}, we will take 
        $\delta_\lam(Z,g,g')=\frac{\lam}{2} \big[L(Z,g)-L(Z,g')\big]^2$.
This choice implies a conservative mechanism: roughly, with high probability,
among functions having low cumulative risk $\Sigma_i$,
$\hg_i$ will be chosen close to $\hg_{i-1}$.
\end{remark}

For any $i\in\{0,\dots,n\}$, the quantities
$S_i$, $\hrho_i$ and $\hg_i$ depend on the training data only through $Z_1^i$, where 
we recall that $Z_1^i$ denotes $(Z_1,\dots,Z_i)$.
Besides they are also random to the extent that they depend on the draws of the functions
$\hg_0,\dots,\hg_{i-1}$.

The SeqRand algorithm produces a prediction function which has three causes of randomness:
the training data, the way $\hg_i$ is obtained (step 2) and the uniform draw (step 3).
For fixed $Z_1^i$ (i.e. conditional to $Z_1^i$), let $\Omega_i$ 
denote the joint distribution of $\hg_0^i = (\hg_0,\dots,\hg_i)$.
The randomizing distribution $\hmu$ of the output prediction function by SeqRand is the distribution
on $\barG$ corresponding to the last two causes of randomness. From the previous definitions,
for any function $h:\barG\rightarrow \R$, we have 
        \begarc
        \expec{g}{\hmu} h(g) = \expec{\hg_0^n}{\Omega_n} \frac{1}{n+1} \sum_{i=0}^{n} h(\hg_i). 
        \endarc
Our main upper bound controls the expected risk $\Ezun \expec{g}{\hmu} R(g)$ 
of the SeqRand procedure.

\begin{theorem} \label{th:1}
Let $\Delta_\lam(g,g') \eqdef \undc{\E}{Z\sim P} \delta_{\lam}(Z,g,g')$ for $g \in G$ and $g'\in \barG$,
where we recall that $\delta_{\lam}$ is a function satisfying the variance inequality (see p.\pageref{varcond}). 
The expected risk of the SeqRand algorithm satisfies
        \begarlab{eq:1gen}
        \Ezun \expec{g'}{\hmu} R(g')
                \le \und{\min}{\rho\in\M} \bigg\{ 
                \expec{g}{\rho} R(g) 
                + \expec{g}{\rho} \Ezun \expec{g'}{\hmu} \Delta_\lam(g,g')
                + \frac{K(\rho,\pi)}{\lam(n+1)} \bigg\}
        \endarlab
In particular, when $\G$ is finite and 
when the loss function $L$ and the set $\calP$ are such that $\delta_\lam \equiv 0$, 
by taking $\pi$ uniform on $\G$, we get
        \begarlab{eq:1app2}
        \Ezun \expec{g}{\hmu} R(g) 
                \le \und{\min}{\G} R
                + \frac{\log |\G|}{\lam(n+1)} 
        \endarlab
\end{theorem}

\begin{proof}
Let $\calE$ denote the expected risk of the SeqRand algorithm:
        \begar 
        \calE \eqdef \Ezun \expec{g}{\hmu} R(g)
                 = \frac{1}{n+1} \sum_{i=0}^{n} \E_{Z_1^i} \expec{\hg_0^i}{\Omega_i} R(\hg_i).
        \endar
We recall that $Z_{n+1}$ is a random variable 
independent of the training set $Z_1^n$ and with the same distribution $P$.
Let $S_{n+1}$ be defined by \eqref{eq:defSi} for $i=n+1$.
To shorten formulae, let $\hpi_i \eqdef \pi_{-\lam S_i}$ so that by definition we have
$\hrho_i=\hpi(\hpi_i)$. The variance inequality implies that
        \begar
        \expec{g'}{\pirho} R(g') \le - \frac{1}{\lam} 
                \E_Z \expec{g'}{\pirho} \log \expec{g}{\rho} e^{-\lam [L(Z,g)+\delta_\lam(Z,g,g')]}.
        \endar
So for any $i\in\{0,\dots,n\},$ for fixed $\hg_0^{i-1}=(\hg_0,\dots,\hg_{i-1})$
and fixed $Z_1^i$, we have
        \begar
        \expe{g'}{\hrho_i} R(g') \le - \frac{1}{\lam} 
                \E_{Z_{i+1}} \expec{g'}{\hrho_i} \log \expec{g}{\hpi_i} 
                e^{-\lam [L(Z_{i+1},g)+\delta_\lam(Z_{i+1},g,g')]}
        \endar
Taking the expectations $\wrt$ $(Z_1^i,\hg_0^{i-1})$, we get
        \begar
        \E_{Z_1^i} \expecc{\hg_0^i}{\Omega_i} R(\hg_i)
                & = & \E_{Z_1^i} \expecc{\hg_0^{i-1}}{\Omega_{i-1}} \expec{g'}{\hrho_i} R(g')\\
                & \le & - \frac{1}{\lam} \E_{Z_1^{i+1}}\expecc{\hg_0^i}{\Omega_i} 
                \log \expec{g}{\hpi_i} e^{-\lam [L(Z_{i+1},g)+\delta_\lam(Z_{i+1},g,\hg_i)]}.
        \endar
Consequently, by the chain rule (i.e. cancellation in the sum of logarithmic terms; \cite{Bar87}) 
and by intensive use of Fubini's theorem, we get        
        \begar
        \calE = \frac{1}{n+1} \sum_{i=0}^{n} \E_{Z_1^i} \expecc{\hg_0^i}{\Omega_i} R(\hg_i)\\
        \, \le - \frac{1}{\lam(n+1)} \sum_{i=0}^{n} 
                \E_{Z_1^{i+1}} \expecc{\hg_0^i}{\Omega_i} \log \expec{g}{\hpi_i} e^{-\lam 
                [L(Z_{i+1},g)+\delta_\lam(Z_{i+1},g,\hg_i)]}\\
        \, = - \frac{1}{\lam(n+1)} 
                \E_{Z_1^{n+1}} \expecc{\hg_0^n}{\Omega_n} \sum_{i=0}^{n} 
                \log \expe{g}{\hpi_i} e^{-\lam [L(Z_{i+1},g)+\delta_\lam(Z_{i+1},g,\hg_i)]}\\
        \, = - \frac{1}{\lam(n+1)} 
                \E_{Z_1^{n+1}} \expecc{\hg_0^n}{\Omega_n} \sum_{i=0}^{n} 
                \log \left( \frac{\expec{g}{\pi} e^{-\lam S_{i+1}(g)}}
                {\expec{g}{\pi} e^{-\lam S_{i}(g)}} \right)\\
        \, = - \frac{1}{\lam(n+1)} 
                \E_{Z_1^{n+1}} \expecc{\hg_0^n}{\Omega_n} 
                \log \left( \frac{\expec{g}{\pi} e^{-\lam S_{n+1}(g)}}
                {\expec{g}{\pi} e^{-\lam S_{0}(g)}} \right)\\
        \, = - \frac{1}{\lam(n+1)} 
                \E_{Z_1^{n+1}} \expecc{\hg_0^n}{\Omega_n} 
                \log \expec{g}{\pi} e^{-\lam S_{n+1}(g)}\\
        \endar
Now from the following lemma, 
we obtain
        \begar
        \calE 
        & \le & - \frac{1}{\lam(n+1)} \log \expec{g}{\pi} e^{-\lam \expecc{Z_1^{n+1}}{P^{\otimes (n+1)}} 
                \expecc{\hg_0^n}{\Omega_n} S_{n+1}(g)}\\
        & = & - \frac{1}{\lam(n+1)} \log \expec{g}{\pi} e^{-\lam \big[ (n+1) R(g) + 
                        \Ezun \expecc{\hg_0^n}{\Omega_n} 
                \sum_{i=0}^{n} \Delta_\lam(g,\hg_i) \big]}\\
        & = & \und{\min}{\rho\in\M} \left\{ 
                \expec{g}{\rho} R(g) 
                + \expec{g}{\rho} \Ezun \expecc{\hg_0^n}{\Omega_n} \frac{\sum_{i=0}^{n} \Delta_\lam(g,\hg_i)}{n+1} 
                + \frac{K(\rho,\pi)}{\lam(n+1)} \right\}.
        \endar  
        
\begin{lemma} \label{eq:concpart}
Let $\W$ be a real-valued measurable function defined on a product space $\A_1\times\A_2$ and
let $\mu_1$ and $\mu_2$ be probability distributions on respectively $\A_1$ and $\A_2$ such
that $\expec{a_1}{\mu_1} \log \expec{a_2}{\mu_2} e^{-\W(a_1,a_2)} < +\infty$. We have
        \begar
        - \expec{a_1}{\mu_1} \log \expec{a_2}{\mu_2} e^{-\W(a_1,a_2)} \le - \log \expec{a_2}{\mu_2} e^{-\expec{a_1}{\mu_1} \W(a_1,a_2)}.
        \endar
\end{lemma}
\begin{proof}
By using twice \eqref{eq:legendre} and Fubini's theorem, we have
        \begar
        - \expecc{a_1}{\mu_1} \log \expec{a_2}{\mu_2} e^{-\W(a_1,a_2)} 
                & = & \expecc{a_1}{\mu_1} \und{\inf}{\rho} \big\{ \expec{a_2}{\rho} \W(a_1,a_2) + K(\rho,\mu_2) \big\} \\
        & \le & \und{\inf}{\rho} \expecc{a_1}{\mu_1} \big\{ \expec{a_2}{\rho} \W(a_1,a_2) + K(\rho,\mu_2) \big\} \\
        & = & - \log \expec{a_2}{\mu_2} e^{-\expecc{a_1}{\mu_1} \W(a_1,a_2)}.
        \endar
\end{proof}
Inequality \eqref{eq:1app2} is a direct consequence of \eqref{eq:1gen}.
\end{proof}

Theorem \ref{th:1} bounds the expected risk of a randomized procedure, where the
expectation is taken $\wrt$ both the training set distribution and the randomizing 
distribution. From the following lemma, for convex loss functions, \eqref{eq:1app2} implies 
        \begarlab{eq:1app3}
        \undc{\E}{Z_1^n} R( \expec{g}{\hmu} g )
                 \le \und{\min}{\G} R
                + \frac{\log |\G|}{\lam(n+1)},
        \endarlab
where we recall that $\hmu$ is the randomizing distribution of 
the SeqRand algorithm and $\lam$ is a parameter
whose typical value is the largest $\lam>0$ such that $\delta_\lam\equiv 0$.

\begin{lemma} \label{le:rand}
For convex loss functions, the doubly expected risk of a randomized algorithm is greater than
the expected risk of the deterministic version of the randomized algorithm, i.e. if
$\hrho$ denotes the randomizing distribution, we have
        $$\E_{Z_1^n} R(\expec{g}{\hrho} g) \le \E_{Z_1^n} \expec{g}{\hrho} R(g).$$ 
\end{lemma}
\begin{proof}
The result is a direct consequence of Jensen's inequality.
\end{proof}

In \cite{Ces97}, the authors rely on worst-case analysis to recover standard-style statistical results 
such as Vapnik's bounds \cite{Vap95}.
Theorem \ref{th:1} can be seen as a complement to this pioneering work. 
Inequality \eqref{eq:1app3} is the model selection bound 
that is well-known for least square regression and entropy loss,
and that has been recently proved for general losses in \cite{Jud06}.

Let us discuss the generalized form of the result.
The $\rhs$ of \eqref{eq:1gen} is a classical regularized risk, which appears naturally in
the PAC-Bayesian approach (see e.g. \cite{Cat02,Cat03b,Aud04,Zha06}).
An advantage of stating the result this way is to be able 
to deal with uncountable infinite $\G$.
Even when $\G$ is countable, this formulation has some benefit to the extent that
for any measurable function $h: \G \rightarrow \R$,
        ${\min}_{\rho\in\M} \{ 
                \expec{g}{\rho} h(g) + K(\rho,\pi) \}
                \le \und{\min}{g\in\G} \{ 
                h(g) + \log \pi^{-1}(g) \}.$

Our generalization error bounds depend on two quantities $\lam$ and $\pi$ which 
are the parameters of our algorithm.  
Their choice depends on the precise setting.
Nevertheless, when $\G$ is finite and with no particular structure a priori, 
a natural choice for $\pi$ is the uniform distribution on $\G$. 

Once the distribution $\pi$ is fixed, an appropriate choice for the parameter $\lam$ is the minimizer of 
the $\rhs$ of \eqref{eq:1gen}. This minimizer is unknown by the statistician, and it is an open
problem to adaptively choose $\lam$ close to it in this general context. Solutions for specific sequential prediction frameworks 
 are known (see \cite[Section 2]{AueCesGen02} and \cite[Lemma 3]{Ces07}). They are 
 based on incremental updating of $\lam$. In appendix, one may found a slight improvement of 
the argument used in the forementioned works, based on Lemma \ref{lem:partition}.

\section{Link with sequential prediction} \label{sec:online}

This section aims at providing examples for which the variance inequality (p.\pageref{varcond}) holds, at 
stating results coming from the online learning community in our batch setting (Section \ref{sec:onl2bat}), and
at providing new results for the sequential prediction setting in which
no probabilistic assumption is made on the way the data are generated (Section \ref{sec:seqpred}).

\subsection{From online to batch} \label{sec:onl2bat}

In \cite{Vov90,Hau98,Vov98}, the loss function is assumed to satisfy:
there are positive numbers $\eta$ and $c$ such that 
        \begarlab{eq:vov}
        \forall \, \rho \in \M \quad \exists \, g_\rho : \X \rightarrow \Y 
                \quad \forall \, x\in\X \quad \forall \, y\in\Y\\
        \qquad\qquad\qquad\qquad
                L[(x,y),g_\rho] \le -\frac{c}{\eta} \log \expec{g}{\rho} e^{-\eta L[(x,y),g]}
        \endarlab

\begin{remark} \label{rem:expcon}
If $g\mapsto e^{-\eta L(z,g)}$ is concave, then \eqref{eq:vov} holds for $c=1$
(and one may take $g_\rho = \expec{g}{\rho} g$).
\end{remark}

Assumption \eqref{eq:vov} implies that the variance inequality is satisfied both 
for $\lam=\eta$ and $\delta_\lam(Z,g,g') = (1-1/c) L(Z,g')$ 
and for $\lam=\eta/c$ and $\delta_\lam(Z,g,g') = (c-1) L(Z,g)$, and we may take
in both cases $\pirho$ as the Dirac distribution at $g_\rho$. 
This leads to the \emph{same} procedure that
is described in the following straightforward corollary of Theorem \ref{th:1}.

\begin{corollary} \label{th:ex1}
Let $g_{\pi_{- \eta \Sigma_i}}$ be defined in the sense of \eqref{eq:vov} (for $\rho=\pi_{- \eta \Sigma_i}$).
Consider the algorithm which predicts by drawing 
a function in $\{g_{\pi_{- \eta \Sigma_0}},\dots,g_{\pi_{- \eta \Sigma_n}}\}$
according to the uniform distribution. 
Under Assumption \eqref{eq:vov},
its expected risk $\Ezun \frac{1}{n+1} \sum_{i=0}^n R(g_{\pi_{- \eta \Sigma_i}})$ 
is upper bounded by
        \begarlab{eq:ex1}
         c \, \und{\min}{\rho\in\M} 
                \big\{ \expec{g}{\rho} R(g) + \frac{K(\rho,\pi)}{\eta(n+1)} \big\}.
        \endarlab
\end{corollary}

This result is not surprising in view of the following two results. The first one
comes from worst-case analysis in sequential prediction.

\begin{theorem}[Haussler et al. \cite{Hau98}, Theorem 3.8] \label{th:hau}
Let $\G$ be countable. For any $g\in\G$, let $\Sigma_i(g)=\sum_{j=1}^{i} L(Z_j,g)$ (still) 
denote the cumulative loss up to time $i$
of the expert which always predicts according to function $g$. 
Under Assumption \eqref{eq:vov}, the cumulative loss on $Z_1^n$ of the strategy in which the prediction at time $i$
is done according to $g_{\pi_{- \eta \Sigma_{i-1}}}$ in the sense of \eqref{eq:vov}
(for $\rho=\pi_{- \eta \Sigma_{i-1}}$) is bounded by 
        \begarlab{eq:thhau}
        \inf_{g\in\G} \{ c \Sigma_n(g) + \frac{c}{\eta} \log \pi^{-1}(g) \}.
        \endarlab
\end{theorem}

The second result shows how the previous bound can be transposed 
into our model selection context by the following lemma.

\begin{lemma} \label{lem:rand}
Let $\A$ be a learning algorithm which produces the prediction function $\A(Z_1^i)$ at time $i+1$, i.e. from the data 
$Z_1^i=(Z_1,\dots,Z_i)$.
Let $\L$ be the randomized algorithm which produces a prediction function $\L(Z_1^n)$ 
drawn according to the uniform
distribution on $\{\A(\emptyset),\A(Z_1),\dots,\A(Z_1^n)\}$. 
The (doubly) expected risk of $\L$ is equal to $\frac{1}{n+1}$ times 
the expectation of the cumulative loss of $\A$ on the sequence $Z_1,\dots,Z_{n+1}$.
\end{lemma}
\begin{proof}
By Fubini's theorem, we have 
        \begar
        \E R[\L(Z_1^n)] & = & \frac{1}{n+1} \sum_{i=0}^n \Ezun R[\A(Z_1^i)] \\
        & = & \frac{1}{n+1} \sum_{i=0}^n \E_{Z_1^{i+1}} L[Z_{i+1},\A(Z_1^i)]\\
        & = & \frac{1}{n+1} \expecc{Z_1^{n+1}}{P} \sum_{i=0}^n L[Z_{i+1},\A(Z_1^i)].
        \endar
\end{proof}

For any $\eta>0$, let $c(\eta)$ denote the infimum of the $c$ for which \eqref{eq:vov} holds.
Under weak assumptions, Vovk (\cite{Vov98}) proved that the infimum exists and studied the
behavior of $c(\eta)$ and $a(\eta)=c(\eta)/\eta$, which are key quantities of \eqref{eq:ex1} and \eqref{eq:thhau}.
Under weak assumptions, and in particular in the examples given in Table \ref{tab:table},
the optimal constants in \eqref{eq:thhau} are $c(\eta)$ and $a(\eta)$ (\cite[Theorem 1]{Vov98}) and
we have $c(\eta) \ge 1$, $\eta \mapsto c(\eta)$ nondecreasing and $\eta \mapsto a(\eta)$ nonincreasing.
From these last properties, we understand the trade-off which occurs to choose the optimal $\eta$.

        \begin{table}[!ht] 
        \footnotesize
        \centering
            \begin{tabular}{|c|c|c|c|}
            \hline
            & Output space & Loss L(Z,g) & $c(\eta)$\\
            \hline
            Entropy loss & $\Y=[0;1]$ & $Y \log\big(\frac{Y}{g(X)}\big)$ & $c(\eta)=1$ if $\eta \le 1$\\
            \cite[Example 4.3]{Hau98} & & $+ (1-Y) \log\big( \frac{1-Y}{1-g(X)}\big)$ & $c(\eta)=\infty$ if $\eta > 1$\\
            \hline
            Absolute loss game & $\Y=[0;1]$  & $|Y-g(X)|$ & $\frac{\eta}{2 \log[2/(1+e^{-\eta})]}$ \\
            \cite[Section 4.2]{Hau98} & & & = $1 + \eta/4 + \pto(\eta)$\\
            \hline
            Square loss & $\Y=[-B,B]$ & $[Y-g(X)]^2$ & $c(\eta)=1$ if $\eta \le 1/(2B^2)$\\
            \cite[Example 4.4]{Hau98} & & & $c(\eta)=+\infty$ if $\eta > 1/(2B^2)$\\
            \hline
            $L_q$-loss & $\Y=[-B,B]$ & $|Y-g(X)|^q$ & $c(\eta)=1\qquad\qquad$ \\
            (see p. \pageref{sec:power}) &  & $q> 1$ & \qquad if $\eta \le \frac{q-1}{q B^q} (1 \wedge 2^{2-q})$\\
            \hline
            \end{tabular}
        \caption{Value of $c(\eta)$ for different loss functions. Here $B$ denotes a positive real.} \label{tab:table}
        \end{table}

Table \ref{tab:table} specifies \eqref{eq:ex1} in different well-known learning tasks. For instance, 
for bounded least square regression (i.e. when $|Y|\le B$ for some $B>0$), the generalization error of the algorithm described in Corollary \ref{th:ex1} 
when $\eta=1/(2B^2)$ is upper bounded by
        \begarlab{eq:lsup}
        {\min}_{\rho\in\M} 
                \big\{ \expec{g}{\rho} R(g) + 2B^2\frac{K(\rho,\pi)}{n+1} \big\}.
        \endarlab
The constant appearing in front of the Kullback-Leibler divergence 
is much smaller than the ones obtained in unbounded regression setting even with gaussian noise
and bounded regression function (see \cite{Bun05,Jud06} and \cite[p.87]{Cat02}). 
The differences between these results partly comes from the absence of boundedness assumptions
on the output and from the weighted average used in the aforementioned works.
Indeed the weighted average prediction function, i.e. $\expec{g}{\rho} g$, does not 
satisfy \eqref{eq:vov} for $c=1$ and $\eta=1/(2B^2)$ as was pointed out in \cite[Example 3.13]{Hau98}.
Nevertheless, it satisfies \eqref{eq:vov} for $c=1$ and $\eta \le 1/(8B^2)$ (by using the concavity of $x\mapsto e^{-x^2}$ 
on $[-1/\sqrt 2;1/\sqrt 2]$ and Remark \ref{rem:expcon}), 
which leads to similar but weaker bound (see \eqref{eq:ex1}).

\paragraph{Case of the $L_q$-losses.} \label{sec:power}

To deal with these losses, we need the following slight generalization of the result 
given in Appendix A of \cite{Kiv99L}. 

\begin{theorem} \label{th:kivwar}
Let $\Y=[a;b]$. 
We consider a non-regularized loss function, i.e.
a loss function such that $L(Z,g)=\ell[Y,g(X)]$ for any
$Z=(X,Y)\in\Z$ and some function $\ell:\Y\times\Y\rightarrow\R$. 
For any $y\in\Y$, let $\ell_y$ be the function 
$\big[ y' \mapsto \ell(y,y') \big]$.
If for any $y\in\Y$
\begin{itemize}
\item $\ell_y$ is continuous on $\Y$
\item $\ell_y$ decreases on $[a;y]$, increases on $[y;b]$ and $\ell_y(y)=0$
\item $\ell_y$ is twice differentiable on the open set $(a;y) \cup (y;b)$,
\end{itemize}
then \eqref{eq:vov} is satisfied for $c=1$ and 
        \begarlab{eq:etacond}
        \eta \le \und{\inf}{a\le y_1 < y < y_2\le b}
                \frac{\ell_{y_1}'(y)\ell_{y_2}''(y)-\ell_{y_1}''(y)\ell_{y_2}'(y)}
                {\ell_{y_1}'(y)[\ell_{y_2}'(y)]^2-[\ell_{y_1}'(y)]^2 \ell_{y_2}'(y)},
        \endarlab
where the infimum is taken $\wrt$ $y_1,y$ and $y_2$.
\end{theorem}

\jyproof{th:kivwar}

\begin{remark}
This result simplifies the original one to the extent that 
$\ell_y$ does not need to be twice differentiable at point $y$
and the range of values for $y$ in the infimum is $(y_1;y_2)$ instead of $(a;b)$.
\end{remark}

\begin{corollary} \label{cor:seqlq}
For the $L_q$-loss, when $\Y=[-B;B]$ for some $B>0$, condition \eqref{eq:vov} is satisfied for $c=1$ and 
        \begar
        \eta \le \frac{q-1}{q B^q} (1 \wedge 2^{2-q})
        \endar
\end{corollary}

\begin{proof}
We apply Theorem \ref{th:kivwar}.
By simple computations, the $\rhs$ of \eqref{eq:etacond} is 
        \begar
        & \und{\inf}{-B\le y_1 < y < y_2\le B} \frac{(q-1)(y_2-y_1)}{q(y-y_1)(y_2-y)[(y-y_1)^{q-1}+(y_2-y)^{q-1}]} \\
        = & \frac{q-1}{q(2B)^q} \und{\inf}{0< t <1} \frac{1}{t(1-t)[t^{q-1}+(1-t)^{q-1}]}\\
        \endar
For $1<q\le 2$, the infimum is reached for $t=1/2$ and \eqref{eq:etacond} can be written as
        $
        \eta \le \frac{q-1}{q B^q}.
        $
For $q\ge 2$, since the previous infimum is larger than $\inf_{0<t<1} \frac{1}{t(1-t)}=4$, 
\eqref{eq:etacond} is satisfied at least when $\eta \le \frac{4(q-1)}{q\;(2B)^{q}}.$
\end{proof}

\subsection{Sequential prediction} \label{sec:seqpred}

First note that using Corollary \ref{cor:seqlq} and Theorem \ref{th:hau}, we obtain 
a new result concerning sequential prediction for $L_q$ loss. Nevertheless this
result is not due to our approach but on a refinement of the argument in \cite[Appendix A]{Kiv99L}.
In this section, we will rather concentrate on giving results for sequential prediction
coming from the arguments underlying Theorem \ref{th:1}.

In the online setting, the data points come one by one and there is no probabilistic assumption on 
the way they are generated. In this case, one should modify the definition of the variance function into:
for any $\lam>0$, let $\delta_\lam$ be a real-valued function defined on $\Z\times\G\times\barG$
that satisfies the following \emph{online variance inequality} \label{ovarcond}
        \begin{multline*}
        \forall \, \rho\in\M \quad \exists \, \pirho \in \D \quad \forall \, Z\in\Z \\
        \qquad\qquad
                \expec{g'}{\pirho} \log \expec{g}{\rho} e^{\lam \big[L(Z,g') - L(Z,g) - \delta_\lam(Z,g,g')\big]} \le 0.
        \end{multline*}
The only difference with the variance inequality (defined in p.\pageref{varcond}) is the removal of the expectation with respect to $Z$.
Naturally if $\delta_\lam$ satisfies the online variance inequality, then it satisfies the variance inequality.
The online version of the SeqRand algorithm is described in Figure~\ref{fig:oseqr}.
It satisfies the following theorem whose proof follows the same line as the one of Theorem \ref{th:1}.

\begin{figure} [h!]
\hspace*{0.5cm}\hbox{\raisebox{0.4em}{\vrule depth 0pt height 0.4pt width 12cm}}
\begin{enumerate}
\item[] Input: $\lam>0$ and $\pi$ a distribution on the set $\G$. 
\item Define $\hrho_0 \eqdef \hpi(\pi)$ in the sense of the online variance inequality (p.\pageref{ovarcond})
and draw a function $\hg_{0}$ according to this distribution. For data $Z_1$, predict according to $\hg_{0}$.
Let $S_0(g)=0$ for any $g\in\G$. 
\item For any $i\in\{1,\dots,n-1\}$, define 
        $$
        S_{i}(g) \eqdef S_{i-1}(g) + L(Z_{i},g) + \delta_\lam(Z_{i},g,\hg_{i-1}) \quad \text{for any } g\in\G.
        $$
and
        $$
        \hrho_i \eqdef \hpi(\pi_{-\lam S_i}) \qquad\text{in the sense of the online variance inequality (p.\pageref{ovarcond})}
        $$
and draw a function $\hg_i$ according to the distribution $\hrho_i$.
For data $Z_{i+1}$, predict according to $\hg_{i}$.
\end{enumerate}
\hspace*{0.5cm}\hbox{\raisebox{0.4em}{\vrule depth 0pt height 0.4pt width 12cm}}
\caption{The online SeqRand algorithm} \label{fig:oseqr}
\end{figure}

\begin{theorem} \label{th:o}
The cumulative loss of the online SeqRand algorithm satisfies
        \begin{multline*}
        \sum_{i=1}^n \E_{\hg_{i-1}} L(Z_i,\hg_{i-1})\\
        \qquad
                \le \und{\min}{\rho\in\M} \bigg\{ 
            	\expec{g}{\rho} \sum_{i=1}^n L(Z_i,g)
                + \expec{g}{\rho} \expecc{\hg_0^{n-1}}{\Omega_n} \sum_{i=1}^n \delta_\lam(Z_{i},g,\hg_{i-1}) 
                + \frac{K(\rho,\pi)}{\lam} \bigg\}
        \end{multline*}
In particular, when $\G$ is finite, by taking $\pi$ uniform on $\G$, we get
        $$
        \sum_{i=1}^n \E_{\hg_{i-1}} L(Z_i,\hg_{i-1})
                \le \und{\min}{g\in\G} \bigg\{ 
                \sum_{i=1}^n L(Z_i,g)
                + \expecc{\hg_0^{n-1}}{\Omega_n} \sum_{i=1}^n \delta_\lam(Z_{i},g,\hg_{i-1}) 
                + \frac{\log|\G|}{\lam} \bigg\}
        $$
\end{theorem}
 
Up to the online variance function $\delta_\lam$, the online variance inequality is the generic algorithm condition of \cite[p.11]{Hau98}.
So cases where $\delta_\lam$ are equal to zero are already known. Now new results can be obtained by using 
that for any loss function $L$ and any $\lam>0$, the online variance inequality is 
satisfied for $\delta_\lam(Z,g,g')=\frac{\lam}{2} \big[L(Z,g)-L(Z,g')\big]^2$ (proof in Section \ref{sec:corvar}).
The associated distribution $\hpi(\rho)$ is then just $\rho$.
This leads to the following corollary.

\begin{corollary} \label{cor:o}
The cumulative loss of the online SeqRand algorithm with $\delta_\lam(Z,g,g')=\frac{\lam}{2} \big[L(Z,g)-L(Z,g')\big]^2$ 
and $\hpi(\rho)=\rho$ for any $\rho\in\M$ satisfies
        \begin{multline} \label{eq:o}
        \sum_{i=1}^n \E_{\hg_{i-1}} L(Z_i,\hg_{i-1})
                \le \und{\min}{\rho\in\M} \bigg\{ 
                \expec{g}{\rho} \sum_{i=1}^n L(Z_i,g)\\
		\qquad\qquad\qquad
                + \frac{\lam}{2} \expec{g}{\rho} \expecc{\hg_0^{n-1}}{\Omega_n} \sum_{i=1}^n \big[L(Z_{i},g)-L(Z_i,\hg_{i-1})\big]^2 
                + \frac{K(\rho,\pi)}{\lam} \bigg\}.
        \end{multline}
\end{corollary}

Note that the prediction functions $\hg_i$ appears in both the left-hand side and the right-hand side of \eqref{eq:o}.
For loss functions taking their values in an interval of range $A$, we have $[L(Z_{i},g)-L(Z_i,\hg_{i-1})]^2\le A^2$.
So when $\G$ is finite, by taking $\pi$ uniform on $\G$ and
$\lam=\sqrt{(2\log|\G|)/(nA^2)}$, we obtain that the cumulative regret satisfies the more explicit cumulative regret bound:
	\bigbegarlab{eq:oo}
	\sum_{i=1}^n \E_{\hg_{i-1}} L(Z_i,\hg_{i-1}) - \und{\min}{g\in\G} \, \sum_{i=1}^n L(Z_i,g) \le A \sqrt{2n \log |\G|}.
	\bigendarlab
This bound is loose by a factor $2$ (see \cite[Theorem 2.2]{CesLug06}). 
Nevertheless an advantage of the online SeqRand algorithm with $\delta_\lam(Z,g,g')=\frac{\lam}{2} \big[L(Z,g)-L(Z,g')\big]^2$ 
is that it will take advantage of situations in which
$$\sum_{i=1}^n \expec{g}{\rho} \expecc{\hg_0^{n-1}}{\Omega_n} \big[L(Z_{i},g)-L(Z_i,\hg_{i-1})\big]^2 \ll n A^2,$$
whereas it is not clear that the exponentially weighted average forecaster does.
The proper tuning of the parameter $\lam$ is a nontrivial task, whereas 
for the exponentially weighted average forecaster with incremental updates, it has been recently proved that one can 
tune this parameter without any prior knowledge on the loss sequences \cite[Theorem 6]{Ces07}.
The argument given in the appendix p.\pageref{sec:ada} can be applied in order to use incremental updates, but for 
the online SeqRand algorithm with $\delta_\lam(Z,g,g')=\frac{\lam}{2} \big[L(Z,g)-L(Z,g')\big]^2$, we do
not know how to choose the updates in order to recover a result similar to \cite[Theorems 5 and 6]{Ces07}.

\section{Model selection aggregation under Juditsky, Rigollet and Tsybakov assumptions (\cite{Jud06})} \label{sec:jud}

The main result of \cite{Jud06}
 relies on the following assumption on the loss
function $L$ and the set $\calP$ of probability distributions on $\Z$ in which
we assume that the true distribution is.
There exist $\lam>0$ and a real-valued function $\psi$ defined on $\G\times\G$
such that for any $P\in\calP$
        \lbegarlab{eq:jud}
        \undc{\E}{Z\sim P} e^{\lam[L(Z,g')-L(Z,g)]} \le \psi(g',g) \qquad \text{ for any } g,g'\in\G\\
        \psi(g,g)=1 \qquad \text{ for any } g\in\G\\
        \text{the function } \big[g\mapsto \psi(g',g)\big] \text{ is concave for any } g'\in\G\\
        \rendarlab

Theorem \ref{th:1} gives the following result.
\begin{corollary} \label{cor:jud}
Consider the algorithm which draws uniformly its prediction function in
the set $\{\expec{g}{\pi_{-\lam \Sigma_0}} g,\dots,\expec{g}{\pi_{-\lam \Sigma_n}} g\}$.
Under Assumption \eqref{eq:jud}, its expected risk 
        $\expecc{Z_1^n}{P^{\otimes n}} \frac{1}{n+1} \sum_{i=0}^n R(\expec{g}{\pi_{-\lam \Sigma_i}}g)$ 
is upper bounded by
        \begarlab{eq:jud1}
        \und{\min}{\rho\in\M} \big\{ \expec{g}{\rho} R(g) + \frac{K(\rho,\pi)}{\lam(n+1)} \big\}.
        \endarlab
\end{corollary}

\begin{proof}
We start by proving that the variance inequality holds with
$\delta_\lam~\equiv~0$, and that we may take 
$\pirho$ as the Dirac distribution at the function
        $\expec{g}{\rho} g$.
By using Jensen's inequality and Fubini's theorem, Assumption~\eqref{eq:jud} implies that
        \begar
        \expec{g'}{\pirho} \undc{\E}{Z\sim P} \log \expec{g}{\rho} 
                e^{\lam [L(Z,g') - L(Z,g)]}\\
        \qquad\qquad\qquad\qquad\qquad
                = \; \undc{\E}{Z\sim P} \log \expec{g}{\rho} 
                e^{\lam [L(Z,\expec{g'}{\rho} g') - L(Z,g)]}\\
        \qquad\qquad\qquad\qquad\qquad
                \le \; \log \expec{g}{\rho} \undc{\E}{Z\sim P} 
                e^{\lam [L(Z,\expec{g'}{\rho} g') - L(Z,g)]}\\
        \qquad\qquad\qquad\qquad\qquad
                \le \; \log \expec{g}{\rho} \psi(\expec{g'}{\rho} g',g)\\
        \qquad\qquad\qquad\qquad\qquad
                \le \; \log \psi(\expec{g'}{\rho} g',\expec{g}{\rho} g)\\
        \qquad\qquad\qquad\qquad\qquad = \; 0,
        \endar
so that we can apply Theorem \ref{th:1}. It remains to note that in this 
context the SeqRand algorithm is the one described in the corollary.

\end{proof}

In this context, the SeqRand algorithm reduces to the randomized version of Algorithm $\aaa$. 
From Lemma \ref{le:rand}, for convex loss functions, \eqref{eq:jud1} 
also holds for the risk of Algorithm $\aaa$.
Corollary \ref{cor:jud} also shows that the risk bounds for Algorithm $\aaa$ 
proved in \cite[Theorem 3.2 and the examples of Section 4.2]{Jud06}
hold with the same constants for the SeqRand algorithm
(provided that the expected risk $\wrt$ the training set distribution 
is replaced by the expected risk $\wrt$ both training set and randomizing distributions).

On Assumption \eqref{eq:jud} we should say that it does not a priori require the function
$L$ to be convex. Nevertheless, any known relevant examples deal with ``strongly'' convex loss functions
and we know that in general the assumption will not hold for the SVM (or hinge) loss function and for 
the absolute loss function.
Indeed, without further assumption, one cannot expect rates better than $1/\sqrt{n}$ for these loss functions 
(see Section \ref{sec:precisex}).

By taking the appropriate variance function $\delta_{\lam}(Z,g,g')$, it 
is possible to prove that the results in \cite[Theorem 3.1 and Section 4.1]{Jud06}
holds for the SeqRand algorithm (provided that the expected risk $\wrt$ the training set distribution 
is replaced by the expected risk $\wrt$ both training set and randomizing distributions). 
The choice of $\delta_{\lam}(Z,g,g')$, which for sake of shortness we do not specify,
is in fact such that the resulting SeqRand algorithm is again the randomized version of Algorithm $\aaa$.


\section{Standard-style statistical bounds} \label{sec:var}

This section proposes new results of a different kind. In the previous sections,
under convexity assumptions, we were able to achieve fast rates.
Here we have assumption neither on the loss function nor on the 
probability generating the data. Nevertheless we show that the SeqRand 
algorithm applied for $\delta_\lam(Z,g,g')= \lam [L(Z,g)-L(Z,g')]^2 /2$
satisfies a sharp standard-style statistical bound.

This section contains two parts: the first one provides results in expectation (as in the preceding sections)
whereas the second part provides deviation inequalities on the risk that requires advances
on the sequential prediction analysis.

\subsection{Bounds on the expected risk}

\subsubsection{Bernstein's type bound} \label{sec:berns}

\begin{theorem} \label{cor:var}
Let $V(g,g') = \expecc{Z}{P} \big\{[L(Z,g)-L(Z,g')]^2\big\}.$ Consider the SeqRand algorithm
(see p.\pageref{fig:seqr})
applied with $\delta_\lam(Z,g,g')= \lam [L(Z,g)-L(Z,g')]^2 /2$ and $\pirho=\rho$.
Its expected risk $\Ezun \expec{g}{\hmu} R(g)$, where we recall that
$\hmu$ denotes the randomizing distribution, satisfies
        \begarlab{eq:corvar}
        \expecc{Z_1^n}{P^{\otimes n}} \expec{g'}{\hmu} R(g') 
                \le \und{\min}{\rho\in\M} \Big\{ 
                \expec{g}{\rho} R(g) 
                + \frac{\lam}{2} \expec{g}{\rho} \expecc{Z_1^n}{P^{\otimes n}} \expec{g'}{\hmu} V(g,g')
                + \frac{K(\rho,\pi)}{\lam(n+1)} \Big\}\\
        \endarlab
\end{theorem}

\begin{proof}
See Section \ref{sec:corvar}.
\end{proof}

To make \eqref{eq:corvar} more explicit and to obtain a generalization error bound in which the randomizing distribution 
does not appear in the $\rhs$ of the bound, the following corollary considers a widely used assumption
relating the variance term to the excess risk (see Mammen and Tsybakov \cite{Mam99,Tsy04}, and also Polonik \cite{Pol95}). 
Precisely, from Theorem \ref{cor:var}, we obtain 

\begin{corollary} \label{cor:mamtsy}
If there exist $0\le \gamma\le 1$ 
and a prediction function $\tildg$ (not necessarily in $\G$) such that 
$V(g,\tildg) \le c [R(g)-R(\tildg)]^{\gamma}$ for any $g\in\G$, 
the expected risk $\calE=\Ezun \expec{g}{\hmu} R(g)$ of the SeqRand algorithm used in Theorem \ref{cor:var} satisfies
\begitem
\item When $\gamma=1$, 
        \begar
        \calE - R(\tildg) \le \und{\min}{\rho\in\M} \Big\{ 
                \frac{1+c\lam}{1-c\lam} \big[ \expec{g}{\rho} R(g)-R(\tildg) \big]
                + \frac{K(\rho,\pi)}{(1-c\lam)\lam(n+1)} \Big\}
        \endar
        In particular, for $\G$ finite, $\pi$ the uniform distribution, $\lam=1/(2c)$, 
        when $\tildg$ belongs to $\G$, we get
        $\calE \le \und{\min}{g\in\G} R(g) + \frac{4c\log|\G|}{n+1}.$

\item When $\gamma<1$, for any $0<\beta<1$ and for $\tildR(g)\eqdef R(g)-R(\tildg)$, 
        \begar
        \calE - R(\tildg) \le \Big\{ \frac{1}{\beta}\,\und{\min}{\rho\in\M} \Big( \expec{g}{\rho} [ \tildR(g) 
                + c \lam \tildR^\gamma(g) ] 
                + \frac{K(\rho,\pi)}{\lam(n+1)} \Big) \Big\} \vee \big( \frac{c\lam}{1-\beta} \big)^{\frac{1}{1-\gamma}}.
        \endar
\enditem
\end{corollary}

\begin{proof}
See Section \ref{sec:cormamtsy}.
\end{proof}
To understand the sharpness of Theorem \ref{cor:var}, we have to compare this result with the following one that comes from 
the traditional ({PAC}-Bayesian) statistical learning approach which relies on supremum of empirical processes.
In the following theorem, we consider the estimator minimizing the uniform bound, i.e. the estimator for which we have the smallest
upper bound on its generalization error. 

\begin{theorem} \label{cor:pac}
We still use $V(g,g') = \E_Z \big\{[L(Z,g)-L(Z,g')]^2\big\}.$
The generalization error of
the algorithm which draws its prediction function according to the Gibbs distribution $\pi_{-\lam \Sigma_n}$
satisfies
        \begarlab{eq:pacexp2}
        \Ezun \expec{g'}{\pi_{-\lam \Sigma_n}} R(g')\\
        \quad
                \le \und{\min}{\rho\in\M} \Big\{ \expec{g}{\rho} R(g) + \frac{K(\rho,\pi)+1}{\lam n} 
                +\lam \expec{g}{\rho} \Ezun \expec{g'}{\pi_{-\lam \Sigma_n}} V(g,g')\\
        \quad
                +\lam \frac{1}{n} \sum_{i=1}^n \expec{g}{\rho} \Ezun \expec{g'}{\pi_{-\lam \Sigma_n}} 
                [L(Z_i,g)-L(Z_i,g')]^2
                \Big\}.
        \endarlab
Let $\varphi$ be the positive convex increasing function defined as
$\varphi(t) \eqdef \frac{e^t-1-t}{t^2}$ and $\varphi(0) = \frac{1}{2}$ by continuity.
When $\sup_{z\in\Z,g\in\G,g'\in\G} |L(z,g')-L(z,g)| \le B$, we also have
        \begarlab{eq:pacexp1}
        \Ezun \expec{g'}{\pi_{-\lam \Sigma_n}} R(g') \le \und{\min}{\rho\in\M}
                \Big\{ \expec{g}{\rho} R(g) \\
        \qquad\quad
                +\lam \varphi(\lam B) \expec{g}{\rho} \Ezun \expec{g'}{\pi_{-\lam \Sigma_n}} V(g,g')
                +\frac{K(\rho,\pi)+1}{\lam n}\Big\}.
        \endarlab

\end{theorem}
\begin{proof}
See Section \ref{sec:corpac}.
\end{proof}


As in Theorem \ref{cor:var}, there is a variance term in which the randomizing distribution is involved. 
As in Corollary \ref{cor:mamtsy}, one can convert \eqref{eq:pacexp1} into a proper generalization error bound,
that is a non trivial bound $\Ezun \E_{g\sim\pi_{-\lam \Sigma_n}} R(g) \le \B(n,\pi,\lam)$ where 
the training data do not appear in $\B(n,\pi,\lam)$.

By comparing \eqref{eq:pacexp1} and \eqref{eq:corvar},
we see that the classical approach requires 
the quantity ${\sup}_{g\in\G,g'\in\G} | L(Z,g')-L(Z,g) |$
to be uniformly bounded and the unpleasing function $\varphi$ appears. 
In fact, using technical small expectations theorems (see e.g. \cite[Lemma 7.1]{Aud03a}), 
exponential moments conditions on the above quantity would be sufficient. 

The symmetrization trick used to prove Theorem \ref{cor:var} is performed in the prediction functions space. 
We do not call on the second virtual training set currently used in
statistical learning theory (see \cite{Vap95}). 
Nevertheless both symmetrization tricks end up to the same nice property: we need 
no boundedness assumption on the loss functions. In our setting, symmetrization on training data leads to 
an unwanted expectation and to a constant four times larger (see the two variance terms of \eqref{eq:pacexp2}
and the discussion in \cite[Section 8.3.3]{Aud03b}).

In particular, deducing from Theorem \ref{cor:pac} a corollary similar to Corollary \ref{cor:mamtsy} is only possible 
through \eqref{eq:pacexp1} and provided that we have a boundedness assumption on
$\sup_{z\in\Z,g\in\G,g'\in\G} |L(z,g')-L(z,g)|$.
Indeed one cannot use \eqref{eq:pacexp2} because of the last variance term in \eqref{eq:pacexp2} 
(since $\Sigma_n$ depends on $Z_i$).

Our approach has nevertheless the following limit: 
the proof of Corollary \ref{cor:mamtsy} does not use a chaining argument. 
As a consequence, in the particular case when the model has polynomial entropies 
(see e.g. \cite{Mam99}) and when the assumption in Corollary \ref{cor:mamtsy} holds for $\gamma<1$ (and not for $\gamma=1$),
Corollary \ref{cor:mamtsy} does not give the minimax optimal convergence rate.
Combining the better variance control presented here with the chaining argument is an open problem.

\subsubsection{Hoeffding's type bound} 

Contrary to generalization error bounds coming from Bernstein's inequality, \eqref{eq:corvar} does not require
any boundedness assumption. For bounded losses, without any variance assumption 
(i.e. roughly when the assumption used in Corollary \eqref{cor:mamtsy} does not hold for $\gamma>0$), 
tighter results are obtained by using Hoeffding's inequality, that is:
for any random variable $W$ satisfying $a\le W \le b$, then for any $\lam>0$
        \begar
        \E e^{\lam(W-\E W)} \le e^{\lam^2(b-a)^2/8}.
        \endar

\begin{theorem} \label{th:hoeff}
Assume that for any $z\in\Z$ and $g\in\G$, we have $a \le L(z,g) \le b$ for some reals $a,b$.
Consider the SeqRand algorithm (see \pageref{fig:seqr})
applied with $\delta_\lam(Z,g,g')= \lam (b-a)^2 /8$ and $\pirho=\rho$.
Its expected risk $\Ezun \expec{g}{\hmu} R(g)$, where we recall that
$\hmu$ denotes the randomizing distribution, satisfies
        \begarlab{eq:corhoeff1}
        \Ezun \expec{g}{\hmu} R(g) 
                \le \und{\min}{\rho\in\M} \Big\{ 
                \expec{g}{\rho} R(g) 
                + \frac{\lam (b-a)^2}{8} 
                + \frac{K(\rho,\pi)}{\lam(n+1)} \Big\}
        \endarlab
In particular, when $\G$ is finite, by taking $\pi$ uniform on $\G$
and $\lam=\sqrt{\frac{8\log|\G|}{(b-a)^2(n+1)}}$, we get
        \begarlab{eq:corhoeff2}
        \Ezun \expec{g}{\hmu} R(g) 
                - \und{\min}{g\in\G} R(g) \le (b-a) \sqrt{\frac{\log|\G|}{2(n+1)}}
        \endarlab
\end{theorem}

\begin{proof}
From Hoeffding's inequality, we have
        \begar
        \expec{g'}{\pirho} \log \expec{g}{\rho} e^{\lam [L(Z,g') - L(Z,g)]}
                & = & \log \expec{g}{\rho} e^{\lam [\expec{g'}{\pirho} L(Z,g') - L(Z,g)]}\\
        & \le & \frac{\lam^2 (b-a)^2}{8},
        \endar
hence the variance inequality holds for $\delta_\lam\equiv \lam (b-a)^2 /8$ and $\pirho=\rho$.
The result directly follows from Theorem \ref{th:1}.
\end{proof}

The standard point of view (see Appendix \ref{app:stdhoeff}) applies Hoeffding's inequality to the random variable
$W=L(Z,g')-L(Z,g)$ for $g$ and $g'$ fixed and $Z$ drawn according to the probability generating the data.
The previous theorem uses it on the random variable
$W=L(Z,g')-\expec{g}{\rho} L(Z,g)$ for fixed $Z$ and fixed probability distribution $\rho$ but
for $g'$ drawn according to $\rho$.
Here the gain is a multiplicative factor equal to $2$ (see Appendix \ref{app:stdhoeff}).

\subsection{Deviation inequalities}

For the comparison between Theorem \ref{cor:var} and Theorem \ref{cor:pac} to be fair, 
one should add that \eqref{eq:pacexp1} and \eqref{eq:pacexp2} come from
deviation inequalities that are not exactly obtainable to the author's knowledge with the arguments developed here.
Precisely, consider the following adaptation of Lemma 5 of \cite{Zha05}. 
\begin{lemma} \label{lem:conczha}
Let $\A$ be a learning algorithm which produces the prediction function $\A(Z_1^i)$ at time $i+1$, i.e. from the data 
$Z_1^i=(Z_1,\dots,Z_i)$.
Let $\L$ be the randomized algorithm which produces a prediction function $\L(Z_1^n)$ 
drawn according to the uniform
distribution on $\{\A(\emptyset),\A(Z_1),\dots,\A(Z_1^n)\}$. 
Assume that $\sup_{z,g,g'} | L(z,g) - L(z,g')| \le B$ for some $B>0$.
Conditionally to $Z_1,\dots,Z_{n+1}$, the expectation of the risk of $\L$ $\wrt$ to the uniform draw is
        $
        \frac{1}{n+1} \sum_{i=0}^n R[\A(Z_1^i)] 
        $
and satisfies: for any $\eta>0$ and $\eps>0$, for any reference prediction function $\tildg$,
with probability at least $1-\eps$ $\wrt$ the distribution of $Z_1,\dots,Z_{n+1}$, 
        \begarlab{eq:devrand2}
        \frac{1}{n+1} \sum_{i=0}^n R[\A(Z_1^i)] -R(\tildg)
                \le \frac{1}{n+1} \sum_{i=0}^n \big\{ L[Z_{i+1},\A(Z_1^i)] - L(Z_{i+1},\tildg) \big\}\\
        \qquad\qquad\qquad\qquad\qquad
                + \eta \varphi(\eta B) \frac{1}{n+1} \sum_{i=0}^n  
                V[\A(Z_1^i),\tildg] + \frac{\logeps}{\eta(n+1)}
        \endarlab
where we still use $V(g,g') = \E_{Z} \big\{[L(Z,g)-L(Z,g')]^2\big\}$ for any prediction functions $g$ and $g'$
and $\varphi(t) \eqdef \frac{e^t-1-t}{t^2}$ for any $t>0$.
\end{lemma}

\begin{proof}
See Section \ref{sec:proofconczha}. 
\end{proof}

We see that two variance terms appear.
The first one comes from the worst-case analysis and is hidden in
$\sum_{i=0}^n \big\{ L[Z_{i+1},\A(Z_1^i)] - L(Z_{i+1},\tildg) \big\}$
and the second one comes from the concentration result (Lemma \ref{lem:zha}).
The presence of this last variance term annihilates the benefits of our approach in 
which we were manipulating variance terms much smaller than the traditional Bernstein's variance term.

To illustrate this point, consider for instance least square regression with boun\-ded outputs: 
from Theorem \ref{th:hau} and Table \ref{tab:table}, 
the hidden variance term is null.
In some situations, the second variance term $\frac{1}{n+1} \sum_{i=0}^n V[\A(Z_1^i),\tildg]$ 
may behave like a positive constant: for instance, this occurs when $\G$ contains two very different functions having 
the optimal risk $\min_{g\in\G} R(g)$. By optimizing $\eta$, this will lead to a 
deviation inequality of order $n^{-1/2}$
even though from \eqref{eq:lsup} the procedure has $n^{-1}$-convergence rate in expectation.
In \cite[Theorem~3]{Aud07b}, in a rather general learning setting, this deviation inequality of order $n^{-1/2}$ is proved to be optimal.

To conclude, for deviation inequalities, we cannot expect to do better
than the standard-style approach since at some point we use a Bernstein's type bound
$\wrt$ the distribution generating the data. Besides procedures based on
worst-case analysis seem to suffer higher fluctuations of the risk than necessary
(see \cite[discussion of Theorem 3]{Aud07b}).

\begin{remark}
Lemma \ref{lem:conczha} should be compared with Lemma \ref{lem:rand}. 
The latter deals with results in expectation while the former
concerns deviation inequalities. Note that Lemma \ref{lem:conczha}
requires the loss function to be bounded and makes a variance term appear.
\end{remark}

\section{Application to $L_q$-regression for unbounded outputs} \label{sec:pow} 

In this section, we consider the $L_q$-loss: $L(Z,g)=|Y-g(X)|^q.$ 
As a warm-up exercise, we tackle the absolute loss setting (i.e. $q=1$). 
The following corollary holds without any assumption on the output 
(except naturally that if $E_Z |Y|<+\infty$ to ensure finite risk).

\begin{corollary} \label{cor:abs}
Let $q=1$. Assume that ${\sup}_{g\in\G} \, E_Z \, g(X)^2 \le b^2$ for some $b>0$. 
There exists an estimator $\hg$ such that
        \begarlab{eq:uq1}
        \E R(\hg) - \und{\min}{g\in\G} R(g) \le 
                 2b\sqrt{\frac{2\log |\G|}{n+1}}.
        \endarlab
\end{corollary}

\begin{proof}
Using
        $
        \expecc{Z}{P} \big\{[|Y-g(X)|-|Y-g'(X)|]^2\big\} 
                \le 4b^2
        $
and Theorem~\ref{cor:var}, the algorithm considered in Theorem~\ref{cor:var} satisfies
        $
        \E R(\hg) - \undc{\min}{g\in\G} R(g) \le 2\lam b^2 +\frac{\log |\G|}{\lam(n+1)},
        $
which gives the desired result by taking $\lam = \sqrt{\frac{\log |\G|}{2b^2(n+1)}}$.
\end{proof}


Now we deal with the strongly convex loss functions (i.e. $q>1$).
By using Theorem~\ref{th:1} jointly with the symmetrization idea developed in the previous section
allows to obtain new convergence rates in heavy noise situation, i.e. when the output is not constrained to 
have a bounded exponential moment. We start with the following theorem concerning general loss functions.

\begin{theorem} \label{th:genpow}
Assume that ${\sup}_{g\in\G,x\in\X} |g(x)| \le b$ for some $b>0$, and that the output space is $\Y=\R$.
Let $B\ge b$. Consider a loss function $L$ which can be written as
        $
        L[(x,y),g] = \ell[y,g(x)], 
        $
where the function $\ell: \R \times \R \rightarrow \R$ satisfies: there exists $\lam_0>0$ such that for any $y\in[-B;B]$,
the function $y' \mapsto e^{-\lam_0 \ell(y,y')}$ is concave on $[-b;b]$. 
Let 
        \begar
        \Delta(y)~=~\und{\sup}{|\alpha|\le b,|\beta|\le b} \big[\ell(y,\alpha) - \ell(y,\beta)\big].
        \endar
For $\lam\in\,(0; \lam_0]$, consider the algorithm that draws uniformly its prediction function in
the set $\{\expec{g}{\pi_{-\lam \Sigma_0}} g,\dots,\expec{g}{\pi_{-\lam \Sigma_n}} g\}$, and consider
the deterministic version of this randomized algorithm. The expected risk of these algorithms satisfy
        \begarlab{eq:genpow}
        \Ezun R(\frac{1}{n+1} \sum_{i=0}^n \expec{g}{\pi_{-\lam \Sigma_i}}g) \\
        \qquad\quad
                \le \Ezun \frac{1}{n+1} \sum_{i=0}^n R(\expec{g}{\pi_{-\lam \Sigma_i}}g)\\
        \qquad\quad
                \le \und{\min}{\rho\in\M} \Big\{ \expec{g}{\rho} R(g) + \frac{K(\rho,\pi)}{\lam(n+1)} \Big\}\\
        \qquad\quad\quad
                + \E\Big\{ \frac{\lam \Delta^2(Y)}{2}  
                \ds1_{\lam \Delta(Y)<1;|Y|> B} + \big[\Delta(Y)-\frac{1}{2\lam}\big] \ds1_{\lam \Delta(Y)\ge1;|Y|> B} \Big\}.
        \endarlab
\end{theorem}

\begin{proof}
See Section \ref{sec:proofgenpow}.
\end{proof}

\begin{remark}
For $y\in[-B;B]$, concavity of $y' \mapsto e^{-\lam_0 \ell(y,y')}$ on $[-b;b]$ for $\lam_0>0$
implies convexity of $y' \mapsto \ell(y,y')$ on $[-b;b]$.
\end{remark}

In particular, for least square regression, Theorem \ref{th:genpow} can be simplified into:

\begin{theorem} \label{th:gen2}
Assume that ${\sup}_{g\in\G,x\in\X} |g(x)| \le b$ for some $b>0$.
For any $0 < \lam \le 1/(8b^2)$, the expected risk of the algorithm that draws uniformly its prediction function among
$\expec{g}{\pi_{-\lam \Sigma_0}} g$,$\dots$,$\expec{g}{\pi_{-\lam \Sigma_n}} g$ is upper bounded by
        \begarlab{eq:gen2}
        \und{\min}{\rho\in\M} \Big\{ \expec{g}{\rho} R(g) 
                + 8 \lam b^2 \E\big( Y^2  \ds1_{|Y|>(8\lam)^{-1/2}} \big)
                + \frac{K(\rho,\pi)}{\lam(n+1)} \Big\}.
        \endarlab
\end{theorem}

\begin{proof}
For any $B\ge b$, for any $y\in[-B;B]$,
straightforward computations show that $y' \mapsto e^{-\lam_0 (y-y')^2}$ is concave on $[-b;b]$
for $\lam_0=\frac{1}{2(B+b)^2}$, so that we can apply Theorem \ref{th:genpow}.
We have $\Delta(y) = 4 b |y|$ for any $|y|\ge b$ so that by optimizing the parameter $B$, we obtain that the 
expected risk of the algorithm is upper bounded by
        \begar
        & \und{\min}{\rho\in\M} \Big\{ \expec{g}{\rho} R(g) + \frac{K(\rho,\pi)}{\lam(n+1)} \Big\}
                + \E\big\{ \big(4b|Y|-\frac{1}{2\lam}\big) \ds1_{|Y|\ge (4b\lam)^{-1}} \big\}\\
        & \qquad\qquad\qquad\qquad\qquad
                + \E\big\{ 8 \lam b^2 Y^2  \ds1_{(2\lam)^{-1/2} - b < |Y|< (4b\lam)^{-1}} \big\}\\
        \le & \und{\min}{\rho\in\M} \Big\{ \expec{g}{\rho} R(g) + \frac{K(\rho,\pi)}{\lam(n+1)} \Big\}
                + \E\big\{ 8 \lam b^2 Y^2 \ds1_{|Y|\ge (4b\lam)^{-1}} \big\}\\
        & \qquad\qquad\qquad\qquad\qquad
                + \E\big\{ 8 \lam b^2 Y^2  \ds1_{(2\lam)^{-1/2} - b < |Y|< (4b\lam)^{-1}} \big\},
        \endar
which gives the desired result.
\end{proof}

Theorem \ref{th:gen2} improves \cite[Theorem 1]{Bun05}.

\begin{corollary} \label{cor:gen2}
Under the assumptions
        \lbegar
        {\sup}_{g\in\G,x\in\X} |g(x)| \le b \qquad \text{for some } b>0\\
        \E |Y|^s\le A \qquad \text{for some } s\ge 2 \text{ and } A>0\\
        \G \text{ finite}
        \rendar
for $\lam=C_1 \big(\frac{\log |\G|}{n}\big)^{2/(s+2)}$ where $C_1>0$ and $\pi$ the uniform
distribution on $\G$, the expected risk of the  
algorithm that draws uniformly its prediction function among
$\expec{g}{\pi_{-\lam \Sigma_0}} g$,$\dots$,$\expec{g}{\pi_{-\lam \Sigma_n}} g$ is upper bounded by
        \begarlab{eq:squ}
         \und{\min}{g\in\G} R(g) + C \big(\frac{\log |\G|}{n}\big)^{s/(s+2)}
        \endarlab
for a quantity $C$ which depends only on $C_1$, $b$, $A$ and $s$.
\end{corollary}

Juditsky, Rigollet and Tsybakov proved that Corollary \ref{cor:gen2} 
can also be obtained through a simple adaptation of their original analysis
(see \cite[Section 4.1]{Jud06}).

\begin{proof}
The moment assumption on $Y$ implies 
        \begarlab{eq:moments}
        \alpha^{s-q} \E |Y|^q\ds1_{|Y|\ge \alpha} \le A$ for any $0 \le q \le s$ and $\alpha\ge 0.
        \endarlab
As a consequence, the second term in \eqref{eq:gen2} is bounded by
        $8\lam b^2 A (2\lam)^{(s-2)/2}$,
so that \eqref{eq:gen2} is upper bounded by
        $
        {\min}_{g\in\G} R(g) + A 2^{2+s/2} b^2 \lam^{s/2} + \frac{\log |\G|}{\lam n},
        $
which gives the desired result.
\end{proof}

In particular, with the minimal assumption $\E Y^2\le A$ (i.e. $s=2$), the convergence 
rate is of order $n^{-1/2}$, and at the opposite, when $s$ goes to infinity, we recover the $n^{-1}$ rate
we have under exponential moment condition on the output.

Using Theorem \ref{th:genpow}, we can generalize Corollary \ref{cor:gen2} to $L_q$-regression
and obtain the following result.

\begin{corollary} \label{cor:genp}
Let $q>1$. Assume that
        \lbegar
        {\sup}_{g\in\G,x\in\X} |g(x)| \le b \qquad \text{for some } b>0\\
        \E |Y|^s\le A \qquad \text{for some } s\ge q \text{ and } A>0\\
        \G \text{ finite}
        \rendar
Let $\pi$ be the uniform distribution on $\G$, $C_1>0$ and 
        \begar
        \lam = \left\{ \begin{array}{lll}
        C_1 \big(\frac{\log |\G|}{n}\big)^{(q-1)/s} & \qquad \text{ when }q \le s < 2q-2\\
        C_1 \big(\frac{\log |\G|}{n}\big)^{q/(s+2)} & \qquad \text{ when } s \ge 2q-2\\
        \end{array} \right..
        \endar
The expected risk of the algorithm which draws uniformly its prediction function among
$\expec{g}{\pi_{-\lam \Sigma_0}} g$, $\dots$, $\expec{g}{\pi_{-\lam \Sigma_n}} g$ is upper bounded by
        \begar
        \left\{ \begin{array}{lll}
        \und{\min}{g\in\G} R(g) + C \big(\frac{\log |\G|}{n}\big)^{1-\frac{q-1}{s}} & \qquad \text{ when }q \le s \le 2q-2\\
        \und{\min}{g\in\G} R(g) + C \big(\frac{\log |\G|}{n}\big)^{1-\frac{q}{s+2}} & \qquad \text{ when }s\ge 2q-2\\
        \end{array} \right..
        \endar
for a quantity $C$ which depends only on $C_1$, $b$, $A$, $q$ and $s$.
\end{corollary}

\begin{proof}
See Section \ref{sec:corgenp}.
\end{proof}

\begin{remark}
For $q>2$, low convergence rates (that is $n^{-\gamma}$ with $\gamma<1/2$) appear when the moment assumption
is weak: $\E|Y|^s\le A$ for some $A>0$ and $q\le s < 2q-2$. 
Convergence rates faster that the standard non parametric rates $n^{-1/2}$ are achieved for 
$s> 2q-2$. Fast convergence rates systematically occurs when $1<q<2$ since for these values of $q$, we have
$s\ge q > 2q-2$. Surprisingly, for $q=1$, the picture is completely different (see Section \ref{sec:precisex} 
for discussion and minimax optimality of the results of this section).
\end{remark}

\begin{remark}
Corollary \ref{cor:genp} assumes that the prediction functions in $\G$ are uniformly bounded. It is an open problem to have
the same kind of results under weaker assumptions such as a finite moment condition similar to the one used in 
Corollary \ref{cor:abs}.
\end{remark}
\fi

\ifRR
\section{Lower bounds} \label{sec:assouad}

The simplest way to assess the quality of an algorithm and of its expected risk upper bound is 
to prove a risk lower bound saying that no algorithm has better convergence rate.
This section provides this kind of assertions.
\fi
The lower bounds developed here have the same spirit as the ones in \cite{Bre79,Ass83,Bir83},
\cite[Chap. 15]{Dev00} and \cite[Section 5]{Aud03d} to the extent that
it relies on the following ideas:
\begitem
\item the supremum of a quantity $\calQ(P)$ when the distribution $P$ belongs to some set $\calP$ is larger than
the supremum over a well chosen finite subset of $\calP$, and consequently is larger 
than the mean of $\calQ(P)$ when the distribution $P$ is drawn uniformly in the finite subset.
\item when the chosen subset is a hypercube of $2^m$ distributions (see Definition \ref{def:hyper}),
the design of a lower bound over the $2^m$ distributions reduces to the design of a lower
bound over two distributions.
\item when a data sequence $Z_1,\dots,Z_n$ has similar likelihoods according to two different probability distributions, then
no estimator will be accurate for both distributions: the maximum over
the two distribution of the risk of any estimator trained on this sequence will be 
all the larger as the Bayes-optimal prediction associated with the two distributions
are `far away'.
\enditem
We refer the reader to \cite{Bir05} and \cite[Chap.~2]{Tsy04c} for
lower bounds not parlicularly based on finding the appropriate hypercube.
Our analysis focuses on hypercubes since in several settings they afford to obtain 
\ifRR lower \fi bounds 
with both the right convergence rate and close to optimal constants. 
Our contribution \ifRR in this section \fi is 
\begitem
\item 
\ifRR to provide results for general non-regularized loss functions
(we recall that non-regularized loss functions are loss functions which can be written as 
$L[(x,y),g]=\ell[y,g(x)]$ for some function $\ell:\Y\times\Y\rightarrow\R$),
\fi
\item to improve the upper bound on the variational distance appearing in Assouad's argument,
\item to generalize the argument to asymmetrical hypercubes 
\ifRR which, to our knowledge, is the only way to find the lower bound matching the
upper bound of Corollary \ref{cor:genp} for $q\le~s\le~2q-2$,
\fi
\item to express the lower bounds in terms of similarity measures between two distributions characterizing the hypercube.
\ifRR
\item to obtain lower bounds matching the upper bounds obtained in the previous sections.
\fi
\enditem

\begin{remark} \label{rem:onlow}
In \cite{Hau98}, the optimality of the constant in front of the $(\log |\G|)/n$ has been proved
by considering the situation when both $|\G|$ and $n$ goes to infinity. 
Note that this worst-case analysis constant is not necessary the same as our batch setting constant.
This section shows that the batch setting constant is not ``far'' from the
worst-case analysis constant.

Besides Lemma \ref{lem:rand}, which can be used to convert any worst-case analysis upper bounds 
into a risk upper bound in our batch setting, also means that
any lower bounds for our batch setting leads to a lower bound in the sequential 
prediction setting (the converse is not true). Indeed the cumulative loss on the worst sequence of data
is bigger than the average cumulative loss when the data are taken i.i.d. from
some probability distribution.
As a consequence, the bounds developed in this section partially solve the open problem
introduced in \cite[Section 3.4]{Hau98} consisting in developing tight non-asymptotical lower bounds.
For least square loss and entropy loss, our bounds are off by a multiplicative factor smaller
than $4$ (see Remarks \refp{rem:common} and \refp{rem:commop}). 
\end{remark}



\ifRR
This section is organized as follows.
Section \ref{sec:hyper} defines the quantities that characterize hypercubes
of probability distributions and details the links between them. 
Section \ref{sec:sim} defines a similarity measure between probability distributions
coming from $f$-divergences (see \cite{Csi67}) and gives their main properties.
We give our main lower bounds  
in Section \ref{sec:genassouad}.
These bounds are illustrated in Section \ref{sec:examples}.
\fi

\ifRR

\subsectionjy{Hypercube of probability distributions} \label{sec:hyper}

Let $m\in \N^*$.
\ifRR Consider a family of $2^m$ probability distributions on $\Z$
        \begar
        \big\{ P_{\sigmav} : \sigmav \triangleq (\sigma_1,\dots,\sigma_m) \in \{-;+\}^m \big\}
        \endar
\else
Consider a family 
        $
        \big\{ P_{\sigmav} : \sigmav \triangleq (\sigma_1,\dots,\sigma_m) \in \{-;+\}^m \big\}
        $
of $2^m$ probability distributions on $\Z$
\fi
having the same first marginal, denoted~$\mu$:
        $$P_{\sigmav}(dX) = P_{(+,\dots,+)}(dX) \triangleq \mu(dX) \text{ for any } 
                        \sigmav \in \{-;+\}^m,$$
and such that there exist
\begitem
\item a partition $\X_0,\dots,\X_m$ of $\X$,
\item functions $\ha$ and $\hb$ defined on $\X-\X_0$ taking their values in $\Y$
\item functions $\pa$ and $\pb$ defined on $\X-\X_0$ taking their values in $[0;1]$
\enditem
for which for any $j \in \{1,\dots,m\},$ for any $x \in \X_j$, we have
                \begarlab{eq:hyp}
                P_{\sigmav}\big( Y = \ha(x) \big| X=x \big) = p_{\sigma_j}(x)
                        = 1 - P_{\sigmav}\big( Y = \hb(x) \big| X=x \big),
                \endarlab
and for any $x \in \X_0$, the distribution of $Y$ knowing $X=x$
is independent of $\sigmav$ (i.e. the $2^m$ conditional distributions are identical).

In particular, \eqref{eq:hyp} means that for any $x\in\X-\X_0$, the conditional 
probability of the output knowing the input $x$ is concentrated on two 
values\ifRR\, and that, under the distribution $P_{\sigmav}$, the disproportion between
the probabilities of these two values is all the larger as $p_{\sigma_j}(x)$ 
is far from $1/2$ for $j$ the integer such that $x\in\X_j$\fi.

\ifRR
\begin{remark} \label{rem:haneqhb}
Equality \eqref{eq:hyp} indirectly implies that for any $x$,
$\ha(x)\neq\hb(x)$. This is not at all restricting since points for which
we would have liked $\ha(x)=\hb(x)$ can be put in the ``garbage'' set $\X_0$.
\end{remark}
\fi

\ifRR
\begin{figure}[ht] 
\begin{center}
\includegraphics[width=10cm]{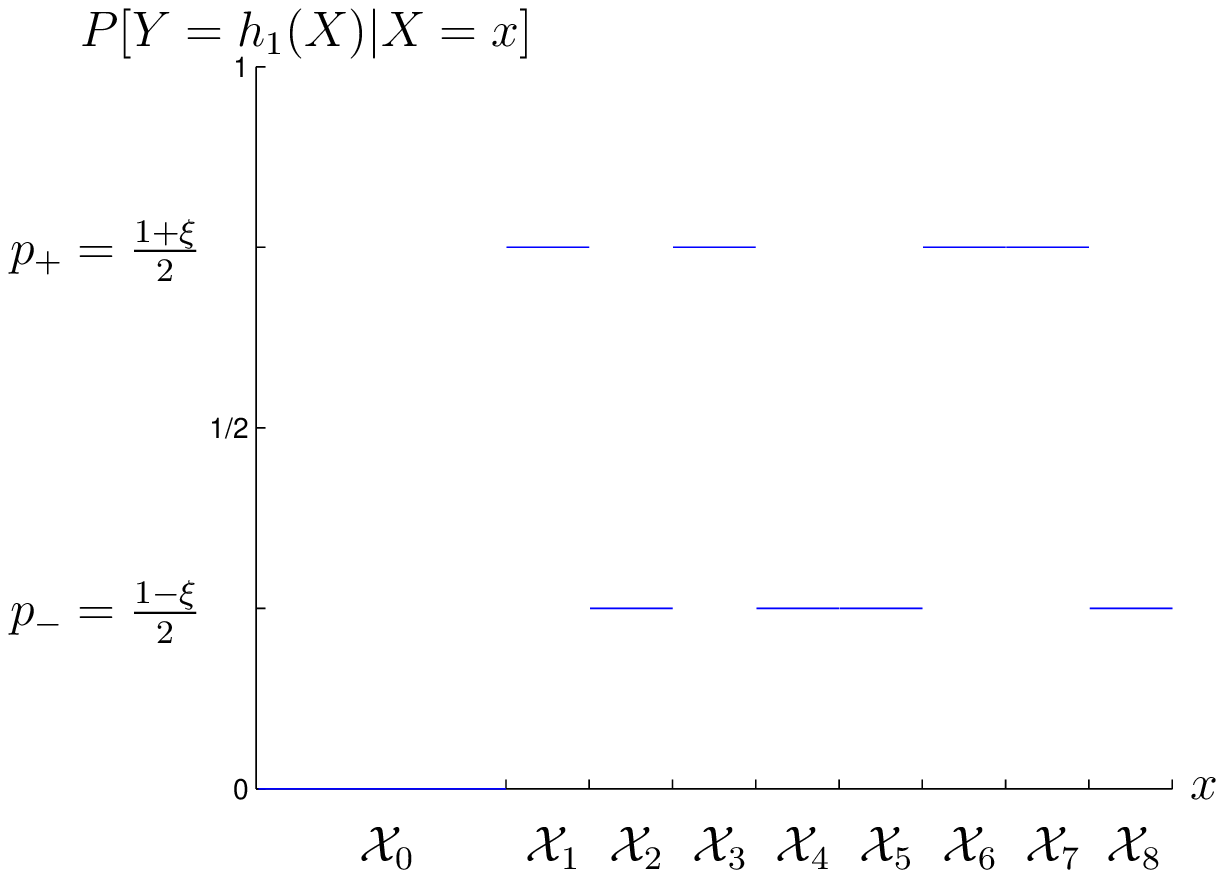}
\vspace{-0.2cm}\caption{Representation of a probability distribution of the hypercube. 
Here the hypercube is a constant and symmetrical one 
(see Definition \ref{def:hyp2}) with $m=8$ and the probability
distribution is characterized by $\bar{\sigma}=(+,-,+,-,-,+,+,-)$.
} \label{fig:hyp}
\end{center}
\end{figure}
\fi

\begin{definition} \label{def:hyper}
The family of $2^m$ probability distributions will be referred to as an 
\emph{hypercube} of distributions
if and only if for any $j\in\{1,\dots,m\},$
        \begin{itemize}
        \item 
        the probability $\mu(\X_j)=\mu(X\in\X_j)$ is independent of~$j$, 
        i.e. $\mu(\X_1)= \dots = \mu(\X_m)$,
        \item the law of $\big(p_+(X),p_-(X),\ha(X),\hb(X)\big)$ 
        when $X$ is drawn according to the conditional distribution 
        $\mu(\bullet|\X_j)\eqdef\mu(\cdot|X\in\X_j)$ is independent of $j$, i.e. the $m$ conditional 
        distributions are identical.
        \end{itemize}
\end{definition}

\ifRR
\begin{remark}
The typical situation in which we encounter hypercubes are when $\X \subseteq \R^d$ 
for some $d\ge1$ and when we have translation invariance to the extent that 
there exist $t_2,\dots,t_d$ in $\R^d$ such that for any $j\in\{2,\dots,m\}$, 
$\X_j=\X_1+t_j$ and for any $x\in\X_1$,
        \begar
        \big(p_+(x+t_j),p_-(x+t_j),\ha(x+t_j),\hb(x+t_j)\big) \\
        \qquad\qquad\qquad\qquad\qquad\qquad\qquad
                = \big(p_+(x),p_-(x),\ha(x),\hb(x)\big).
        \endar
A special hypercube is illustrated in Figure \ref{fig:hyp}.
\end{remark}
\fi

\ifRR
For any $p\in[0;1]$, $\ya\in\Y$, $\yb\in\Y$ and $y\in\Y,$ consider
        \begarlab{eq:varphi}
        \varphi_{p,\ya,\yb}(y) \eqdef p \ell(\ya,y) + (1-p) \ell(\yb,y)
        \endarlab
\ifRR When $y_1\neq y_2$, this is the risk of the prediction function identically equal to $y$ when 
the distribution generating the data satisfies
        $P[Y=\ya]=p=1-P[Y=\yb].$ 
The case $y_1= y_2$ corresponds to $P[Y=\ya=\yb]=1$\ifRR
and will not be of interest to us (since we will use this function
for $\ya=\ha(x)\neq \hb(x)=\yb$). 
\fi.
        
Through this distribution, the quantity
\else and \fi
        \begarlab{eq:phi}
        \phi_{\ya,\yb}(p) \eqdef \und{\inf}{y\in\Y} \varphi_{p,\ya,\yb}(y)\ifRR\else.\fi
        \endarlab
\ifRR can be viewed as the risk of the best constant prediction function.\fi

\ifRR
\begin{remark}
In the binary classification setting, when $\Y=\{-1;+1\}$, from the Bayes rule, the function $\big[x \mapsto 
\text{a minimizer of }\varphi_{P(Y=1|X=x),-1,+1}\big]$
is the best prediction function to the extent that it minimizes the risk $R$.
Section \ref{sec:edgedisc} provides other typical examples of loss functions and \cite{BarJorMcA06}
gives an exhaustive study of their links.
\end{remark}
\fi

For any $q_+$ and $q_-$ in $[0;1]$, introduce
        \begarlab{eq:psidef}
        \psi_{q_+,q_-,\ya,\yb}(\alpha)\\
        \qquad
                \eqdef\phi_{\ya,\yb}[\alpha q_+ +(1-\alpha)q_-]
                -\alpha \phi_{\ya,\yb}(q_+) 
                - (1-\alpha)\phi_{\ya,\yb}(q_-) 
        \endarlab

\ifRR\else
The function $\phi_{\ya,\yb}$ is concave since it is the infimum of concave (affine) functions.
As a direct consequence, the function $\psi_{q_+,q_-,\ya,\yb}$ is concave and non-negative.
By integration by parts, one can prove:
\fi
\else
For any $p,\qa,\qb,\alpha\in[0;1]$ and any $u,\ya,\yb\in\Y$, let
        \lbegarlab{eq:def3}
        \varphi_{p,\ya,\yb}(y) & \eqdef & p \ell(\ya,y) + (1-p) \ell(\yb,y)\\
        \phi_{\ya,\yb}(p) & \eqdef & \und{\inf}{y\in\Y} \varphi_{p,\ya,\yb}(y)\\
        \psi_{q_+,q_-,\ya,\yb}(\alpha) & \eqdef & \phi_{\ya,\yb}[\alpha q_+ +(1-\alpha)q_-]\\
        & & \qquad
                -\alpha \phi_{\ya,\yb}(q_+) 
                - (1-\alpha)\phi_{\ya,\yb}(q_-) 
        \rendarlab

\fi

\ifRR
\begin{lemma} \label{lem:phi}
\begin{enumerate}
\item For any $\ya\in\Y$, $\yb\in\Y$, $q_+\in[0;1]$ and $q_-\in[0;1]$, the functions 
$\phi_{\ya,\yb}$ and $\psi_{q_+,q_-,\ya,\yb}$ are concave\ifRR, and consequently
admit one-sided derivatives everywhere\fi. 
The function $\psi_{q_+,q_-,\ya,\yb}$ is non-negative.
\item Define the function $K_\alpha$ as $K_{\alpha}(t) = [(1-\alpha)t] \wedge [\alpha(1-t)].$ 
Let $q_-=p_- \wedge p_+$ and $q_+=p_- \vee p_+$.
Assume that the function $\phi_{\ya,\yb}$ is twice differentiable by parts on $[q_-;q_+]$
to the extent that there exist
$q_-=\beta_0<\beta_1<\cdots<\beta_d<\beta_{d+1}=q_+$ such that for any $\ell\in\{0,\dots,d\}$,
$\phi_{\ya,\yb}$ is twice differentiable on $]\beta_{\ell};\beta_{\ell+1}[$. 
For any $\ell\in\{1,\dots,d\}$, let $\Delta_\ell$ denote the difference between 
the right-sided and left-sided derivatives of $\phi_{\ya,\yb}$ at point $\beta_\ell$. We have
\ifRR   \begar
        \psi_{p_+,p_-,\ha,\hb}(\alpha) = - (p_+-p_-)^2 \int_0^1 K_\alpha(t) \phi''_{\ha,\hb}\big[tp_++(1-t)p_-\big] dt\\
        \qquad\qquad\qquad\qquad\qquad
                - |p_+-p_-| \sum_{\ell=1}^d K_\alpha\big(\frac{\beta_\ell-p_-}{p_+-p_-}\big) \Delta_\ell 
        \endar
In particular, we have \fi
        \begarlab{eq:psidev}
        \psi_{p_+,p_-,\ya,\yb}(1/2)\\
        \qquad\quad
                =  \frac{(p_+-p_-)^2}{2} \int_0^1 [t\wedge(1-t)] 
                \big| \phi''_{\ya,\yb}\big[tp_++(1-t)p_-\big] \big| dt\\
        \qquad\qquad
                + \demi \sum_{\ell=1}^d \big( |p_+-\beta_\ell| \wedge |\beta_\ell-p_-| \big) | \Delta_\ell |
        \endarlab
\end{enumerate}
\end{lemma}
\fi

\ifRR
\begin{proof}
\begin{enumerate}
\item 
The function $\phi_{\ya,\yb}$ is concave since it is the infimum of concave (affine) functions.
As a direct consequence, the function $\psi_{p_+,p_-,\ya,\yb}$ is concave and non-negative.
\item It suffices to apply the following lemma, 
that is proved in Section \ref{sec:prooflemconvex},
to the function $f:t \mapsto \phi_{\ya,\yb}[t p_+ + (1-t)p_-]$, which has critical points
on $t_\ell$ defined as $t_\ell p_+ + (1-t_\ell)p_-=\beta_\ell$.
\begin{lemma} \label{lem:convex}
Let $f:[0;1]\rightarrow \R$ be a function twice differentiable by parts to the extent that there exist
$0=t_0<t_1<\cdots<t_d<t_{d+1}=1$ such that for any $\ell\in\{0,\dots,d\}$,
$f$ is twice differentiable on $]t_{\ell};t_{\ell+1}[$. Assume that $f$ is continuous and admits left-sided derivatives
$f'_l$ and right-sided derivatives $f'_r(t_\ell)$ at the critical points $t_\ell$. For any $\alpha\in[0;1]$, we have
        \begar
        f(\alpha) -  \alpha f(1) - (1-\alpha) f(0)
                = - \int_0^1 K_\alpha(t) f''(t) dt\\
        \qquad\qquad\qquad\qquad\qquad\qquad\quad
                \; - \sum_{\ell=1}^d K_\alpha(t_\ell) [f'_r(t_\ell)-f'_l(t_\ell)]\\
        \endar
\end{lemma}
\end{enumerate}
\end{proof}
\fi

\begin{definition} \label{def:hyp2}
Let $\big\{ P_{\sigmav} : \sigmav \triangleq (\sigma_1,\dots,\sigma_m) \in \{-;+\}^m \big\}$ be a hypercube of distributions.
\begin{enumerate}
\item The positive integer $m$ is called the \emph{dimension} of the hypercube.
\item The probability $w\eqdef\mu(\X_1)= \dots = \mu(\X_m)$ is called the \emph{edge probability}.
\item The \emph{characteristic function of the hypercube} is the function $\tpsi:\Rp \rightarrow\Rp$
defined as for any $u\in\Rp$
        \begar 
        \tpsi(u) & = & \demi m (u + 1) \expec{X}{\mu}\Big\{ \ds1_{X\in\X_1} 
                \psi_{p_+(X),p_-(X),\ha(X),\hb(X)}\big( \frac{u}{u+1} \big) \Big\}\\
	& = & \frac{mw}{2} (u + 1) \expecd{X}{\mu(\bullet|\X_1)} 
		\psi_{p_+,p_-,\ha,\hb}\big( \frac{u}{u+1} \big)  \qquad \text{in short.}
        \endar 
\item The \emph{edge discrepancies of type I} of the hypercube are 
        \ifRR \lbegarlab{eq:defed1a}
        \da & \eqdef & \frac{\tpsi(1)}{m w} = \expecd{X}{\mu(\bullet|\X_1)} \psi_{p_+,p_-,\ha,\hb}( 1/2 )\\
        \da'& \eqdef & \expecd{X}{\mu(\bullet|\X_1)} (p_+-p_-)^2
        \rendarlab
        \else
        \begarlab{eq:defed1a}
        \da \eqdef \frac{\tpsi(1)}{m w} = \expecd{X}{\mu(\bullet|\X_1)} \big[ \psi_{p_+,p_-,\ha,\hb}( 1/2 ) \big]
        \endarlab
        \fi
\item The \emph{edge discrepancy of type II} of the hypercube is defined as 
        \begarlab{eq:defed1}
        \ifRR
        \db \eqdef \expecd{X}{\mu(\bullet|\X_1)}\big[
                \sqrt{\pa(1-\pb)}-\sqrt{(1-\pa)\pb} \,\big]^2.
        \else 
        \db \eqdef \expecd{X}{\mu(\bullet|\X_1)}\Big\{ 
                \Big( \sqrt{\pa(1-\pb)}-\sqrt{[1-\pa]\pb} \Big)^2\Big\}.
        \fi
        \endarlab
\item 
A probability distribution $P_0$ on $\Z$ satisfying
        $P_0(dX)=\mu(dX)$ and for any $x\in\X-\X_0$,
        $P_0\big[Y=\ha(x)|X=x\big]=\demi=P_0\big[Y=\hb(x)|X=x\big]$
will be referred to as a \emph{base of the hypercube}.
\item Let $P_0$ be a base of the hypercube.
Consider distributions $P_{[\sigma]}, \sigma\in\{-,+\}$ admitting
the following density $\wrt$ $P_0$:
        \begar
        \frac{P_{[\sigma]}}{P_0}(x,y) = \left\{ \begin{array}{lll}
        2 p_\sigma(x) & \text{when } x\in\X_1 \text{ and } y=\ha(x)\\
        2 [ 1-p_\sigma(x)] & \text{when } x\in\X_1 \text{ and } y=\hb(x)\\
        1 & \text{otherwise}
        \end{array} \right.
        \endar
The distributions $P_{[-]}$ and $P_{[+]}$ will be referred to as the \emph{representatives of the hypercube}.
\item When the functions $\pa$ and $\pb$ are constant on $\X_1$, the hypercube will be said
\emph{constant}.
\item When the functions $\pa$ and $\pb$ satisfies $p_{+}=1-p_{-}$ on $\X-\X_0$, 
the hypercube will be said \emph{symmetrical}. In this case, the function 
$2\pa-1$ will be denoted $\xi$ so that
        \ifRR
        \lbegarlab{eq:xi}
        \pa& =& \frac{1+\xi}{2}\\
        \pb& =& \frac{1-\xi}{2}
        \rendarlab
        \else
        $\pa = \frac{1+\xi}{2}$
        and $\pb = \frac{1-\xi}{2}$.
        \fi
Otherwise it will be said \emph{asymmetrical}.
\ifRR\else\item A $(\tilde{m},\tilde{w},\tilde{\db})$-hypercube is a constant and symmetrical $\tilde{m}$-dimensional 
hypercube with edge probability $\tilde{w}$ and edge discrepancy of type II equal to $\tilde{\db}$,
and for which $\pa > 1/2$ and $\ha$ and $\hb$ are constant functions. 
\fi
\end{enumerate}
\end{definition}

The edge discrepancies are non-negative quantities that are 
all the smaller as $p_-$ and $p_+$ become closer. 
\ifRR\else
For a $(\tilde{m},\tilde{w},\tilde{\db})$-hypercube, we have
$\xi = \sqrt{\db}$, $\pb = \frac{1-\sqrt{\db}}{2}$ and $\pa~=~\frac{1+\sqrt{\db}}{2}$
and from \eqref{eq:psidev}, when the function $\phi_{\ha,\hb}$ is twice differentiable on
        $]\pb;\pa[$,
        \begarlab{eq:dacstsym}
        \da = \frac{\db}{2} 
                \int_{0}^1 [t\wedge (1-t)]\Big| \phi_{\ha,\hb}''\Big(\frac{1-\sqrt{\db}}{2}+\sqrt{\db}t\Big)\Big| \, dt.
        \endarlab               
\fi
\ifRR
Let us introduce the following assumption.

\vspace{0.3cm}
\noindent {\bf Differentiability assumption.}  
For any $x\in\X_1$, the function $\phi_{\ha(x),\hb(x)}$ is twice differentiable 
and satisfies for any $t\in[p_-(x) \wedge p_+(x);p_-(x) \vee p_+(x)]$,
        \begarlab{eq:assump1}
        \big| \phi''_{\ha(x),\hb(x)}(t) \big| \ge \zeta
        \endarlab
for some $\zeta>0$.
\vspace{0.3cm}

When $\Y\subseteq\R$, the differentiability assumption is typically fulfilled when
for any $\ya\neq\yb$, the functions
$y\mapsto \ell(\ya,y)$ and $y\mapsto \ell(\yb,y)$ admit second derivatives
lower bounded by a positive constant
and when these functions are minimum for respectively $y=\ya$ and $y=\yb$.
This is the case for least square loss and entropy loss, but it is not the case for
hinge loss, absolute loss or classification loss.
\fi
\ifRR The following result gives the main properties of the characteristic function~$\tpsi$
 and useful lower bounds of it.
\else
The function $\phi_{\ya,\yb}$ is concave since it is the infimum of concave (affine) functions.
Consequently, $\psi_{p_+,p_-,\ya,\yb}$ are concave and non-negative.
As a concave and non-negative

\fi

\begin{lemma} \label{lem:tpsi}
The characteristic function of the hypercube is a concave nondecreasing function
and satisfies 
\ifRR 
\begin{itemize}
\item $\tpsi(0)=0$ 
\item $\tpsi(u) \ge (u\wedge 1)\tpsi(1)$
\item Under the differentiability assumption (see \eqref{eq:assump1}), we have 
        \begar
        \tpsi(u) \ge \frac{mw\zeta}{4} \da' \frac{u}{u+1} \ge \frac{mw\zeta}{8} \da' (u\wedge 1)\ifRR\else.\fi
        \endar 
\ifRR where we recall that 
        $
        \da' \eqdef \expecd{X}{\mu(\bullet|\X_1)} (p_+-p_-)^2.
        $ \fi
In particular, we have
        \begarlab{eq:dasimp}
        \da \ge \frac{\zeta}{8} \da'.
        \endarlab
\end{itemize}
\else
$\tpsi(0)=0$ and $\tpsi(u) \ge (u\wedge 1)\tpsi(1)$.
\fi
\end{lemma}

\begin{proof}
\ifRR From Lemma \ref{lem:phi}, the function $\psi$ is non-negative and concave. Therefore 
the characteristic function is also concave and non-negative on $\Rp$. Consequently,
it is nondecreasing. The remaining assertions of the lemma are then straightforward.
\else Since the function $\psi$ is non-negative and concave, so is  
the characteristic function. Consequently,
$\tpsi$ is nondecreasing. The remaining assertions of the lemma are then straightforward.
\fi 
\end{proof}

\ifRR
To underline the link between the discrepancies of types I and II, 
one may consider \eqref{eq:dasimp} jointly with the following result

\begin{lemma} \label{lem:hypcst}
When a hypercube is constant and symmetrical, i.e. when on $\X_1$
$\pa = \frac{1+\xi}{2}$ and $\pb =\frac{1-\xi}{2}$
for $\xi$ constant, we have
        \begar
        \da'=\db=\xi^2.
        \endar
\end{lemma}
\fi

Finally, since the design of constant and symmetrical hypercubes is the key of numerous lower bounds,
we use the following:

\begin{definition} \label{def:mwxi}
A $(\tilde{m},\tilde{w},\tilde{\db})$-hypercube is a constant and symmetrical $\tilde{m}$-dimensional 
hypercube with edge probability $\tilde{w}$ and edge discrepancy of type II equal to $\tilde{\db}$,
and for which $\pa > 1/2$ and $\ha$ and $\hb$ are constant functions. 
\end{definition}

For these hypercubes, we have $m=\tilde{m}$, $w=\tw$, $\db=\tdb$ and
        \lbegar
        \xi & \equiv & \sqrt{\db}\\
        \pb & \equiv & \frac{1-\sqrt{\db}}{2}\\
        \pa & \equiv & \frac{1+\sqrt{\db}}{2}
        \rendar
and from \eqref{eq:psidev}, when the function $\phi_{\ha,\hb}$ is twice differentiable on
        $]\pb;\pa[$,
        \begarlab{eq:dacstsym}
        \da = \frac{\db}{2} 
                \int_{0}^1 [t\wedge (1-t)]\Big| \phi_{\ha,\hb}''\Big(\frac{1-\sqrt{\db}}{2}+\sqrt{\db}t\Big)\Big| \, dt.
        \endarlab               

\fi

\ifRR

\subsectionjy{$f$-similarity} \label{sec:sim}

Let us introduce a similarity measure between probability distributions.
When a probability distribution $\P$ is absolutely continuous $\wrt$ another probability distribution $\Q$,
i.e. $\P\ll \Q$, $\frac{\P}{\Q}$ denotes the density of $\P$ $\wrt$~$\Q$.

\begin{definition}
Let $f:\Rp \rightarrow \Rp$ be a concave function.
The \emph{$f$-similarity} between two probability distributions is defined as 
        \begarlab{eq:fsim}
        \S_f(\P,\Q) = \left\{ \begin{array}{lll}
                \int f\big(\frac{\P}{\Q}\big) d\Q & \text{if }\P\ll\Q \\
                f(0) & \text{otherwise}
        \end{array} \right.
        \endarlab
\end{definition}

Equivalently, if $p$ and $q$ denote the density of $\P$ and $\Q$ $\wrt$ the probability distribution $(\P+\Q)/2$,
one can define the $f$-similarity as
        \begar
        \S_f(\P,\Q) = \int_{q>0} f\big(\frac{p}{q}\big) d\Q
        \endar

It is called $f$-similarity in reference to $f$-divergence (see \cite{Csi67}) to which it is closely related.
Precisely, 
\ifRR
introduce the function $\tildf =f(1)-f$. $\tildf$ is convex and satisfies $\tildf(1)=0$.
Thus it is associated with an $\tilde{f}$-divergence, which is defined as
        \begarlab{eq:fdiv}
        D_{\tildf}(\P,\Q) = \left\{ \begin{array}{lll}
                \int \tildf\big(\frac{\P}{\Q}\big) d\Q & \text{if }\P\ll\Q \\
                \tildf(0) & \text{otherwise}
        \end{array} \right.
        \endarlab
Then we have $\S_f(\P,\Q) = f(1) - D_{\tildf}(\P,\Q).$ 
\else
we have $\S_f(\P,\Q) = f(1) - D_{\tildf}(\P,\Q)$ where $D_{\tildf}$ is the divergence
associated with the function $\tildf =f(1)-f$.
\fi

Here we use $f$-similarities since they are the quantities that naturally appear when developing 
our lower bounds.
\ifRR As the $f$-divergence, the $f$-similarity is in general asymmetric in $\P$ and $\Q$. 
Nevertheless for a concave function $f:\Rp \rightarrow \Rp$, one may define a concave function
$f^*:\Rp \rightarrow \Rp$ as $f^*(u)=uf(1/u)$ and $f^*(0)=\lim_{u\rightarrow 0} uf(1/u)$,
and we have (see \cite{Csi67} for the equivalent result for $f$-divergence): when $\P\ll\Q$ and $\Q\ll\P$,
        \begarlab{eq:simasymm}
        \S_f(\P,\Q)=\S_{f^*}(\Q,\P).
        \endarlab
\fi
We will use the following properties of $f$-similarities\ifRR\else \;(proof omitted)\fi.

\begin{lemma} \label{lem:sim}
Let $\P$ and $\Q$ be two probability distributions on a measurable space $(\calE,\B)$
such that $\P\ll\Q$. 
\begin{enumerate}
\item Let $f$ and $g$ be non-negative concave functions defined on $\Rp$ and let $a$ and $b$ be
non-negative real numbers. For any probability distributions $\P$ and $\Q$, we have
        $\S_{af+bg}(\P,\Q) = a \S_f(\P,\Q) + b \S_g(\P,\Q)$. 
Besides if $f\le g$, then
        $\S_f(\P,\Q) \le \S_g(\P,\Q)$. 
\item \label{item:2} Let $A\in\B$ such that $\frac{\P}{\Q}=1$ on $A$.
Let $A^c=\calE-A$. Let $f:\Rp \rightarrow \Rp$ be a concave function.
Let $\P'$ and $\Q'$ be probability distributions
on $\calE$ such that 
\begin{itemize}
\item $\P'\ll\Q'$ 
\item $\frac{\P'}{\Q'}=1$ on $A$.
\item $\P'=\P$ and $\Q'=\Q$ on $A^c$, i.e.
for any $B\in\B$, $\P'(B \cap A^c) = \P(B \cap A^c)$ 
and $\Q'(B \cap A^c) = \Q(B \cap A^c)$.
\end{itemize}
We have
        \begar
        \S_{f}(\P',\Q')=\S_{f}(\P,\Q).
        \endar
\item Let $\calE'$ be a measurable space and $\mu$ be a $\sigma$-finite positive measure on $\calE'$.
Let $\{f_v\}_{v\in\calE'}$ be a family of non-negative concave functions defined on $\Rp$ such 
that $(u,v) \mapsto f_v\big(\frac{\P}{\Q}(u)\big)$ is measurable. We have
        \begarlab{eq:sfubini}
        \int_{\calE'} \S_{f_v}(\P,\Q) \mu(dv)= \S_{\int_{\calE'} f_v \mu(dv)}(\P,\Q).
        \endarlab
\end{enumerate}
\end{lemma}
\ifRR
\begin{proof}
\begin{enumerate}
\item It directly follows from the definition of $f$-similarities.
\item We have
        \begar
        \S_{f}(\P',\Q') & = & \int_A f\big(\frac{\P'}{\Q'}\big) d\Q'
                + \int_{A^c} f\big(\frac{\P'}{\Q'}\big) d\Q'\\
        & = & f(1) \Q'(A) + \int_{A^c} f\big(\frac{\P}{\Q}\big) d\Q\\
        & = & f(1) \Q(A) + \int_{A^c} f\big(\frac{\P}{\Q}\big) d\Q\\
        & = & \S_{f}(\P,\Q) 
        \endar

\item This Fubini's type result follows from the definition of the integral of non-negative functions
on a product space.
\end{enumerate}
\end{proof}
\fi
\fi

\subsectionjy{Generalized Assouad's lemma} \label{sec:genassouad}

\ifRR We recall that the n-fold product of a distribution $P$ is denoted~$P^{\otimes n}$. \fi 
We start this section with a general lower bound for hypercubes of 
distributions\ifRR \;(as defined in Section \ref{sec:hyper})\fi. This lower bound is expressed in terms 
of a similarity \ifRR(as defined in Section \ref{sec:sim}) \fi between 
$n$-fold products of representatives of the hypercube.

\begin{theorem} \label{th:assouad}
Let $\calP$ be a set of probability distributions containing a hypercube of distributions 
of characteristic function $\tpsi$ and representatives $P_{[-]}$ and $P_{[+]}$.
For any training set size $n\in\N^*$ and any estimator $\hg$, we have
        \begarlab{eq:thass1}
        \underset{P \in \calP}{\sup} \big\{ \E R( \hg ) 
                - \und{\min}{g} R(g) \big\} \ge 
                \S_{\tpsi} \big( P_{[+]}^{\otimes n} , P_{[-]}^{\otimes n} \big)
        \endarlab
where the minimum is taken over the space of all prediction functions
and $\E R(\hg)$ denotes
the expected risk of the estimator $\hg$ trained on a sample of size $n$:
$\E R(\hg)=\expec{Z_1^n}{P^{\otimes n}} R(\hg_{Z_1^n})
=\expec{Z_1^n}{P^{\otimes n}} \undc{\E}{(X,Y)\sim P}\ell[Y,\hg_{Z_1^n}(X)].$
\end{theorem}

\jyproof{th:assouad}

This theorem provides a lower bound holding for 
any estimator and expressed in terms of the hypercube structure.
\ifRR To obtain a tight lower bound associated with a particular learning task,
it then suffices to find
the hypercube in $\calP$ for which the $\rhs$ of \eqref{eq:thass1} is the largest possible.
\fi
By providing lower bounds of $\S_{\tpsi} \big( P_{[+]}^{\otimes n} , P_{[-]}^{\otimes n} \big)$
that are more explicit $\wrt$ the hypercube parameters, we obtain 
\ifRR 
the following results that are more in a ready-to-use form than Theorem \ref{th:assouad}.
\fi

\begin{theorem} \label{th:appassouad}
Let $\calP$ be a set of probability distributions containing a hypercube of distributions 
characterized by its dimension $m$, its edge probability $w$ and 
its edge discrepancies $\da$ and $\db$ (see Definition \ref{def:hyp2}).
For any estimator $\hg$ and training set size $n\in\N^*$, the following assertions hold.
\begin{enumerate}
\item We have
        \begarlab{eq:w1a}
        \underset{P \in \calP}{\sup} \big\{ \E R( \hg ) 
                - \und{\min}{g} R(g) \big\} & \ge & mw \da \big(1-\sqrt{1 - [1-\db]^{nw}}\big)\\
                & \ge & mw \da \big(1-\sqrt{nw\db}\big).
        \endarlab
\item When the hypercube is constant and symmetrical\ifRR (see Definition \ref{def:hyp2}) \fi,
we also have  
        \begarlab{eq:w1b}
        \underset{P \in \calP}{\sup} \big\{ \E R( \hg ) 
                - \und{\min}{g} R(g) \big\} 
                \ge mw \da \Big\{ \P\Big(|N| > \sqrt{\frac{nw\db}{1-\db}} \Big) - \db^{1/4} \Big\}
        \endarlab
for $N$ a centered gaussian random variable with variance $1$.
\item When the hypercube satisfies $\pa\equiv 1 \equiv 1-\pb$, we also have
        \begarlab{eq:w1c}
        \underset{P \in \calP}{\sup} \big\{ \E R( \hg ) 
                - \und{\min}{g} R(g) \big\} 
                \ge mw \da (1-w)^n
        \endarlab
\ifRR
\item When the hypercube is constant and symmetrical and when the differentiability assumption (see \eqref{eq:assump1})
holds, we also have
        \begarlab{eq:w2}
        \underset{P \in \calP}{\sup} \big\{ \E R( \hg ) 
                - \und{\min}{g} R(g) \big\} 
                \ge \frac{mw \zeta \db}{8}\\ 
         \qquad\qquad
                \times \Big\{ 1+\demi \big[1-\big(1-\sqrt{1-\db}\big)w\big]^n 
                - \demi \big[1+\big(\frac{1+\db}{\sqrt{1-\db}}-1\big)w\big]^n \Big\}
        \endarlab
\fi
\end{enumerate}
\end{theorem}

\begin{proof}
See Section \ref{sec:proofappassouad}.
\end{proof}

\ifRR
The lower bounds \eqref{eq:w1a}, \eqref{eq:w1b} and \eqref{eq:w2} are of the same nature.
\eqref{eq:w1a} is the general lower bound having the simplest form. For constant and symmetrical
hypercubes, it can be refined into \eqref{eq:w1b} and, when the differentiability assumption holds, into \eqref{eq:w2}.
These refinements mainly concern constants as we will see in Section \ref{sec:numerical}.
Finally, \eqref{eq:w1c} is less general but provide results with tight constants when 
convergence rate of order $n^{-1}$ has to be proven 
(see Remarks \refp{rem:common} and \refp{rem:commop}).


To better understand the link between \eqref{eq:w1a}, \eqref{eq:w1b} and \eqref{eq:w2}, the following corollary considers
an asymptotic setting in which $n$ goes to infinity and the parameters of the hypercube varies with $n$ (which
is the typical situation even for finite sample lower bounds).
\fi
\begin{corollary} \label{cor:appassouad}
\ifRR Let $a>0$ and $N$ be a centered gaussian random variable with variance $1$.
Under the assumptions of item 4 of the previous theorem, we have
        \begarlab{eq:dadb}
        \da \ge \frac{\zeta}{8} \db 
        \endarlab
and \eqref{eq:w1a}, \eqref{eq:w1b} and \eqref{eq:w2} respectively lead to
\fi
        \begarlab{eq:asw1a}
        \und{\liminf}{\db \rightarrow 0,nw\db \rightarrow a} \; \frac{1}{mw \da} \;
                \underset{P \in \calP}{\sup} \big\{ \E R( \hg ) 
                - \und{\min}{g} R(g) \big\} & \ge & 1-\sqrt{1 - e^{-a}} \\& \ge & 1-\sqrt{a},
        \endarlab
        \begarlab{eq:asw1b}
        \und{\liminf}{\db \rightarrow 0,nw\db \rightarrow a} \; \frac{1}{mw \da} \;
                \underset{P \in \calP}{\sup} \big\{ \E R( \hg ) 
                - \und{\min}{g} R(g) \big\} 
                \ge \P(|N| > \sqrt{a} ) 
        \endarlab
\ifRR and
        \begarlab{eq:asw2}
        \und{\liminf}{\db \rightarrow 0,nw\db \rightarrow a} \; \frac{8}{mw \zeta \db} \;
                \underset{P \in \calP}{\sup} \big\{ \E R( \hg ) - \und{\min}{g} R(g) \big\} \ge 
                1+\frac{e^{-a/2}}{2} - \frac{e^{3a/2}}{2}
        \endarlab
\fi
\end{corollary}

\begin{proof}
It follows from Theorem \ref{th:appassouad} and $(1+x)^{1/x} \rightarrow e$
when $x\rightarrow 0$.
\end{proof}

\ifRR Inequality \eqref{eq:dadb} leads to (slightly) weakened versions of 
\eqref{eq:asw1a} and \eqref{eq:asw1b} that can be directly compared with \eqref{eq:asw2} 
(see Figure~\ref{fig:lb}). 
A numerical comparison of these bounds is given in Section \ref{sec:numerical}.
\fi

\ifRR
\begin{figure}[ht] 
\begin{center}
\includegraphics{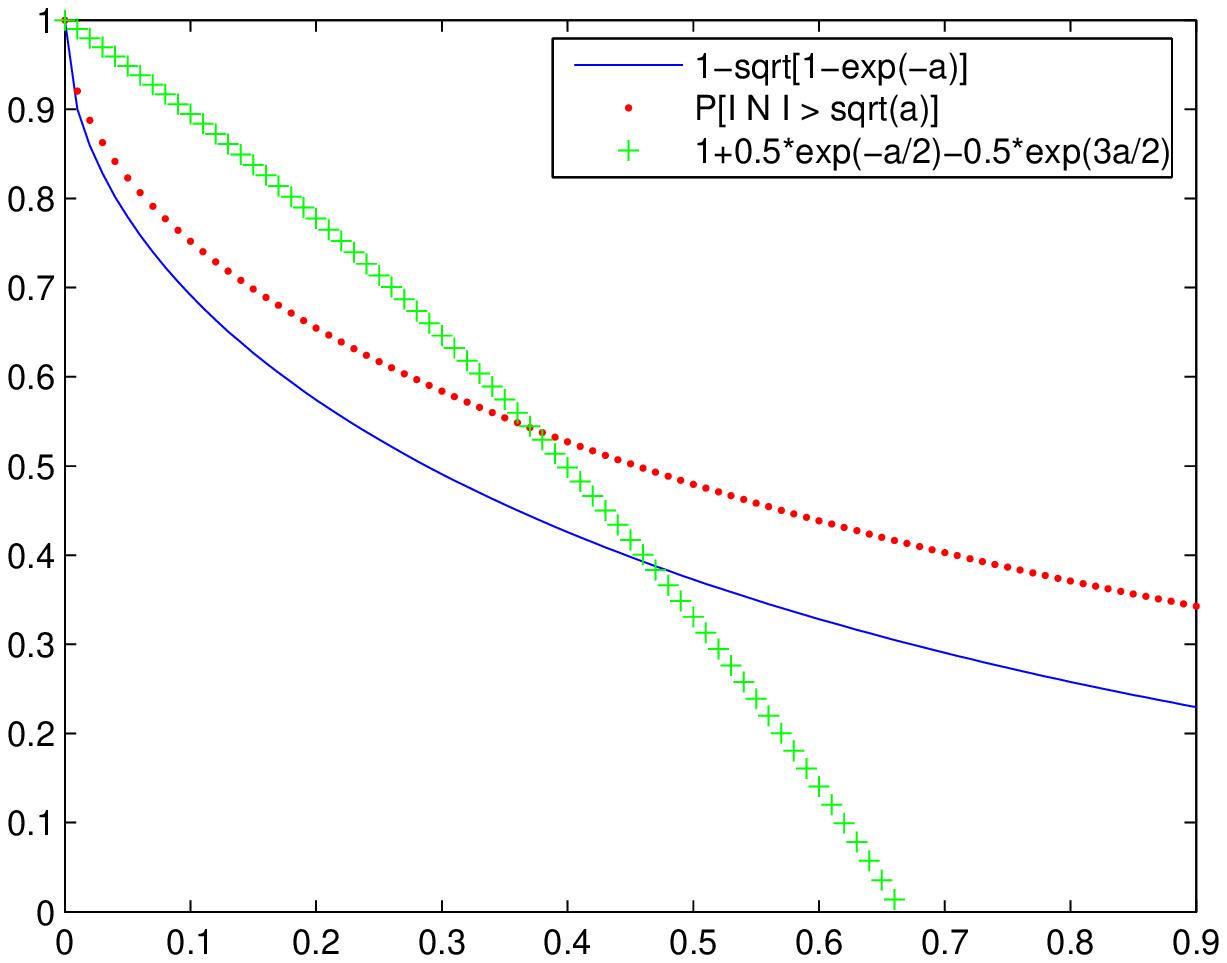}
\caption{Comparison of the $\rhs$ of \eqref{eq:asw1a}, \eqref{eq:asw1b} and \eqref{eq:asw2}} \label{fig:lb}
\end{center}
\end{figure}
\fi

\ifRR
\begin{remark}
The previous lower bounds consider deterministic estimators (or algorithms), i.e. functions from
the training set space $\cup_{n \ge 0} \Z^n$ to the prediction function space 
$\barG$. They still hold for randomized estimators, i.e. functions from the training set space 
to the set $\D$ of probability distributions on $\barG$. 
\end{remark}
\fi

\subsectionjy{Examples} \label{sec:examples}

Theorem \ref{th:appassouad} motivates the following simple strategy to obtain a lower bound
for a given set $\calP$ of probability distributions and a reference set $\G$ of prediction functions: it consists 
in looking for the hypercube contained in the
set $\calP$ and for which 
        \begin{itemize}
        \item the lower bound is maximized,
        \item 
        for any distribution of the hypercube, $\G$ contains a best prediction function, i.e.
        ${\min}_{g} R(g)= {\min}_{g\in \G} R(g)$.
        \end{itemize}
In general, the order of the bound is given by the quantity $m w \da$ 
\ifRR (or $mw\zeta \db$ in the case of \eqref{eq:w2}) \fi
and the quantities $w$ and $\db$ are taken such that $n w \db$ is of order $1$. 
\ifRR 

In this section, we apply this strategy in different learning tasks.
Before giving these lower bounds (Section \ref{sec:precisex}),
Section \ref{sec:edgedisc} stresses on the influence of the loss function in the computations
of the edge discrepancy $\da$ and the constant $\zeta$ of the differentiability assumption
\eqref{eq:assump1}.

\subsubsection{Edge discrepancy $\da$ and constant $\zeta$ of the differentiability assumption} \label{sec:edgedisc}

All the previous lower bounds rely on either the edge discrepancy $\da$ or the constant $\zeta$ 
of the differentiability assumption (that has been introduced to control $\da$ in a simple way).
The aim of this section is to provide more explicit formulas of these quantities 
for different loss functions.

To obtain the formula for $\da$, we essentially use \eqref{eq:defed1a} and \eqref{eq:psidev}
jointly with the explicit computation of the second derivative of the function $\phi$.

\paragraph{\bf Entropy loss.} \label{ex:entropy}
Here we consider $\Y=[0;1]$ and the loss for prediction $y'$ instead of $y$
is $\ell(y,y')=K(y,y')$, where $K(y,y')$ is the Kullback-Leibler divergence between Bernoulli distributions 
with respective parameters $y$ and $y'$, i.e.
        $K(y,y') = y \log \big( \frac{y}{y'} \big) + (1-y) \log \big( \frac{1-y}{1-y'} \big).$
Let $H(y)$ denote the Shannon's entropy of the Bernoulli distribution with parameter $y$, i.e.
        \begarlab{eq:shannon}
        H(y)=- y \log y - (1-y) \log (1-y).
        \endarlab
Computations lead to: for any $p\in[0;1]$,
        \begarlab{eq:exphientropy}
        \phi_{\ya,\yb}(p) = H\big( p y_1 + (1-p) y_2 \big) - p H(y_1) - (1-p) H(y_2),
        \endarlab
hence 
        \begar 
        \phi''_{\ya,\yb}(p) = -\frac{(y_1-y_2)^2}{[py_1+(1-p)y_2][(p(1-y_1)+(1-p)(1-y_2)]}.
        \endar 
This last equality is useful when one wants to compute $\zeta$ satisfying \eqref{eq:assump1}.

\paragraph{\bf Classification loss.} 
In this setting, we have $|\Y|<+\infty$ and the loss incurred by predicting
$y'$ instead of the true value $y$ is $\ell(y,y')=\ds1_{y\neq y'}$.
In this learning task, we have $\phi_{\ya,\yb}(p) = [p \wedge (1-p)] \ds1_{\ya\neq\yb}$
and $\phi''_{\ya,\yb} \equiv 0$ on $[0;1]-\{\demic\}$.
Then \eqref{eq:defed1a}, \eqref{eq:psidev} and Remark \refp{rem:haneqhb} lead to
        \begar
        \da = \E_{\mu(\bullet|\X_1)}\Big\{\big[ \big|p_+-\demi\big| \wedge \big|p_--\demi\big|\big]
                \ds1_{(p_+-\demi) (p_--\demi) < 0} \Big\}.
        \endar

\paragraph{\bf Binary classification losses (or regression losses when the output is binary).} 
In this setting, we have $\Y=\R \cup \{-\infty;+\infty\}$, but we know that $P(Y\in\{-1;+1\})=1$.
So we are only interested in hypercubes of distributions satisfying
this constraint, i.e. such that for any $x\in\X$, $\ha(x)$ and $\hb(x)$ belong to $\{-1;+1\}$.
In this setting, a best prediction function $g$, i.e. a measurable function from $\X$ to $\Y$ minimizing
$R(g) = \E \, \ell[Y,g(X)]$, is determined by the regression function: 
        \begar
        \eta(x) = P(Y=+1|X=x).
        \endar
Let $\sign(x)=\ds1_{x\ge 0}-\ds1_{x< 0}$ be the sign function on $\R$.
\begin{itemize}
\item {\bf $\R$-Classification loss}. The loss function is $\ell(y,y')=\ds1_{yy'<0}$
and a best prediction function is $g^*(x) = \sign(\eta(x)-\demic).$
Without surprise, we recover the same formulae as for the classification loss. 
\item {\bf Hinge loss}. The loss function is $\ell(y,y')=(1-yy')_+=\max\{0;1-yy'\}$
and a best prediction function is $g^*(x) = \sign(\eta(x)-\demic)$.
For any $\ya,\yb\in\{-1;+1\}$, we have
        $\phi_{\ya,\yb}(p) = 2 [p \wedge (1-p)] \ds1_{\ya\neq\yb}$
and $\phi''_{\ya,\yb} \equiv 0$ on $[0;1]-\{\demic\}$.
Then \eqref{eq:defed1a}, \eqref{eq:psidev} and Remark \refp{rem:haneqhb} lead to
        \begar
        \da = 2 \E_{\mu(\bullet|\X_1)}\Big\{\big[ \big|p_+-\demi\big| \wedge \big|p_--\demi\big|\big]
                \ds1_{(p_+-\demi) (p_--\demi) < 0} \Big\},
        \endar
\item {\bf Exponential loss (or AdaBoost loss)}. The loss function is $\ell(y,y')=e^{-yy'}$.
For any $\ya\neq\yb\in\{-1;+1\}$ and any $p\in [0;1]$, the function $\varphi_{p,\ya,\yb}$ is minimized
for $y=\frac{y_1}{2}\log\big(\frac{p}{1-p}\big)$, so a best prediction function is 
        $g^*(x) =\demi\log\big(\frac{\eta(x)}{1-\eta(x)}\big).$
We obtain
        $\phi_{\ya,\yb}(p) = 2 \sqrt{p(1-p)} \ds1_{\ya\neq\yb}$
and 
        \begar
        \phi''_{\ya,\yb}(p) = - \frac{1}{2 [p(1-p)]^{3/2}}\ds1_{\ya\neq\yb}.
        \endar
In this setting, to obtain a lower bound of $\da$, one has typically to compute $\zeta$
satisfying \eqref{eq:assump1} and to use \eqref{eq:dasimp}.
\item {\bf Logit loss}. The loss function is $\ell(y,y')=\log(1+e^{-yy'})$.
For any $\ya\neq\yb\in\{-1;+1\}$ and any $p\in [0;1]$, the function $\varphi_{p,\ya,\yb}$ is minimized
for $y=y_1\log\big(\frac{p}{1-p}\big)$, so a best prediction function is 
        $g^*(x) = \log\big(\frac{\eta(x)}{1-\eta(x)}\big).$
We obtain
        $\phi_{\ya,\yb}(p) = H(p)\ds1_{\ya\neq\yb}$,
where $H(p)$ denote the Shannon's entropy of the Bernoulli distribution with parameter $p$ (see \eqref{eq:shannon}).
We get 
        \begar
        \phi''_{\ya,\yb}(p) = - \frac{1}{p(1-p)}\ds1_{\ya\neq\yb}.
        \endar
Once more, to obtain a lower bound of $\da$, one has typically to compute $\zeta$
satisfying \eqref{eq:assump1} and to use \eqref{eq:dasimp}.
\end{itemize}

\paragraph{\bf $L_q$-loss.} \label{ex:qpower} 
We consider $\Y=\R$ and the loss function is $\ell(y,y')=|y-y'|^q$
with $q\ge 1$. The values $q=1$ and $q=2$ respectively correspond to the absolute loss and 
the least square loss. 
\begin{itemize}
\item {\it Case $q=1$ :} 
Due to the lack of strong convexity of the loss function, the absolute loss setting differs completely from what occurs
for $q>1$ and appears to be similar to the classification and hinge losses settings. Indeed computations lead to 
        $\phi_{\ya,\yb}(p) = [p \wedge (1-p)] |\yb-\ya|$
and $\phi''_{\ya,\yb} \equiv 0$ on $[0;1]-\{\demic\}$.
Then \eqref{eq:defed1a} and \eqref{eq:psidev} lead to
        \begarlab{eq:1power}
        \da = \E_{\mu(\bullet|\X_1)}\Big\{|\hb-\ha|\big[ \big|p_+-\demi\big| \wedge \big|p_--\demi\big|\big]
                \ds1_{(p_+-\demi) (p_--\demi) < 0} \Big\}.
        \endarlab
\item {\it Case $q>1$ :} 
Tedious computations put in Appendix \ref{app:proofphisecond} lead to: for any $p\in[0;1]$,
        \begarlab{eq:qphi}
        \phi_{y_1,y_2}(p) = p(1-p) \frac{|\yb-\ya|^q}{\big[p^{\frac{1}{q-1}}+(1-p)^{\frac{1}{q-1}}\big]^{q-1}}
        \endarlab
and 
        \begarlab{eq:qpower}
        \phi''_{y_1,y_2}(p) = - \frac{q}{q-1} [p(1-p)]^{\frac{2-q}{q-1}}
                \frac{|\yb-\ya|^q}{\big[p^{\frac{1}{q-1}}+(1-p)^{\frac{1}{q-1}}\big]^{q+1}}
        \endarlab
To obtain a lower bound of $\da$, as for the entropy loss, one has to compute $\zeta$
satisfying \eqref{eq:assump1} and to use \eqref{eq:dasimp}.
\item {\it Special case $q=2$ :} \label{ex:ls}
For the least square setting, the formulae simplify into $\phi_{y_1,y_2}(p) = p(1-p) |\yb-\ya|^2$
and $\phi''_{y_1,y_2}(p) = -2 |\yb-\ya|^2$. Then the edge discrepancy $\da$ can be written explicitly as
        \begarlab{eq:ls}
        \da = \frac{1}{4} \E_{\mu(\bullet|\X_1)}\big\{(p_+-p_-)^2 (\hb-\ha)^2 \big\}.
        \endarlab
\end{itemize}

\subsubsection{Various learning lower bounds} \label{sec:precisex}

Before giving learning lower bounds matching up to multiplicative constants the upper bounds
developed in the previous sections, we will start with two standard problems: 
classification lower bounds for Vapnik-Cervonenkis classes and uniform universal consistency. 

\paragraph{Binary classification.} 

We consider $\Y=\{0;1\}$ and $l(y,y')=\ds1_{y\neq y'}$. 
Since the work of Vapnik-Cervonenkis \cite{Vap74}, 
several lower bounds have been proposed and the most achieved ones are given in \cite[Chapter 14]{Dev96}. 
The following theorem provides a significant improvement of the constants of these bounds.

\begin{theorem} \label{th:vc}
Let $L \in [0;1/2]$, $n\in\N$ and $\G$ be a set of prediction functions of VC-dimension $V \ge 2$.
Consider the set $\calP_L$ of probability distributions on $\X\times\{0;1\}$
such that
	$\inf_{g\in\G} R(g) = L$.
For any estimator $\hg$: \\
$\bullet$ when $L = 0$, there exists $\P\in\calP_0$ for which	
	\begarlab{eq:vc1}
	\E R(\hg) - \und{\inf}{g\in\G} R(g) 
		\geq \left\{ \begin{array}{lll}
		\frac{V-1}{2e(n+1)} & \text{ when } n \geq V-2\\
		\frac{1}{2} \big( 1 - \frac{1}{V} \big)^n & \text{ otherwise} \\
		\end{array} \right..
	\endarlab
$\bullet$ when $0 < L \leq 1/2$, there exists $\P\in\calP_L$ for which	
	\begarlab{eq:vc2}
	\E R(\hg) - \und{\inf}{g\in\G} R(g) 
		\geq \left\{ \begin{array}{lll}
		\sqrt{\frac{L(V-1)}{32n}} \vee \frac{2(V-1)}{27 n} & \text{ when } 
			\frac{(1-2L)^2 n}{V} \geq \frac{4}{9}\\
		\frac{1-2L}{6} & \text{ otherwise}
		\end{array} \right..
	\endarlab
$\bullet$ there exists a probability distribution for which
	\begarlab{eq:vcnew}
	\E R(\hg) - \und{\inf}{g\in\G} R(g) 
		\geq \frac{1}{8} \sqrt{\frac{V}{n}} 
	\endarlab
\end{theorem}

\begin{proof}[Sketch]
We have $\phi_{\ya,\yb}(p) = [p \wedge (1-p)] \ds1_{\ya\neq\yb}$
and for constant symmetrical hypercubes $\da=\sqrt{\db}/2$.
Then \eqref{eq:vc1} comes from \eqref{eq:w1c} and the use of a $(V-1,1/(n+1),1)$-hypercube
and a $(V,1/V,1)$-hypercube.

To prove \eqref{eq:vc2}, from \eqref{eq:w1a} and the use of a 
	$\big(V-1,\frac{2L}{V-1},\frac{V-1}{8nL}\big)$-hypercube, a 
	$\big(V-1,\frac{4}{9n(1-2L)^2},(1-2L)^2\big)$-hypercube
and a $(V,1/V,(1-2L)^2)$-hypercube, we obtain
	\begar
	\E R(\hg) - \und{\inf}{g\in\G} R(g) 
		\geq \left\{ \begin{array}{lll}
		\sqrt{\frac{L(V-1)}{32n}} \text{ when } 
			\frac{(1-2L)^2 n}{V-1} \ge \frac{L}{2} \vee \frac{(1-2L)^2}{8L}\\
		\frac{2(V-1)}{27 n(1-2L)} \text{ when } 
			\frac{(1-2L)^2 n}{V-1} \ge \frac{4}{9}\\
		\frac{1-2L}{2} \Big( 1 -  \sqrt{\frac{(1-2L)^2 n}{V}} \Big) \text{ always} 
		\end{array} \right.,
	\endar
which can be weakened into \eqref{eq:vc2}. Finally, \eqref{eq:vcnew} comes from the last inequality and by choosing $L$
such that $1-2L=\frac12\sqrt{\fracc{V}n}$.
\end{proof}

\paragraph{No uniform universal consistency for general losses.}

This type of results is well known and tells that there is no guarantee of doing well on finite samples.
In classification setting, when the input space is infinite, i.e. $|\X|=+\infty$, 
by using a $(\lfloor n\alpha \rfloor , 1/\lfloor n\alpha \rfloor , 1 )$-hypercube with $\alpha$ tending to infinity, one can 
recover that: for any training sample size $n$, 
``any discrimination rule can have an arbitrarily bad probability of error for finite sample size'' (\cite{Dev82}), precisely:
	$$
	\inf_{\hg} \sup_{\P} \big\{ \P[Y \neq \hg(X)] 
                - \und{\min}{g} \P[Y \neq g(X)]  \big\} = 1/2,
    $$
where the infimum is taken over all (possibly randomized) classification rules.
For general loss functions, as soon as $|\X|=+\infty$, we can use
$(\lfloor n\alpha \rfloor , 1/\lfloor n\alpha \rfloor , 1 )$-hypercubes with $\alpha$ tending to infinity and obtain
	\bigbegarlab{eq:nounif}
	\inf_{\hg} \sup_{\P} \big\{ \E R(\hg) - \und{\inf}{g\in\G} R(g) \big\} \ge \und{\sup}{y_1,y_2\in\Y} \psi_{1,0,y_1,y_2}(1/2),
	\bigendarlab
where $\psi$ is the function defined in \eqref{eq:psidef}.

\paragraph{Entropy loss setting.}

We consider $\Y=[0;1]$ and $\ell(y,y')=K(y,y')$ (see p.\pageref{ex:entropy}).
We have seen in Section \ref{sec:online} that there exists an estimator $\hg$ such that
        \begarlab{eq:uentropy}
        \E R(\hg) - \und{\min}{g\in\G} R(g) \le \frac{\log |\G|}{n}
        \endarlab
The following consequence of \eqref{eq:w1c} shows that this result is tight.

\begin{theorem} \label{th:entropy}
For any training set size $n\in\N^*$, positive integer $d$ and 
input space $\X$ containing at least $\lfloor\log_2 (2d)\rfloor$ points, 
there exists a set $\G$ of $d$ prediction functions such that:
for any estimator $\hg$ there exists a probability distribution on the data space 
$\X\times[0;1]$ for which
        \begar
        \E R(\hg) - \und{\min}{g\in\G} R(g) \ge e^{-1} (\log 2) \big( 1 \wedge \frac{\lfloor \log_2 |\G| \rfloor}{n+1} \big)
        \endar
\end{theorem}

\begin{proof}
We use a $\big( \tm , \frac{1}{n+1} \wedge \frac{1}{\tm} , 1 \big)$-hypercube
with $\tm=\lfloor \log_2 |\G| \rfloor= \big\lfloor \frac{\log |\G|}{\log 2} \big\rfloor$, $\ha\equiv 0$ and $\hb\equiv 1$.
From \eqref{eq:psidef}, \eqref{eq:defed1a} and \eqref{eq:exphientropy}, we have
        \begar
        \da=\psi_{1,0,0,1}(1/2)
        =\phi_{0,1}(1/2)=H(1/2)=\log 2.
        \endar
From \eqref{eq:w1c}, we obtain
        \begar
        \E R(\hg) - \und{\min}{g\in\G} R(g) 
                \ge \big( \frac{\lfloor \log_2 |\G| \rfloor}{n+1} \wedge 1\big) 
                                (\log 2) \big(1-\frac{1}{n+1} \wedge \frac{1}{\lfloor \log_2 |\G| \rfloor}\big)^n 
                
        \endar
Then the result follows from $[1-1/(n+1)]^{n} \searrow e^{-1}$.
\end{proof}

\begin{remark} \label{rem:common}
For $|\G| < 2^{n+2}$, the lower bound matches the upper bound \eqref{eq:uentropy}
up to the multiplicative factor $e \approx 2.718$ . For $|\G| \ge 2^{n+2}$,
the size of the model is too large and, without any extra assumption,
no estimator can learn from the data.
To prove the result, we consider distributions for which the output is deterministic when knowing the input.
So the lower bound does not come from noisy situations but from situations in which different prediction
functions are not separated by the data to the extent that no input data falls into the (small) subset on which they are different.
\end{remark}

\paragraph{$L_q$-regression with bounded outputs.}   \label{sec:exlqb}
We consider $\Y=[-B;B]$ and $\ell(y,y')=|y-y'|^q$ (see p.\pageref{ex:qpower}).
The following two theorems are roughly summed up in Figure~\ref{fig:lq}
that represents the optimal convergence rate for $L_q$-regression.

\begin{figure}[ht]
\begin{center}
\includegraphics[width=9cm]{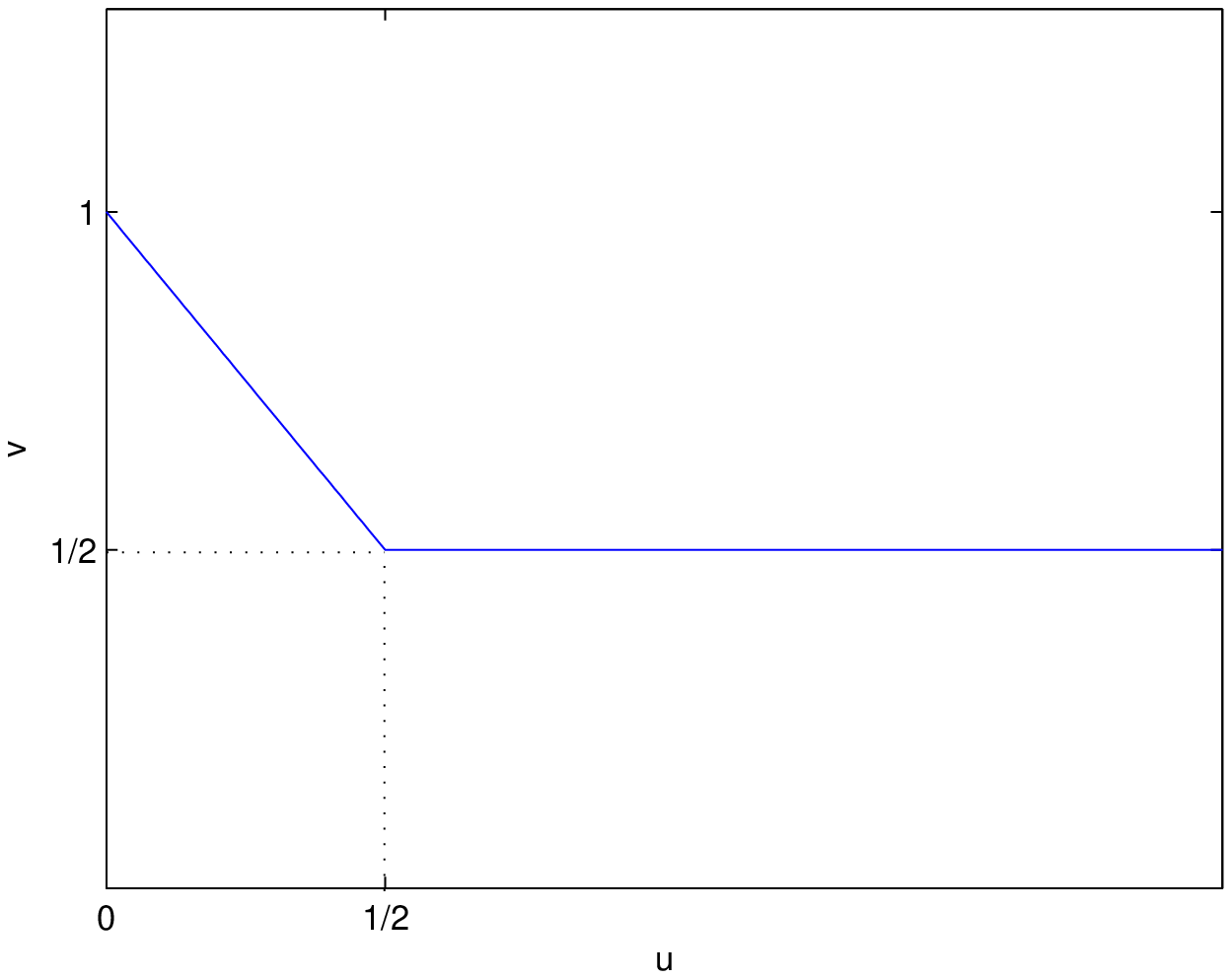}
\caption{{\it Influence of the convexity of the loss on the optimal convergence rate.}
Let $c>0$. 
We consider $L_q$-losses with $q=1+c \big(\frac{\log|\G|}{n}\big)^{u}$
for $u\ge 0.$ 
For such values of $q$, the optimal convergence rate of the associated learning task is of order 
$\big(\frac{\log|\G|}{n}\big)^{v}$ with $1/2 \le v \le 1$.
This figure represents the value of $u$ in abscissa and the value of $v$ in ordinate.
The value $u=0$ corresponds to constant $q$ greater than $1$.
For these $q$, the optimal convergence rate is of order $n^{-1}$
while for $q=1$ or ``very close'' to $1$, the convergence rate is
of order $n^{-1/2}$.
}  \label{fig:lq}
\end{center}
\end{figure}

\begin{itemize}
\item {\it Case $1\le q \le 1+\sqrt{\frac{\lfloor \log_2 |\G| \rfloor}{4n}\wedge 1}$: }
From \eqref{eq:corhoeff2}, there exists an estimator $\hg$ such that
        \begarlab{eq:u1b}
        \E R(\hg) - \und{\min}{g\in\G} R(g) \le 
                2^{\frac{2q-1}{2}} B^q \, \sqrt{\frac{\log |\G|}{n}}    
        \endarlab
The following corollary of Theorem \ref{th:appassouad} shows that this result is tight.

\begin{theorem} \label{th:lqpt}
Let $B>0$ and $d\in\N^*.$ For any training set size $n\in\N^*$ and 
any input space $\X$ containing at least $\lfloor\log_2 d\rfloor$ points, 
there exists a set $\G$ of $d$ prediction functions such that:
for any estimator $\hg$ there exists a probability distribution on the data space
$\X\times[-B;B]$ for which
        \begar
        \E R(\hg) - \und{\min}{g\in\G} R(g) 
                \ge \left\{ \begin{array}{lll}
                c_q B^q \sqrt{\frac{\lfloor \log_2 |\G| \rfloor}{n}} & \text{ if $\,|\G| < 2^{4 n+1}$}\\
                2 c_q B^q & \text{ otherwise}
                \end{array} \right.,
        \endar
where 
        \begar
        c_q = \left\{ \begin{array}{lll}
        1/4 & \text{\qquad if \, $q=1$}\\
        q/40 & \text{\qquad if \, $1< q \le 1+\sqrt{\frac{\lfloor \log_2 |\G| \rfloor}{4n}\wedge 1}$}
        \end{array} \right.
        \endar
\end{theorem}

\begin{proof}
See Section \ref{sec:proofthlqpg}.
\end{proof}

\item {\it Case $q> 1+\sqrt{\frac{\lfloor \log_2 |\G| \rfloor}{4n}\wedge 1}$ : }
We have seen in Section \ref{sec:online} that there exists an estimator $\hg$ such that
        \begarlab{eq:u2b}
        \E R(\hg) - \und{\min}{g\in\G} R(g) \le \frac{q (1 \wedge 2^{q-2})B^q}{q-1} (\log 2) \, \frac{\log_2 |\G|}{n}
        \endarlab
The following corollary of Theorem \ref{th:appassouad} shows that this result is tight.

\begin{theorem} \label{th:qsup1}
Let $B>0$ and $d\in\N^*.$ For any training set size $n\in\N^*$ and 
input space $\X$ containing at least $\lfloor\log_2 (2d)\rfloor$ points, 
there exists a set $\G$ of $d$ prediction functions such that:
for any estimator $\hg$ there exists a probability distribution on the data space
$\X\times[-B;B]$ for which
        \begar
        \E R(\hg) - \und{\min}{g\in\G} R(g) 
                \ge \big( \frac{q}{90(q-1)} \vee e^{-1} \big)
                        B^q \big( \frac{\lfloor \log_2 |\G| \rfloor}{n+1} \wedge 1 \big).
        \endar
\end{theorem}

\begin{proof}
See Section \ref{sec:proofthlqpg}.
\end{proof}

\begin{remark} \label{rem:commop}
For least square regression (i.e. q=2), Remark \refp{rem:common} holds provided that  
the multiplicative factor becomes $2e\log 2 \approx 3.77$. 
More generally, the method
used here gives close to optimal constants but not the exact ones.
We believe that this limit is due to the use of the hypercube structure. Indeed,
the reader may check that for hypercubes of distributions, the upper bounds
used in this section are not constant-optimal since the simplifying step consisting in using 
$\min_{\rho\in\M} \cdots \le \min_{g\in\G} \cdots$ is loose.
\end{remark}
\end{itemize}

The reader may recover that there are essentially two classes of bounded losses: the ones which are not convex
or not enough convex (typical examples are the classification loss, the hinge loss and the absolute loss)
and the ones which are sufficiently convex (typical examples are the least square loss, the entropy loss, 
the logit loss and the exponential loss). For the first class of losses, the edge discrepancy of type I is 
proportional to $\sqrt{\db}$ for constant and symmetrical hypercubes and \eqref{eq:w1a} leads to a convergence rate 
of $\sqrt{\fracl{\log |\G|}{n}}$. For the second class, the convergence rate 
is $\fracl{\log |\G|}{n}$ and the lower bound can be explained by the fact that,
when two prediction functions are different on a set with low probability (typically $n^{-1}$),
it often happens that the training data has no input points in this set. For such training data,
it is impossible to consistently choose the right prediction function.

This picture of convergence rates for finite models is rather well-known,
since 
\begin{itemize}
\item similar bounds (with looser constants) were known before for some cases (e.g. in classification, see \cite{Vap74,Dev96}). 
\item mutatis mutandis, the picture exactly matches the picture in the individual sequence prediction literature:
for mixable loss functions (similar to ``sufficiently convex''), the minimax regret is 
$\gdo(\log |\G|)/n$, whereas for $0/1$-type loss functions, it is $\gdo\big(\sqrt{(\log |\G|)/n}\big)$ (see e.g. \cite{Hau98}). 
\end{itemize}


\paragraph{$L_q$-regression for unbounded outputs having finite moments.}   \label{sec:exlqu}

\begin{itemize}
\item {\it Case $q=1$ :} From \eqref{eq:uq1}, when ${\sup}_{g\in\G} E_Z g(X)^2\le b^2$ for some $b>0$,
there exists an estimator for which
        \begar
        \E R(\hg) - \undc{\min}{g\in\G} R(g) \le 2b\sqrt{\fracl{2\log |\G|}{n}}.
        \endar
        
The following corollary of Theorem \ref{th:appassouad} shows that this result is tight.

\begin{theorem} 
For any training set size $n\in\N^*$, positive integer $d$, positive real number $b$ and 
input space $\X$ containing at least $\lfloor\log_2 d\rfloor$ points, 
there exists a set $\G$ of $d$ prediction functions uniformly bounded by $b$ such that:
for any estimator $\hg$ there exists a probability distribution for which
$\E|Y| < +\infty$ and 
       \begar
       \E R(\hg) - \und{\min}{g\in\G} R(g) 
               \ge \frac{b}{4}\sqrt{\frac{\lfloor \log_2 |\G| \rfloor}{n} \wedge \frac{1}{4}}
       \endar
\end{theorem}

\begin{proof}
Let $\tm = \lfloor \log_2 |\G| \rfloor$.
We consider a $\big(\tm,1/\tm,\sqrt{\frac{\tm}{4n}\wedge 1}\big)$-hypercube
with $\ha\equiv -b$ and $\hb\equiv b$.
One may check that
$\da= b \sqrt{\db}$ so that \eqref{eq:w1a} gives that
for any estimator there exists a probability distribution for which
$\E|Y| < +\infty$ and 
       \begar
       \E R(\hg) - \und{\min}{g\in\G} R(g) 
               \ge b\sqrt{\frac{\tm}{4n} \wedge 1} \Big( 1 - \sqrt{\frac{1}{4} \wedge \frac{n}{\tm}} \Big),
       \endar
hence the desired result.
\end{proof}

\item {\it Case $q>1$ :} 
First let us recall the upper bound. In Corollary \ref{cor:genp}, 
under the assumptions
        \lbegar
        {\sup}_{g\in\G,x\in\X} |g(x)| \le b \qquad \text{for some } b>0\\
        \E |Y|^s\le A \qquad \text{for some } s\ge q \text{ and } A>0\\
        \G \text{ finite}
        \rendar
we have proposed an algorithm satisfying        
        \begar
        R(\hg) - \und{\min}{g\in\G} R(g) \le \left\{ \begin{array}{lll}
        C \big(\frac{\log |\G|}{n}\big)^{1-\frac{q-1}{s}} & \qquad \text{ when }q \le s \le 2q-2\\
        C \big(\frac{\log |\G|}{n}\big)^{1-\frac{q}{s+2}} & \qquad \text{ when }s\ge 2q-2\\
        \end{array} \right..
        \endar
for a quantity $C$ which depends only on $b$, $A$, $q$ and $s$.

The following corollary of Theorem \ref{th:appassouad} shows that this result is tight
and is illustrated by Figure \ref{fig:lqunbounded}.

\begin{theorem} \label{th:lqunbounded}
Let $d\in\N^*$, $s\ge q>1$, $b>0$ and $A>0$.
For any training set size $n\in\N^*$ and input space $\X$ containing at least $\lfloor\log_2 (2d)\rfloor$ points, 
there exists a set $\G$ of $d$ prediction functions uniformly bounded by $b$ such that:
for any estimator $\hg$ there exists a probability distribution on the data space $\X\times \R$ for which
$\E|Y|^s \le A$ and 
        \begar
        \E R(\hg) - \und{\min}{g\in\G} R(g) 
                \ge \left\{ \begin{array}{lll}
        C \big(\frac{\log |\G|}{n} \wedge 1\big)^{1-\frac{q-1}{s}}\\ 
        C \big(\frac{\log |\G|}{n} \wedge 1\big)^{1-\frac{q}{s+2}} 
        \end{array} \right.,
        \endar
for a quantity $C$ which depends only on the real numbers $b$, $A$, $q$ and $s$.
\end{theorem}

Both inequalities simultaneously hold but the first one is tight for $q \le s \le 2q-2$
while the second one is tight for $s\ge 2q-2$.
They are both based on \eqref{eq:w1a} applied to a $\lfloor \log_2 |\G| \rfloor$-dimensional
hypercubes.

Contrary to other lower bounds obtained in this work, the first inequality is based on asymmetrical hypercubes. 
The use of this kind of hypercubes can be partially explained by the fact that the learning task
is asymmetrical. Indeed all values of the output space do not have the same status since
predictions are constrained to be in $[-b;b]$ while outputs are allowed to be in the whole real space
(see the constraints on the hypercube in the proof given in Section \ref{sec:proofth:lqunbounded}). 

\end{itemize}

\begin{figure}[ht] 
\begin{minipage}[t]{7cm}
\begin{center}
\includegraphics[width=7cm,clip]{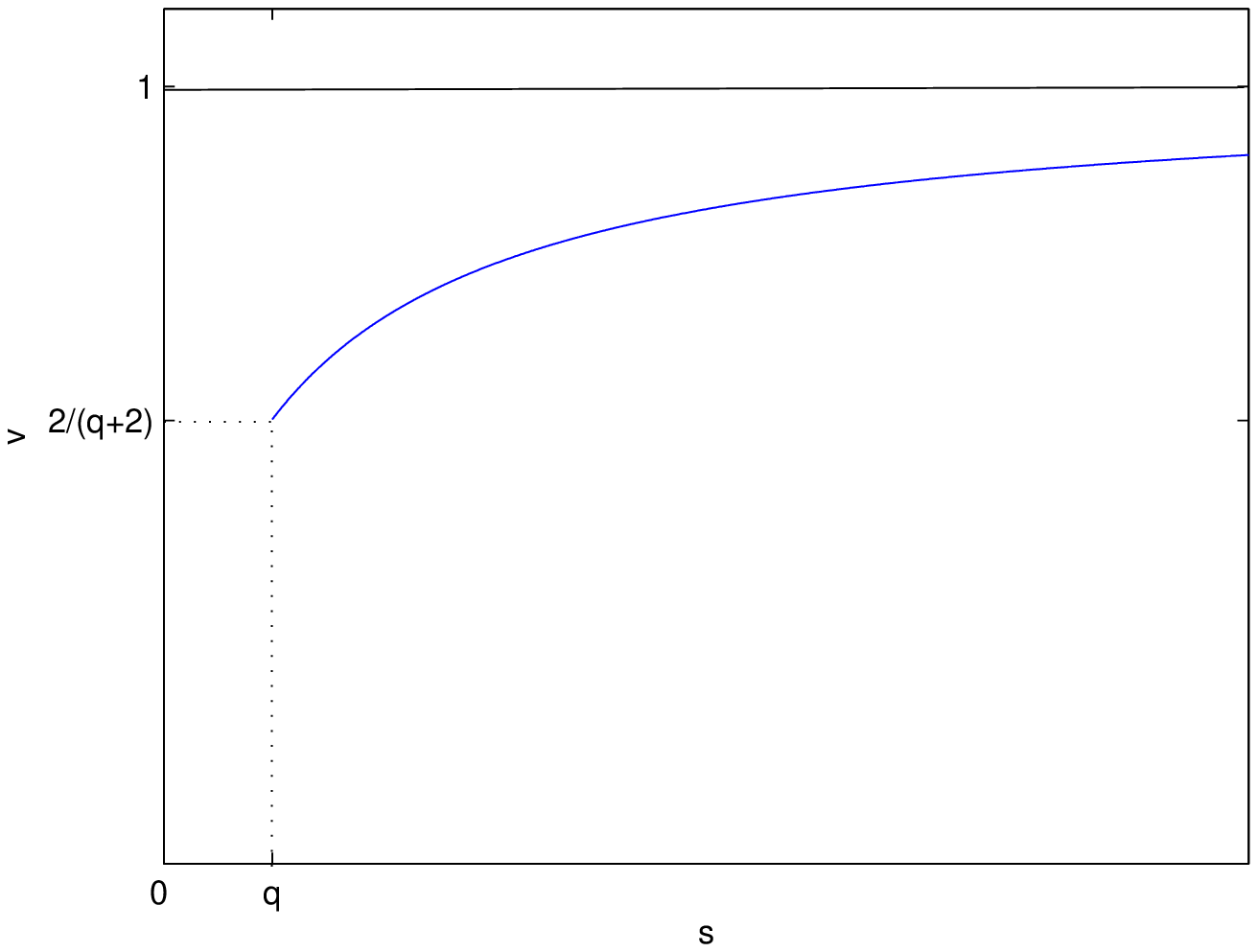}
\end{center}
\end{minipage}
\hfill
\begin{minipage}[t]{7cm}
\begin{center}
\includegraphics[width=7cm,clip]{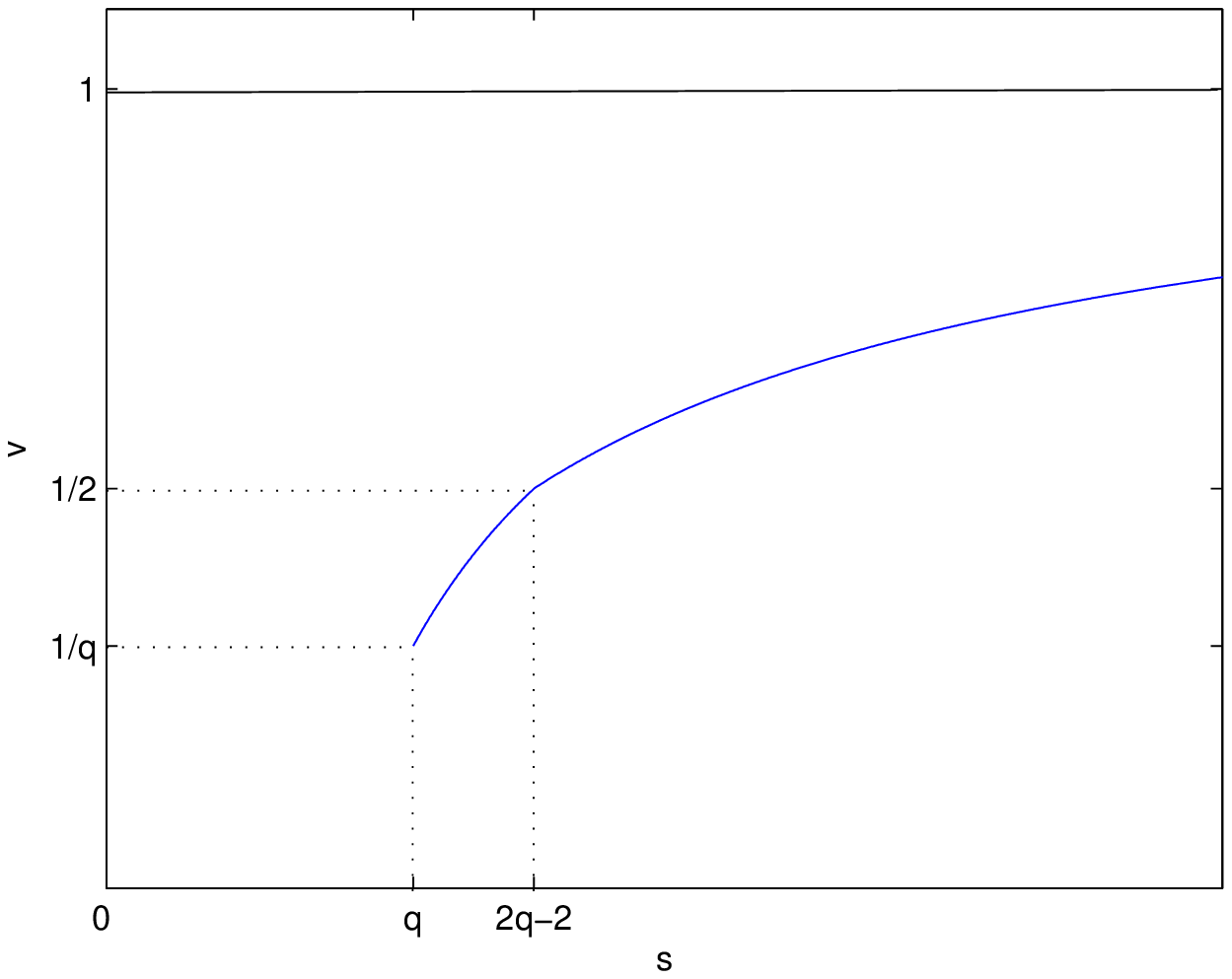}
\end{center}
\end{minipage}
\caption{{\it Optimal convergence rates in $L_q$-regression when the output has a finite moment of order $s$ 
(see Theorem \ref{th:lqunbounded}).}
The convergence rate is of order 
        $\big(\frac{log|\G|}{n}\big)^{v}$ with $0 < v \le 1$.
The figure represents the value of $s$ in abscissa and the value of $v$ in ordinate.
Two cases have to be distinguished. For $1<q\le 2$ (figure on the left), $v$ depends smoothly on $q$.
For $q>2$ (figure on the right), two stages are observed depending whether $s$ is larger than $2q-2$.
} \label{fig:lqunbounded}
\end{figure}

\subsubsectionjy{Numerical comparison of the lower bounds \eqref{eq:w1a}, \eqref{eq:w1b} and \eqref{eq:w2}} \label{sec:numerical}

To compare \eqref{eq:w1a}, \eqref{eq:w1b} and \eqref{eq:w2},
we will compare their asymptotical version, i.e. \eqref{eq:asw1a}, \eqref{eq:asw1b} and \eqref{eq:asw2}. 

\paragraph{Classification in VC classes.} \label{sec:exvc}

We consider $|\Y|=2$ and $\ell(y,y')=\ds1_{y\neq y'}$. 
Since the differentiability assumption does not hold in this setting, 
we only compare \eqref{eq:w1a} and \eqref{eq:w1b}.

\begin{theorem} \label{th:asvc}
Let $\calP$ be the set of all probability distributions on the data space~$\Z$.
Let $(\G_n)_{n\in\N}$ be a family of prediction function spaces of VC-dimension $V_n \ge 2$
satisfying $n / V_n \und{\rightarrow}{n\rightarrow +\infty} +\infty$.
For any algorithm~$\hg$:
        \begar
        \und{\liminf}{n \rightarrow +\infty} \; \sqrt{\frac{n}{V_n}} \;
                \underset{P \in \calP}{\sup} \big\{ \Ezun R( \hg_{Z_1^n} ) 
                - \und{\inf}{g\in\G_n} R(g) \big\} 
                \ge \left\{ \begin{array}{lll}
        \alpha_1 & \text{ from \eqref{eq:asw1a}}\\
        \alpha_2 & \text{ from \eqref{eq:asw1b}}
        \end{array} \right.
        \endar
        $
        \text{with }\left\{ \begin{array}{lll}
        \alpha_1 & = & \und{\max}{a>0} \frac{\sqrt{a} (1-\sqrt{1-e^{-a}}) }{2} \approx 0.135 \\ 
        \alpha_2 & = & \und{\max}{a>0} \sqrt{\frac{a}{2\pi}} \int_{a}^{+\infty} e^{-t^2/2} dt \approx 0.170 
        \end{array} \right..
        $\\
In particular for a given set $\G$ of finite VC-dimension $V$, for $n$ sufficiently large, any estimator $\hg$ satisfies
        \begar
        \underset{P \in \calP}{\sup} \big\{ \E R(\hg) - \undc{\inf}{g\in\G} R(g) \big\}
                \ge \frac{1}{6}\sqrt{V/n}.
        \endar
\end{theorem}
\begin{proof}
It suffices to apply Corollary \ref{cor:appassouad} to $(V_n,1/V_n, a V_n/n)$-hypercubes, use that 
$\da=\sqrt{\db}/2$ and choose the real number $a>0$ to maximize the lower bound. 
\end{proof}

The two inequalities, coming from \eqref{eq:asw1a} and \eqref{eq:asw1b},
simultaneously hold. They only differ by a multiplicative constant.

\paragraph{Least square regression with unbounded outputs satisfying $\E Y^2 \le A$ for some $A>0$.}

We consider the context of Corollary \ref{cor:gen2} with $s=2$.
The best explicit constant in \eqref{eq:squ} is obtained from \eqref{eq:gen2}
for  $\lam=\sqrt{\frac{\log |\G|}{8b^2 A (n+1)}}\wedge (8b^2)^{-1}$.
When $\log |\G| \le An/(8b^2)$, we get
        \begar
        R(\hg) - \und{\min}{g\in\G} R(g) \le \sqrt{32} \, b \sqrt{A}\sqrt{\frac{\log |\G|}{n}}
        \endar
Now let us give the associated lower bounds coming 
from \eqref{eq:asw1a}, \eqref{eq:asw1b} and \eqref{eq:asw2}.

\begin{theorem} \label{th:l2ub}
Let $\X$ be an infinite input space.
Let $A>0$ and $\calP$ be the set of probability distributions 
on $\X\times \R$ such that $\E Y^2 \le A$.
There exists a family of prediction function spaces $(\G_n)_{n\in\N}$ such that
any prediction function in these sets is uniformly bounded by $b$,
their sizes at most grows subexponentially, i.e. $n/\log |\G_n|$ goes to infinity
when $n$ goes to infinity, and for any algorithm $\hg$ 
        \begar
        \und{\liminf}{n \rightarrow +\infty} \; \sqrt\frac{n}{\log |\G_n|} \;
                \underset{P \in \calP}{\sup} \big\{ \Ezun R( \hg_{Z_1^n} ) 
                - \und{\min}{g} R(g) \big\} 
                \ge \left\{ \begin{array}{lll}
        \beta_1 \, b \sqrt{A} \text{  \, from \eqref{eq:asw1a}}\\
        \beta_2 \, b \sqrt{A} \text{  \, from \eqref{eq:asw1b}}\\
        \beta_3 \, b \sqrt{A} \text{  \, from \eqref{eq:asw2}}  
        \end{array} \right.
        \endar
with 
        \lbegar
        \beta_1 & = & (\log 2)^{-1} \max_{a>0} \sqrt{a} (1-\sqrt{1-e^{-a}}) 
                \approx 0.3897\\
        \beta_2 & = & (\log 2)^{-1} \max_{a>0} \sqrt{\frac{2a}{\pi}} \int_{a}^{+\infty} e^{-t^2/2} dt
                \approx 0.4904\\
        \beta_3 & = & (\log 2)^{-1} \max_{a>0} \sqrt{a} \big(1+\demi e^{-a/2}-\demi e^{3a/2}\big) 
                \approx 0.5154\\
        \rendarp
\end{theorem}

\jyproof{th:l2ub}

The three inequalities, coming from \eqref{eq:asw1a}, \eqref{eq:asw1b} and \eqref{eq:asw2},
only differ by a multiplicative constant. The one coming from \eqref{eq:asw2} gives the tightest 
result. The difference between the upper bound and this lower bound is a multiplicative factor smaller than $11$.

\begin{remark}
In all these examples (i.e. the ones of Section \ref{sec:examples}), we have only considered constant hypercubes. The use of non-constant
hypercubes can be required when smoothness assumptions are put on the regression function
$\eta:x\mapsto P(Y=1|X=x)$. This is typically the case in works on plug-in classifiers (\cite{Ant02,AudTsy05}).
For instance, the proof of \cite[Theorems 3.5 and 4.1]{AudTsy05} relies on 
non-constant symmetrical hypercubes for which the function $\xi$ (see Definition \ref{def:hyp2}) is chosen 
such that it vanishes on the border of the partition cells, which ensures the regularity of the regression function $\eta$.
\end{remark}

\section{Summary of contributions and open problems} \label{sec:sum}

This work has developed minimax optimal risk bounds for the 
general learning task consisting in predicting as well as the
best function in a reference set.
It has proposed to summarize this learning problem by the variance function
appearing in the variance inequality (p.\pageref{varcond}).
The SeqRand algorithm (Figure \ref{fig:seqr}) based on this variance function
leads to minimax optimal convergence rates in the model selection aggregation problem,
and our analysis gives a nice unified view to results coming from different communities.

In particular, results coming from the online learning literature are recovered in Section \ref{sec:onl2bat}.
Corollary \ref{cor:o} gives a new bound in the online learning setting (sequential prediction with expert advice).
The generalization error bounds obtained by Juditsky, Rigollet and Tsybakov in \cite{Jud06}
are recovered for a slightly different algorithm in Section \ref{sec:jud}.

Without any extra assumption on the learning task, we have obtained a Bernstein's type 
bound which has no known equivalent form when the loss function is not assumed to be bounded (Section \ref{sec:berns}).
When the loss function is bounded, the use of Hoeffding's inequality $\wrt$ 
Gibbs distributions on the prediction function space 
instead of the distribution generating the data leads to an improvement by a factor $2$
of the standard-style risk bound (Theorem \ref{th:hoeff}).

To prove that our bounds are minimax optimal, we have refined Assouad's lemma particularly by taking into account 
the properties of the loss function. 
Theorem \ref{th:appassouad} is tighter than previous versions of Assouad's lemma and
easier to apply to a learning setting than Fano's lemma (see e.g. \cite{Tsy04c}), 
besides the latter leads in general to very loose constants.
It improves the constants of lower bounds related to Vapnik-Cervonenkis classes by a factor greater than $1000$.
We have also illustrated our upper and lower bounds by studying 
the influence of the noise of the output and of the convexity of the loss function.

For the $L_q$-loss with $q \ge 1$, new matching upper and lower bounds are given:
in the online learning framework under boundedness assumption
(Corollary \ref{cor:seqlq} and Section \ref{sec:exlqb} jointly with Remark \ref{rem:onlow}), 
in the batch learning setting under boundedness assumption (Section \ref{sec:onl2bat} and Section \ref{sec:exlqb}), 
in the batch learning setting for unbounded observations under moment assumptions (Section \ref{sec:pow} and \ref{sec:exlqu}).
In the latter setting, we still do assume that the prediction functions are bounded. It is an open problem to replace this 
boundedness assumption with a moment condition.

Finally this work has the following limits.
Most of our results concern expected risks and it is an open problem to provide 
corresponding tight exponential inequalities. 
Besides we should emphasize that our expected risk upper bounds
hold only for our algorithm. This is quite different 
from the classical point of view that simultaneously gives 
upper bounds on the risk of any prediction function in the model. 
To our current knowledge, this classical approach has a flexibility
that is not recovered in our approach.
For instance, in several learning tasks, Dudley's chaining trick \cite{Dud78}
is the only way to prove risk convergence with the optimal rate. So a natural
question and another open problem is whether it is possible to combine the better 
variance control presented here with the chaining argument (or other localization argument
used while exponential inequalities are available).

\fi
\section{Proofs} \label{sec:proofs}
\ifRR
\secproofth{kivwar}

First, by a scaling argument, it suffices to prove the result for $a=0$ and $b=1$.
For $\Y=[0;1]$, we modify the proof in Appendix A of \cite{Kiv99L}. Precisely, claims $1$ and $2$, with 
the notation used there, become:
\begin{enumerate}
\item 
If the function $f$ is concave in $\alpha([p;q])$ then we have $A_t(q)\le B_t(p)$,
\item If $c\ge R(z,p,q)$ for any $z\in(p;q)$, then the function $f$ is 
concave in $\alpha([p;q])$.
\end{enumerate}
Up to the missing $\alpha$ (typo), the difference is that we restrict ourselves to values of 
$z$ in $[p;q]$.
The proof of Claim $2$ has no new argument. 
For claim 1, it suffices to modify the definition
of $x'_{t,i}$ into $x'_{t,i}=q \wedge G^{-1}[\ell(p,x_{t,i})]\in[p;q]$. Then we have
$L(p,x'_{t,i})\le L(p,x_{t,i})$ and $L(q,x'_{t,i})\le L(p,x_{t,i})$, hence
$\alpha(x'_{t,i})\ge \alpha(x_{t,i})$ and $\gamma(x'_{t,i})\ge \gamma(x_{t,i})$.
Now one can prove that $f$ is decreasing on $\alpha([p;q])$.
By using Jensen's inequality, we get
        \begar
        \Delta_{t}(q) & = & - c \log \sum_{i=1}^n v_{t,i}\gamma(x_{t,i})\\
        & \ge & - c \log \sum_{i=1}^n v_{t,i}\gamma(x'_{t,i})\\
        & = & - c \log \sum_{i=1}^n v_{t,i}f[\alpha(x'_{t,i})]\\
        & \ge & - c \log f\big[\sum_{i=1}^n v_{t,i}\alpha(x'_{t,i})\big]\\
        & \ge & - c \log f\big[\sum_{i=1}^n v_{t,i}\alpha(x_{t,i})\big]\\
        & = & L[q,G^{-1}(\Delta_t(p))]\\
        \endar
The end of the proof of claim $1$ is then identical.

\subsection{Proof of Corollary \ref{cor:jud}} \label{sec:corjud}

We start by proving that the variance inequality holds with
$\delta_\lam \equiv 0$, and that we may take 
$\pirho$ be the Dirac distribution at the function
        $\expec{g}{\rho} g$.
By using Jensen's inequality and Fubini's theorem, Assumption \eqref{eq:jud} implies that
        \begar
        \expec{g'}{\pirho} \undc{\E}{Z\sim P} \log \expec{g}{\rho} 
                e^{\lam [L(Z,g') - L(Z,g)]}\\
        \qquad\qquad\qquad\qquad\qquad
                = \; \undc{\E}{Z\sim P} \log \expec{g}{\rho} 
                e^{\lam [L(Z,\expec{g'}{\rho} g') - L(Z,g)]}\\
        \qquad\qquad\qquad\qquad\qquad
                \le \; \log \expec{g}{\rho} \undc{\E}{Z\sim P} 
                e^{\lam [L(Z,\expec{g'}{\rho} g') - L(Z,g)]}\\
        \qquad\qquad\qquad\qquad\qquad
                \le \; \log \expec{g}{\rho} \psi(\expec{g'}{\rho} g',g)\\
        \qquad\qquad\qquad\qquad\qquad
                \le \; \log \psi(\expec{g'}{\rho} g',\expec{g}{\rho} g)\\
        \qquad\qquad\qquad\qquad\qquad = \; 0,
        \endar
so that we can apply Theorem \ref{th:1}. It remains to note that in this 
context the SeqRand algorithm is the one described in the corollary.

\subsection{Proof of Theorem \ref{cor:var}} \label{sec:corvar}

To check that the variance inequality holds, it suffices to prove that for any $z\in\Z$
        \begarlab{eq:before}
        \expec{g'}{\rho} \log \expec{g}{\rho} e^{\lam[L(z,g')-L(z,g)]-\frac{\lam^2}{2}[L(z,g')-L(z,g)]^2}
                \le 0.
        \endarlab
To shorten formulae, let $\alpha(g',g) \eqdef \lam[L(z,g')-L(z,g)]$.
By Jensen's inequality and the following symmetrization trick, \eqref{eq:before} holds.
        \begarlab{eq:sym}
        \expec{g'}{\rho} \expec{g}{\rho} e^{\alpha(g',g) - \frac{\alpha^2(g',g)}{2}}\\
        \quad
                \le \frac{1}{2} \expec{g'}{\rho} \expec{g}{\rho} e^{\alpha(g',g) - \frac{\alpha^2(g',g)}{2}}
                + \frac{1}{2} \expec{g'}{\rho} \expec{g}{\rho} e^{-\alpha(g',g) - \frac{\alpha^2(g',g)}{2}}\\
        \quad
                \le \expec{g'}{\rho} \expec{g}{\rho} \cosh\big(\alpha(g,g')\big) e^{- \frac{\alpha^2(g',g)}{2}}\\
        \quad
                \le 1
        \endarlab
where in the last inequality we used the inequality $\cosh(t)\le e^{t^2/2}$ for any $t\in\R$.
The result then follows from Theorem \ref{th:1}.

\subsection{Proof of Corollary \ref{cor:mamtsy}} \label{sec:cormamtsy}

To shorten the following formula, let $\mu$ denote the law of the prediction function
produced by the SeqRand algorithm ($\wrt$ simultaneously the training set and the randomizing procedure).
Then \eqref{eq:corvar} can be written as: for any $\rho\in\M$,
        \begarlab{eq:pr1}
        \expec{g'}{\mu} R(g')
                \le 
                \expec{g}{\rho} R(g) 
                + \frac{\lam}{2} \expec{g}{\rho} \expec{g'}{\mu} V(g,g') 
                + \frac{K(\rho,\pi)}{\lam(n+1)}
        \endarlab
Define $\tildR(g) = R(g)-R(\tildg)$ for any $g\in\G$.
Under the generalized Mammen and Tsybakov assumption, for any $g,g'\in\G$, we have
        \begar
        \frac{1}{2}V(g,g') & \le & \undc{\E}{Z} \big\{[L(Z,g)-L(Z,\tildg)]^2\big\}+\undc{\E}{Z} \big\{[L(Z,g')-L(Z,\tildg)]^2\big\}\\
        & \le & c \tildR^{\gamma}(g)+c \tildR^{\gamma}(g'),
        \endar
so that \eqref{eq:pr1} leads to 
        \begarlab{eq:pr2}
        \expec{g'}{\mu} [ \tildR(g') - c \lam \tildR^\gamma(g') ]
                \le 
                \expec{g}{\rho} [ \tildR(g) 
                + c \lam \tildR^\gamma(g) ] 
                + \frac{K(\rho,\pi)}{\lam(n+1)}.
        \endarlab 
This gives the first assertion. For the second statement, let ${\tilde u}\eqdef \expec{g'}{\mu} \tildR(g')$ 
and $\chi(u)\eqdef u-c \lam u^\gamma$. By Jensen's inequality, the $\lhs$ of \eqref{eq:pr2} is lower bounded
by $\chi({\tilde u})$. By straightforward computations, for any $0<\beta<1$, when 
$u\ge \big( \frac{c\lam}{1-\beta} \big)^{\frac{1}{1-\gamma}}$,
$\chi(u)$ is lower bounded by $\beta u$, which implies the desired result.

\subsection{Proof of Theorem \ref{cor:pac}} \label{sec:corpac}

Let us prove \eqref{eq:pacexp1}.
Let $r(g)$ denote the empirical risk of $g\in\G$, that is
        $r(g)=\frac{\Sigma_n(g)}{n}.$
Let $\rho\in\M$ be some fixed distribution on $\G$. 
From \cite[Section 8.1]{Aud03b}, with probability at least $1-\eps$ $\wrt$
the training set distribution, for any $\mu \in \M$, we have
        \begar
        \expec{g'}{\mu} R(g') - \expec{g}{\rho} R(g)\\
        \;
                \le \expec{g'}{\mu} r(g') - \expec{g}{\rho} r(g)
                +\lam \varphi(\lam B) \expec{g'}{\mu} \expec{g}{\rho} V(g,g')
                +\frac{K(\mu,\pi)+\logeps}{\lam n}.
        \endar
Since the Gibbs distribution $\pi_{-\lam \Sigma_n}$  minimizes $\mu\mapsto\expec{g'}{\mu} r(g') +\frac{K(\mu,\pi)}{\lam n}$,
we have 
        \begar 
        \expec{g'}{\pi_{-\lam \Sigma_n}} R(g')\\
        \quad\qquad
                \le \expec{g}{\rho} R(g) 
                +\lam \varphi(\lam B) \expec{g'}{\pi_{-\lam \Sigma_n}} \expec{g}{\rho} V(g,g')
                +\frac{K(\rho,\pi)+\logeps}{\lam n}.
        \endar
Then we apply the following inequality 
        \begar
        \E W \le \E (W \vee 0) = \int_{0}^{+\infty} \P(W > u) du =\int_{0}^{1} \eps^{-1} \P(W > \logeps) d\eps
        \endar
to the random variable 
        \begar
        W=\lam n \big[ \expec{g'}{\pi_{-\lam \Sigma_n}} R(g')
                - \expec{g}{\rho} R(g) \\
        \qquad\qquad\qquad\qquad
                - \lam \varphi(\lam B) \expec{g'}{\pi_{-\lam \Sigma_n}} \expec{g}{\rho} V(g,g') \big]
                - K(\rho,\pi).
        \endar
We get $\E W \le 1$.
At last we may choose the distribution $\rho$ minimizing the upper bound to obtain \eqref{eq:pacexp1}.
Similarly using \cite[Section 8.3]{Aud03b}, we may prove \eqref{eq:pacexp2}.

\subsection{Proof of Lemma \ref{lem:conczha}} \label{sec:proofconczha}

It suffices to apply the following adaptation of 
Lemma 5 of \cite{Zha05} to \begar\xi_i(Z_1,\dots,Z_i)=L[Z_i,\A(Z_1^{i-1})]-L(Z_i,\tildg).\endar

\begin{lemma} \label{lem:zha}
Let $\varphi$ still denote the positive convex increasing function defined as
$\varphi(t) \eqdef \frac{e^t-1-t}{t^2}$. Let $b$ be a real number.
For $i=1,\dots,n+1$, let $\xi_i:\Z^i \rightarrow \R$ be a function uniformly upper bounded by $b$.
For any $\eta>0$, $\eps>0$, with probability at least $1-\eps$ $\wrt$ the distribution of $Z_1,\dots,Z_{n+1}$, we have
        \begarlab{eq:zha}
        \sum_{i=1}^{n+1} \xi_i(Z_1,\dots,Z_i) \le \sum_{i=1}^{n+1} \E_{Z_i}\xi_i(Z_1,\dots,Z_i) \\
        \qquad\qquad\qquad\qquad\quad   
                + \eta \varphi(\eta b) \sum_{i=1}^{n+1} \E_{Z_i} \xi_i^2(Z_1,\dots,Z_i) + \frac{\logeps}{\eta},
        \endarlab
where $\E_{Z_i}$ denotes the expectation $\wrt$ the distribution of $Z_i$ only.
\end{lemma}

\begin{remark}
The same type of bounds without variance control can be found in~\cite{Ces04}.
\end{remark}

\begin{proof}
For any $i\in\{0,\dots,n+1\}$, define
        \begar
        \psi_i=\psi_i(Z_1,\dots,Z_i)\eqdef\sum_{j=1}^i \xi_j - \sum_{j=1}^i \E_{Z_j}\xi_j - \eta \varphi(\eta b) 
                \sum_{j=1}^i \E_{Z_j}\xi_j^2. 
        \endar
where $\xi_j$ is the short version of $\xi_j(Z_1,\dots,Z_j)$.
For any $i\in\{0,\dots,n\}$, we trivially have 
        \begarlab{eq:psis}
        \psi_{i+1}-\psi_{i}=\xi_{i+1} - \E_{Z_{i+1}} \xi_{i+1} - \eta \varphi(\eta b) \E_{Z_{i+1}} \xi_{i+1}^2.
        \endarlab
Now for any $b\in\R$, $\eta>0$ and any random variable $W$ such that $W\le b$ a.s., we have
        \begarlab{eq:jy71}
        \E e^{\eta(W - \E W -\eta \varphi(\eta b) \E W^2)} \le 1.
        \endarlab
\begin{remark}
The proof of \eqref{eq:jy71} is standard and can be found e.g. in \cite[Section~7.1.1]{Aud03a}.
We use \eqref{eq:jy71} instead of the inequality used to prove Lemma 5 of \cite{Zha05},
i.e. $\E e^{\eta[W - \E W -\eta \varphi(\eta b') \E( W -\E W)^2]} \le 1$ for
$W -\E W\le b'$ since we are interested in excess risk bounds. Precisely,
we will take $W$ of the form 
        $W=L(Z,g)-L(Z,g')$ 
for fixed functions $g$ and $g'$. Then we have
$W \le \sup_{z,g} L - \inf_{z,g} L$
while we only have 
$W-\E W \le 2\big( \sup_{z,g} L - \inf_{z,g} L\big)$. Besides the gain of having 
$\E( W -\E W)^2$ instead of $\E W^2$ is useless in the applications we develop here.
\end{remark}

By combining \eqref{eq:jy71} and \eqref{eq:psis}, we obtain
        \begarlab{eq:epsis}
        \E_{Z_{i+1}} e^{\eta(\psi_{i+1}-\psi_{i})}\le 1.
        \endarlab
By using Markov's inequality, we upper bound the following probability
$\wrt$ the distribution of $Z_1,\dots,Z_{n+1}$:
        \begar
        \P\Big( \sum_{i=1}^{n+1} \xi_i
                > \sum_{i=1}^{n+1} \E_{Z_i} \xi_i 
                + \eta \varphi(\eta b) \sum_{i=1}^{n+1} \E_{Z_i} \xi_i^2 + \frac{\logeps}{\eta} \Big)\\
        \qquad\qquad
                = 
        \P\big( \eta \psi_{n+1} > \logeps \big)\\
        \qquad\qquad
                = 
        \P\big( \eps e^{\eta \psi_{n+1}} > 1 \big)\\
        \qquad\qquad
                \le
        \eps \E e^{\eta \psi_{n+1}} \\
        \qquad\qquad
                \le
        \eps \E_{Z_{1}} \big( e^{\eta (\psi_{1}-\psi_{0})} \E_{Z_{2}} \big( \cdots 
                e^{\eta (\psi_{n}-\psi_{n-1})} \E_{Z_{n+1}} e^{\eta (\psi_{n+1}-\psi_{n})} \big)\big)\\
        \qquad\qquad
                \le \eps
        \endar  
where the last inequality follows from recursive use of \eqref{eq:epsis}.
\end{proof}

\subsection{Proof of Theorem \ref{th:genpow}}  \label{sec:proofgenpow}

The first inequality follows from Jensen's inequality. Let us prove the second.
According to Theorem \ref{th:1}, it suffices to check that the variance inequality holds for 
$0< \lam \le \lam_0$, $\pirho$ the Dirac distribution at $\expec{g}{\rho} g$ and 
        \begar
        \delta_{\lam}[(x,y),g,g'] = \delta_{\lam}(y) \eqdef 
                \und{\min}{0\le \zeta \le 1} 
                \Big[ \zeta \Delta(y) + \frac{(1-\zeta)^2\lam \Delta^2(y)}{2} \Big] \ds1_{|y|> B}\\
        \qquad\qquad\qquad
                = \frac{\lam \Delta^2(y)}{2}  
                \ds1_{\lam \Delta(y)<1;|y|> B} + \big[\Delta(y)-\frac{1}{2\lam}\big] \ds1_{\lam \Delta(y)\ge1;|y|> B}.
        \endar

\begitem
\item
For any $z=(x,y)\in\Z$ such that $|y|\le B$, for any probability distribution $\rho$ and for the above values of 
$\lam$ and $\delta_{\lam}$, by Jensen's inequality, we have
        \begar
        \expec{g}{\rho} e^{\lam[L(z,\expec{g'}{\rho} g') - L(z,g) - \delta_{\lam}(z,g,g')]}\\
        \qquad\qquad\qquad
                = e^{\lam L(z,\expec{g'}{\rho} g')} \expec{g}{\rho}  e^{- \lam \ell[y,g(x)]} \\
        \qquad\qquad\qquad
                \le e^{\lam L(z,\expec{g'}{\rho} g')} \Big( \expec{g}{\rho}  e^{- \lam_0 \ell[y,g(x)]} \Big)^{\lam/\lam_0}\\
        \qquad\qquad\qquad
                \le e^{\lam \ell[y,\expec{g'}{\rho} g'(x)] - \lam \ell[y,\expec{g}{\rho}g(x)]} \\
        \qquad\qquad\qquad
                = 1,
        \endar
where the last inequality comes from the concavity of $y' \mapsto e^{-\lam_0 \ell(y,y')}$. 
This concavity argument goes back to \cite[Section 4]{Kiv99}, and was also used in 
\cite{Bun05} and in some of the examples given in \cite{Jud06}.
\item
For any $z=(x,y)\in\Z$ such that $|y|> B$, for any $0 \le \zeta \le 1$, by using twice Jensen's inequality 
and then by using the symmetrization trick presented in Section \ref{sec:var}, we have
        \begar
        \expec{g}{\rho} e^{\lam[L(z,\expec{g'}{\rho} g') - L(z,g) - \delta_{\lam}(z,g,g')]}\\
        \qquad
                = e^{- \delta_{\lam}(y)} \expec{g}{\rho} e^{\lam[L(z,\expec{g'}{\rho} g') - L(z,g)]}\\
        \qquad
                \le e^{- \delta_{\lam}(y)} \expec{g}{\rho} e^{\lam[\expec{g'}{\rho} L(z,g') - L(z,g)]}\\
        \qquad
                \le e^{- \delta_{\lam}(y)} \expec{g}{\rho} \expec{g'}{\rho} e^{\lam[L(z,g') - L(z,g)]}\\
        \qquad
                = e^{- \delta_{\lam}(y)} \expec{g}{\rho} \expec{g'}{\rho} \Big\{ e^{\lam(1-\zeta)[L(z,g') - L(z,g)]
                - \frac{1}{2} \lam^2(1-\zeta)^2[L(z,g') - L(z,g)]^2}\\
        \qquad\qquad\qquad\qquad\qquad\qquad 
                \times \; e^{
                \lam \zeta [L(z,g') - L(z,g)]+ \frac{1}{2} \lam^2(1-\zeta)^2[L(z,g') - L(z,g)]^2}\Big\}\\
        \qquad
                \le e^{- \delta_{\lam}(y)} \expec{g}{\rho} \expec{g'}{\rho} \Big\{ e^{\lam(1-\zeta)[L(z,g') - L(z,g)]
                - \frac{1}{2} \lam^2(1-\zeta)^2[L(z,g') - L(z,g)]^2}\\
        \qquad\qquad\qquad\qquad\qquad\qquad 
                \times \; e^{\lam \zeta \Delta(y)
                + \frac{1}{2} \lam^2(1-\zeta)^2\Delta^2(y)}\Big\}\\
        \qquad
                \le e^{- \delta_{\lam}(y)} e^{ \lam \zeta \Delta(y)
                + \frac{1}{2} \lam^2(1-\zeta)^2\Delta^2(y)}\\
        \endar
Taking $\zeta\in[0;1]$ minimizing the last $\rhs$, we obtain that
        \begar
        \expec{g}{\rho} e^{\lam[L(z,\expec{g'}{\rho} g') - L(z,g) - \delta_{\lam}(z,g,g')]} \le 1
        \endar
\enditem
From the two previous computations, we obtain that for any $z\in\Z$,
        \begar
        \log \expec{g}{\rho} e^{\lam[L(z,\expec{g'}{\rho} g') - L(z,g) - \delta_{\lam}(z,g,g')]} \le 0,
        \endar
so that the variance inequality holds for the above values of $\lam$, $\pirho$ and $\delta_{\lam}$,
and the result follows from Theorem \ref{th:1}.

\subsection{Proof of Corollary \ref{cor:genp}} \label{sec:corgenp}

To apply Theorem \ref{th:genpow}, we will first determine $\lam_0$ for which
the function $\zeta: y' \mapsto e^{-\lam_0 |y-y'|^q}$ is concave.
For any given $y\in[-B;B]$, for any $q>1$, straightforward computations give
        \begar
        \zeta''(y') = \big[\lam_0 q |y'-y|^{q} - (q-1) \big] \lam_0 q |y'-y|^{q-2} e^{-\lam_0 |y-y'|^q}
        \endar
for $y'\neq y$, hence $\zeta'' \le 0$ on $[-b;b]-\{y\}$ for $\lam_0=\frac{q-1}{q(B+b)^q}.$
Now since the derivative $\zeta'$ is defined at the point $y$, we conclude that 
the function $\zeta$ is concave on $[-b;b]$, so that we may use Theorem \ref{th:genpow}
with $\lam_0=\frac{q-1}{q(B+b)^q}.$

Contrary to the least square setting, we do not have a simple close formula for $\Delta(y)$, but 
for any $|y|\ge b$, we have 
        \begar
        2bq(|y|-b)^{q-1} \le \Delta(y) \le 2bq (|y|+b)^{q-1}.
        \endar
As a consequence, when $|y| \ge b+(2bq\lam)^{-1/(q-1)}$, we have $\lam \Delta(y) \ge 1$ and 
$\Delta(y) - 1/(2\lam)$ can be upper bounded by $C' |y|^{q-1}$,
where the quantity $C'$ depends only on $b$ and $q$.

For other values of $|y|$, i.e. when $b \le |y| < b+(2bq\lam)^{-1/(q-1)}$, we have
        \begar
        \frac{\lam \Delta^2(y)}{2}  
                \ds1_{\lam \Delta(y)<1;|y|> B} + \big[\Delta(y)-\frac{1}{2\lam}\big] \ds1_{\lam \Delta(y)\ge1;|y|> B} \\
        \qquad\qquad\qquad\qquad
                = \und{\min}{0\le \zeta \le 1} 
                \Big[ \zeta \Delta(y) + \frac{(1-\zeta)^2\lam \Delta^2(y)}{2} \Big] \ds1_{|y|> B}\\
        \qquad\qquad\qquad\qquad
                \le 
                \demi \lam \Delta^2(y) \ds1_{|y|> B}\\
        \qquad\qquad\qquad\qquad
                \le 2\lam b^2 q^2 (|y|+b)^{2q-2} \ds1_{|y|> B}\\
        \qquad\qquad\qquad\qquad
                \le C'' \lam |y|^{2q-2}\ds1_{|y|> B},
        \endar
where $C''$ depends only on $b$ and $q$.

Therefore, from \eqref{eq:genpow}, for any $0<b\le B$ and $\lam>0$ satisfying $\lam \le \frac{q-1}{q(B+b)^q}$, 
the expected risk is upper bounded by
        \begarlab{eq:prp}
        \und{\min}{\rho\in\M} \Big\{ \expec{g}{\rho} R(g) + \frac{K(\rho,\pi)}{\lam(n+1)} \Big\}
                + \E\big\{ C' |Y|^{q-1} \ds1_{|Y|\ge b+(2bq\lam)^{-1/(q-1)} ; |Y| > B} \big\}\\
        \qquad\qquad\qquad\qquad\qquad\quad \;\;
                + \E\big\{ C'' \lam |Y|^{2q-2} \ds1_{B < |Y|< b+(2bq\lam)^{-1/(q-1)}} \big\}.
        \endarlab
Let us take $B=\big(\frac{q-1}{q\lam}\big)^{1/q}-b$ with $\lam$ small enough to ensure 
that $b\le B \le b+ (2bq\lam)^{-1/(q-1)}$. This means that $\lam$ should be taken smaller than
some positive constant depending only on $b$ and $q$. Then \eqref{eq:prp} can be written as 
        \begar
        \und{\min}{\rho\in\M} \Big\{ \expec{g}{\rho} R(g) + \frac{K(\rho,\pi)}{\lam(n+1)} \Big\}
                + \E\big\{ C' |Y|^{q-1} \ds1_{|Y|\ge b+(2bq\lam)^{-1/(q-1)}} \big\}\\
        \qquad\qquad\qquad\qquad\qquad\quad \;\;
                + \E\big\{ C'' \lam |Y|^{2q-2} \ds1_{(\frac{q-1}{q\lam})^{1/q}-b < |Y|< b+(2bq\lam)^{-1/(q-1)}} \big\}.
        \endar
Now using \eqref{eq:moments}, we can upper bound \eqref{eq:prp} with
        \begar
        \und{\min}{g\in\G} R(g) + \frac{\log |\G|}{\lam n} 
                + C \lam^{\frac{s+1-q}{q-1}}
                + C \lam \Big( \lam^{\frac{s-2q+2}{q}} \ds1_{s\ge 2q-2} + \lam^{\frac{2-2q+s}{q-1}} \ds1_{s<2q-2} \Big)
        \endar
where $C$ depends only on $b$, $A$, $q$ and $s$. So we get
        \begar
        \Ezun \frac{1}{n+1} \sum_{i=0}^n R(\expec{g}{\pi_{-\lam \Sigma_i}}g)\\
        \qquad\qquad\qquad
                \le \und{\min}{g\in\G} R(g) + \frac{\log |\G|}{\lam n} 
                + C \lam^{\frac{s+1-q}{q-1}} + C \lam^{\frac{s-q+2}{q}} \ds1_{s\ge 2q-2}\\
        \qquad\qquad\qquad
                \le \und{\min}{g\in\G} R(g) + \frac{\log |\G|}{\lam n} 
                + C \lam^{\frac{s+1-q}{q-1}} \ds1_{s< 2q-2} + C \lam^{\frac{s-q+2}{q}} \ds1_{s\ge 2q-2},
        \endar
since $\frac{s+1-q}{q-1} \ge \frac{s-q+2}{q}$ is equivalent to $s\ge 2q-2$. By taking 
$\lam$ of order of the minimum of the $\rhs$ (which implies that $\lam$ goes to $0$ when $n/\log|\G|$ goes to infinity), we obtain the desired result.

\else
\fi

\ifRR
\subsection{Proof of Lemma \ref{lem:convex}} \label{sec:prooflemconvex}

Let $\Delta_\ell= f'_r(t_\ell)-f'_l(t_\ell)$.
We have
        \begin{equation*}
        \begin{split}
        & f(\alpha) -  \alpha f(1) - (1-\alpha) f(0)\\
        & \qquad\qquad
                = f(\alpha)-f(0) - \alpha [f(1)-f(0)]\\
        & \qquad\qquad
                = \int_0^\alpha f'(u)du - \alpha \int_0^1 f'(u)du\\
        & \qquad\qquad
                = \int_0^\alpha \big( \int_0^u f''(t) dt + \sum_{\ell:t_\ell\in]0;u[} \Delta_\ell \big) du \\
        & \qquad\qquad\qquad\qquad\qquad
                - \alpha \int_0^1 \big( \int_0^u f''(t) dt + \sum_{\ell:t_\ell\in]0;u[} \Delta_\ell \big) du\\
        & \qquad\qquad
                = \int_{[0;1]^2} \big( \ds1_{0<t<u<\alpha} - \alpha \ds1_{0<t<u<1} \big) f''(t) dtdu\\
        & \qquad\qquad\qquad\qquad\qquad
                + \sum_{\ell:t_\ell\in (0;1)} [ (\alpha-t_\ell)\ds1_{t_\ell<\alpha} - \alpha ( 1 - t_\ell ) ] \Delta_\ell\\
        & \qquad\qquad
                = - \int_{[0;1]} K_\alpha(t) f''(t) dt
                - \sum_{\ell:t_\ell\in]0;1[} K_{\alpha}( t_\ell ) \Delta_\ell
        \end{split}
        \end{equation*}
\fi

\secproofth{assouad} 
The symbols $\sigma_1,\dots,\sigma_m$ still denote the coordinates of $\sigmav \in \{-;+\}^m$.
For any $r\in\{-;0;+\}$, define 
        $\sigmav_{j,r} \eqdef 
        (\sigma_1,\dots,\sigma_{j-1},r,\sigma_{j+1},\dots,\sigma_m)$
as the vector deduced from $\sigmav$ by fixing its $j$-th coordinate to $r$.
Since $\sigmav_{j,+}$ and $\sigmav_{j,-}$ belong to $\{-;+\}^m$, we have already 
defined $P_{\sigmav_{j,+}}$ and $P_{\sigmav_{j,-}}$.
Now we define the distribution $P_{\sigmav_{j,0}}$ as
        $P_{\sigmav_{j,0}}(dX) = \mu(dX)$ and 
        \begar
        1-P_{\sigmav_{j,0}}( Y = \hb(X) | X )\\
        \qquad\qquad\qquad
                = P_{\sigmav_{j,0}}( Y = \ha(X) | X ) = \left\{ \begin{array}{lll}
                \frac{1}{2} \text{ for any } X \in \X_j\\
                P_{\sigmav}( Y = \ha(X) | X ) \text{ otherwise}
                \end{array} \right..
        \endar
The distribution $P_{\sigmav_{j,0}}$ differs from $P_{\sigmav}$ only by the
conditional law of the output knowing that the input is in $\X_j$.
\ifRR
We recall that $P^{\otimes n}$ denotes the $n$-fold product of a distribution $P$.
\fi
For any $r\in\{-;+\}$, introduce \ifRR the likelihood ratios for the data $Z_1^n=(Z_1,\dots,Z_n)$:\fi
        \begar
        \pi_{r,j} (Z_1^n) \triangleq \frac{ P_{\sigmav_{j,r}}^{\otimes n}}
                { P_{\sigmav_{j,0}}^{\otimes n}}(Z_1^n)
        \endar
Note that this quantity is independent of the value of $\sigmav$.
In the following, to shorten the notation,
we will sometimes use $\ha$ for $\ha(X)$, $\hb$ for $\hb(X)$, $p_+$ for $p_{+}(X)$,
$p_-$ for $p_{-}(X)$.
Let $\nu$ be the uniform distribution on $\{-,+\}$, i.e. 
        $\nu\big(\{+\}\big) = 1/2 = 1-\nu\big(\{-\}\big).$
In the following, $\E_{\sigmav}$ denotes the expectation when $\sigmav$ is drawn according
to the $m$-fold product distribution of $\nu$, and $\E_X=\E_{X\sim\mu}$.
We have
        \begarlab{eq:assouadarg}
        & \underset{P \in \calP}{\sup} \Big\{ \PN R(\hg)
                - \min_g R(g) \Big\}\\ 
        \geq & \underset{\sigmav \in \{-;+\}^m}{\sup} \Big\{
                \undc{\E}{Z_1^n\sim P_{\sigmav}^{\otimes n}} 
                \expec{Z}{P_\sigmav} \ell[Y,\hg(X)] - \min_g \expec{Z}{P_\sigmav} \ell[Y,g(X)] \Big\}\\ 
        = & \underset{\sigmav \in \{-;+\}^m}{\sup} \Big\{
                \undc{\E}{Z_1^n\sim P_{\sigmav}^{\otimes n}} 
                \expec{X}{P_\sigmav(dX)} \Big[ \expec{Y}{P_\sigmav(dY|X)} \ell[Y,\hg(X)]\\ 
        & \qquad\qquad\qquad\qquad\qquad\qquad\qquad\qquad      
                - \und{\min}{y\in\Y} \, \expec{Y}{P_\sigmav(dY|X)} \ell(Y,y) \Big] \Big\}\\ 
        = & \und{\sup}{\sigmav \in \{-;+\}^m} \Big\{
                \undc{\E}{Z_1^n\sim P_{\sigmav}^{\otimes n}} 
                \expecc{X}{\mu} \Big[ \sum_{j=0}^m \ds1_{X\in\X_j}\\
        & \qquad\qquad
                \times \Big( \varphi_{p_{\sigma_j},\ha,\hb}[\hg(X)] 
                - \phi_{\ha,\hb}[p_{\sigma_j}] \Big) \Big] \Big\}\\ 
        \geq & \expecc{\sigmav}{\nuv} \undc{\E}{Z_1^n\sim P_{\sigmav}^{\otimes n}} 
                \expecc{X}{\mu} \Big[ \sum_{j=1}^m \ds1_{X\in\X_j}\\
        & \qquad\qquad
                \times \Big( \varphi_{p_{\sigma_j},\ha,\hb}[\hg(X)] 
                -\phi_{\ha,\hb} [p_{\sigma_j}]\Big) \Big]\\ 
        = & \sum_{j=1}^m \expecc{X}{\mu} \Big\{ \ds1_{X\in\X_j} \expecc{\sigmav}{\nuv} 
                \expec{Z_1^n}{P_{\sigmav_{j,0}}^{\otimes n}} \Big[ \frac{ P_{\sigmav}^{\otimes n}}
                {P_{\sigmav_{j,0}}^{\otimes n}}(Z_1^n)
                \\
        & \qquad\qquad
                \times \Big( \varphi_{p_{\sigma_j},\ha,\hb}[\hg(X)] 
                -\phi_{\ha,\hb}[p_{\sigma_j}] \Big) \Big] \Big\}\\ 
        = & \sum_{j=1}^m \expecc{X}{\mu} \Big\{ \ds1_{X\in\X_j} \expecc{\sigma_1,\dots,\sigma_{j-1},
                \sigma_{j+1},\dots,\sigma_m}{\nu^{\otimes (m-1)}} 
                \expec{Z_1^n}{P_{\sigmav_{j,0}}^{\otimes n}} \\
        & \qquad\qquad
                \expec{\sigma_j}{\nu} \pi_{\sigma_j,j}(Z_1^n) \Big(
                 \varphi_{p_{\sigma_j},\ha,\hb}[\hg(X)] 
                -\phi_{\ha,\hb}[p_{\sigma_j}] \Big) \Big\}
        \endarlab
The two inequalities in \eqref{eq:assouadarg} are Assouad's argument (\cite{Ass83}).
For any $x\in\X$, introduce 
        \ifRR
        \begar
        \tpsi_{x}(u) \eqdef \demi ( u + 1) \psi_{p_+(x),p_-(x),\ha(x),\hb(x)}\big(\frac{u}{ u + 1} \big).
        \endar  
        \fi
Introduce
\ifRR
        \begar
        \alpha_j(Z^n_1) = \frac{\pi_{+,j} (Z_1^n)}{\pi_{+,j} (Z_1^n) + \pi_{-,j} (Z_1^n)}.
        \endar
\fi
The last expectation in \eqref{eq:assouadarg} is
        \begarlab{eq:ass1}
        \expec{\sigma}{\nu} \pi_{\sigma,j}(Z_1^n) \Big(
                 \varphi_{p_{\sigma}(X),\ha(X),\hb(X)}[\hg(X)] 
                -\phi_{\ha(X),\hb(X)}[p_{\sigma}(X)] \Big) \\
        \qquad\qquad
                = \demi \big[ \pi_{+,j}(Z_1^n) + \pi_{-,j}(Z_1^n)\big]\\
        \qquad\qquad\qquad\qquad 
                \times \Big\{ \alpha_j(Z_1^n) \varphi_{p_+,\ha,\hb}[\hat{g}(X)]
                + [1-\alpha_j(Z_1^n)] \varphi_{p_-,\ha,\hb}[\hat{g}(X)]\\
        \qquad\qquad\qquad\qquad
                \quad - \alpha_j(Z_1^n) \phi_{\ha,\hb}(p_+)
                - [1-\alpha_j(Z_1^n)] \phi_{\ha,\hb}(p_-)\Big\}\\
        \qquad\qquad
                = \demi \big[ \pi_{+,j}(Z_1^n) + \pi_{-,j}(Z_1^n)\big]
                \Big\{ \varphi_{\alpha_j(Z_1^n) p_++ [1-\alpha_j(Z_1^n)] p_- ,\ha,\hb}[\hat{g}(X)]\\
        \qquad\qquad\qquad\qquad
                \quad - \alpha_j(Z_1^n) \phi_{\ha,\hb}(p_+)
                - [1-\alpha_j(Z_1^n)] \phi_{\ha,\hb}(p_-)\Big\}\\
        \qquad\qquad
                \ge \demi \big[ \pi_{+,j}(Z_1^n) + \pi_{-,j}(Z_1^n)\big]
                \Big\{ \phi_{\ha,\hb}\Big(\alpha_j(Z_1^n) p_+ + [1-\alpha_j(Z_1^n)] p_-\Big)\\
        \qquad\qquad\qquad\qquad
                \quad - \alpha_j(Z_1^n) \phi_{\ha,\hb}(p_+)
                - [1-\alpha_j(Z_1^n)] \phi_{\ha,\hb}(p_-)\Big\}\\
        \qquad\qquad
                = \demi \big[ \pi_{+,j}(Z_1^n) + \pi_{-,j}(Z_1^n)\big]
                \psi_{p_+,p_-,\ha,\hb}[\alpha_j(Z_1^n)]\\
        \qquad\qquad
                = \pi_{-,j}(Z_1^n) \tpsi_{X} \Big( \frac{\pi_{+,j}(Z_1^n)}{\pi_{-,j}(Z_1^n)} \Big)
        \endarlab
\ifRR so that \fi
        \begar
        \underset{P \in \calP}{\sup} \Big\{ \PN R(\hg)
                - \min_g R(g) \Big\}\\ 
        \qquad\qquad\qquad
        \ge \sum_{j=1}^m \expecc{X}{\mu} \Big\{ \ds1_{X\in\X_j} \expecc{\sigmav}{\nu^{\otimes (m-1)}} 
                \expec{Z_1^n}{P_{\sigmav_{j,0}}^{\otimes n}} 
                \Big[ \pi_{-,j}(Z_1^n) \tpsi_{X} \Big( \frac{\pi_{+,j}(Z_1^n)}{\pi_{-,j}(Z_1^n)} \Big)\Big]\Big\}\\
        \qquad\qquad\qquad
        = \sum_{j=1}^m \expecc{X}{\mu} \Big\{ \ds1_{X\in\X_j} \expecc{\sigmav}{\nu^{\otimes (m-1)}} 
                \S_{\tpsi_{X}}\big( P_{\sigmav_{j,+}}^{\otimes n} , P_{\sigmav_{j,-}}^{\otimes n} \big) \Big\}.
        \endar
Now since we consider a hypercube, for any $j\in\{1,\dots,m\}$, all the terms in the sum are equal. 
Besides 
\ifRR from part \ref{item:2} of Lemma \ref{lem:sim}, 
\fi
the last $f$-similarity does not depend on $\sigmav$, and in particular
for $j=1$, the $f$-similarity is equal to $\S_{\tpsi_{X}}\big( P_{[+]}^{\otimes n} , P_{[-]}^{\otimes n} \big)$\ifRR,
where we recall that $P_{[+]}$ and $P_{[-]}$ denote the representatives of the hypercube
(see Definition \ref{def:hyp2})\fi.
Therefore we obtain
        \begar  
        \underset{P \in \calP}{\sup} \Big\{ \E R(\hg)
                - \min_g R(g) \Big\}
        & \ge & m \expecc{X}{\mu} \Big\{ \ds1_{X\in\X_1} 
                \S_{\tpsi_{X}}\big( P_{[+]}^{\otimes n} , P_{[-]}^{\otimes n} \big) \Big\}\\
        & = & m \expecc{X}{\mu} \S_{\tpsi_{X} \ds1_{X\in\X_1} }\big( P_{[+]}^{\otimes n} , P_{[-]}^{\otimes n} \big) \\
        \ifRR & = & m \S_{\expecc{X}{\mu} ( \ds1_{X\in\X_1} \tpsi_{X} )}
                \big( P_{[+]}^{\otimes n} , P_{[-]}^{\otimes n} \big) \\\fi
        & = & \S_{\tpsi}\big( P_{[+]}^{\otimes n} , P_{[-]}^{\otimes n} \big)
        \endar
\ifRR where the second to last equality comes from the second part of Lemma~\ref{lem:sim}.\fi


\ifRR
\subsection{Proof of Theorem \ref{th:appassouad}} \label{sec:proofappassouad}

First, when the hypercube satisfies $\pa\equiv 1 \equiv 1-\pb$, from the definition of $\da$
given in \eqref{eq:defed1a}, we have
        $\S_{\tpsi} \big( P_{[+]}^{\otimes n} , P_{[-]}^{\otimes n} \big)
                = mw \da (1-w)^n$
so that Theorem \ref{th:assouad} implies \eqref{eq:w1c}.

Inequalities \eqref{eq:w1a}, \eqref{eq:w1b} and \eqref{eq:w2} are deduced from 
Theorem \ref{th:assouad} by lower bounding the $\tpsi$-similarity in different ways.

Since $u \mapsto u\wedge 1$ and $u \mapsto \frac{u}{u + 1}$
are non-negative concave functions defined on $\Rp$,
we may define the similarities
        \lbegar
        \S_{\wedge}(\P,\Q) \eqdef \int \big(\frac{\P}{\Q} \wedge 1 \big)d\Q = \int \big(d\P \wedge d\Q \big) \\
        \S_{\point}(\P,\Q) \eqdef \int \frac{\P}{\Q} \frac{1}{\frac{\P}{\Q}+1} d\Q 
                = \int \frac{d\P d\Q}{d\P+d\Q} 
        \rendarp
where the second equality of both formulas introduces a formal (but intuitive) notation.

From Theorem \ref{th:assouad}, Lemma \ref{lem:tpsi} and item 1 of Lemma \ref{lem:sim}, by using
$\tpsi(1)=mw \da$, we obtain

\begin{corollary}
Let $\calP$ be a set of probability distributions containing a hypercube of distributions 
of characteristic function $\tpsi$ and representatives $P_{[-]}$ and $P_{[+]}$.
For any estimator $\hg$, we have
        \begarlab{eq:thass2}
        \underset{P \in \calP}{\sup} \big\{ \E R( \hg ) 
                - \und{\min}{g} R(g) \big\} 
        & \ge & mw \da \S_{\wedge} \big( P_{[+]}^{\otimes n} , P_{[-]}^{\otimes n} \big) 
        \endarlab
where the minimum is taken over the space of prediction functions.
Besides if 
for any $x\in\X_1$ the function $\phi_{\ha(x),\hb(x)}$ is twice differentiable 
and satisfies for any $t\in[p_-(x) \wedge p_+(x);p_-(x) \vee p_+(x)]$,
        $
        -\phi''_{\ha(x),\hb(x)}(t) \ge \zeta
        $
for some $\zeta>0$, then we have
        \begarlab{eq:thass5}
        \underset{P \in \calP}{\sup} \big\{ \E R( \hg ) 
                - \und{\min}{g} R(g) \big\} 
                \ge \frac{m w \zeta}{4} \da' \; 
                \S_{\point} \big( P_{[+]}^{\otimes n} , P_{[-]}^{\otimes n} \big) ;
        \endarlab
\end{corollary}
The following lemma and \eqref{eq:thass2} imply \eqref{eq:w1a} and \eqref{eq:w1b}. 

\begin{lemma} \label{lem:wun}
We have
        \begarlab{eq:swedgen1a}
        \S_{\wedge} \big( P_{[+]}^{\otimes n} , P_{[-]}^{\otimes n} \big) 
                \ge 1-\sqrt{1 - [1-\db]^{nw}}
                \ge 1-\sqrt{nw \db }.
        \endarlab
When the hypercube is symmetrical and constant, for $N$ 
a centered gaussian random variable with variance $1$, we have
        \begarlab{eq:swedgen2}
        \S_{\wedge} \big( P_{[+]}^{\otimes n} , P_{[-]}^{\otimes n} \big) 
                \ge \P\Big(|N| > \sqrt{\frac{nw\db}{1-\db}} \Big) - \db^{1/4}
        \endarlab

\end{lemma} 

\jyproof{lem:wun}

\begin{remark}
It is interesting to note that \eqref{eq:swedgen2} 
is asymptotically optimal to the extent that for a $(m,w,\db)$-hypercube (see Definition \ref{def:mwxi}),
we have
        \begarlab{eq:swedgen3}
        \Big| \S_{\wedge} \big( P_{[+]}^{\otimes n} , P_{[-]}^{\otimes n} \big) 
                - \P\big( |N| > \sqrt{n w\db} \big) \Big|
        \und{\longrightarrow}{nw\rightarrow+\infty,\db\rightarrow 0} 0,
        \endarlab
[Proof in Appendix \ref{app:proofasympt}.]
\end{remark}

The following lemma and \eqref{eq:thass5} imply \eqref{eq:w2}.
\begin{lemma} \label{lem:wdeux}
When the hypercube is symmetrical and constant, we have
        \begar
        \S_{\bullet} \Big( P_{[+]}^{\otimes n} , P_{[-]}^{\otimes n} \Big) 
                \ge \demi \Big\{ 1+\demi \big[1-\big(1-\sqrt{1-\db}\big)w\big]^n 
                - \demi \big[1+\big(\frac{1+\db}{\sqrt{1-\db}}-1\big)w\big]^n \Big\}
        \endar
\end{lemma}
\jyproof{lem:wdeux}

\subsecprooflem{wun}

For $\sigma\in\{-,+\}$, the conditional law of $(X,Y)$ knowing $X\in\X_1$, 
when $(X,Y)$ follows the law $P_{[\sigma]}$, is denoted
$P_{\X_1,\sigma}$ and is called the \emph{restricted representatives of the hypercube}.
More explicitly, the probability distribution $P_{\X_1,\sigma}$ is such that
its first marginal $P_{\X_1,\sigma}(dX)$ is $\mu(\bullet|\X_1)$ and for any $x\in\X_1$
        \begar
        P_{\X_1,\sigma}\big( Y = \ha(x) \big| X=x \big) = p_{\sigma}(x)
                        = 1 - P_{\X_1,\sigma}\big( Y = \hb(x) \big| X=x \big).
        \endar

The following lemma relates the similarity between representatives 
of the hypercube and the similarity between restricted representatives.

\begin{lemma} \label{lem:simconv}
Consider a convex function $\gamma:\R_+\rightarrow\R_+$ such
that 
        \begar
        \gamma(k)\le \S_{\wedge} \big( P_{\X_1,+}^{\otimes k} , P_{\X_1,-}^{\otimes k} \big) 
        \endar
for any $k\in\{0,\dots,n\}$,
where by convention $\S_{\wedge} \big( P_{\X_1,+}^{\otimes 0} , P_{\X_1,-}^{\otimes 0} \big)=1$.
For any estimator $\hg$, we have
        \begar
         \S_{\wedge} \big( P_{[+]}^{\otimes n} , P_{[-]}^{\otimes n} \big) 
                \ge \gamma(n w).
        \endar
\end{lemma}

\begin{proof}
For any points $z_1=(x_1,y_1),\dots,z_n=(x_n,y_n)$ in $\X\times\{h_1,h_2\}$, let $\C(z_1,\dots,z_n)$ denotes the number of $z_i$ for which $x_i\in\X_1$.
For any $k\in\{0,\dots,n\}$, let $B_k=\C^{-1}(\{k\})$ denote the subset of $\big(\X\times\{h_1,h_2\}\big)^n$ for which
exactly $k$ points are in $\X_1\times\{h_1,h_2\}$. We recall that there are $\binom{n}{k}$ possibilities of taking $k$ elements among $n$ and
the probability of $X\in\X_1$ when $X$ is drawn according to $\mu$ is $w=\mu(\X_1)$.
Let $\Z_1=\X_1 \times \{h_1,h_2\}$ and let $\Z_1^c$ denote the complement of $\Z_1$.
We have
        \begarlab{eq:simdecomp}
        & \S_{\wedge} \big( P_{[+]}^{\otimes n} , P_{[-]}^{\otimes n} \big)\\
        = & \int 1\wedge\Big( \frac{P_{[+]}^{\otimes n}}{P_{[-]}^{\otimes n}} (z_1,\dots,z_n) \Big) 
                dP_{[-]}^{\otimes n}(z_1,\dots,z_n)\\
        = & \sum_{k=0}^n
                \int_{B_k} 1\wedge\Big( \frac{P_{[+]}}{P_{[-]}} (z_1)\cdots \frac{P_{[+]}}{P_{[-]}}(z_n) \Big) 
                dP_{[-]}(z_1)\cdots dP_{[-]}(z_n)\\
		= & \sum_{k=0}^n \binom{n}{k} 
                \int_{(\Z_1)^k\times(\Z_1^c)^{n-k}} 1\wedge\Big( \frac{P_{[+]}}{P_{[-]}} (z_1)\cdots 
                \frac{P_{[+]}}{P_{[-]}}(z_n) \Big) dP_{[-]}^{\otimes n}(z_1,\dots,z_n)\\
        = & \sum_{k=0}^n \binom{n}{k} 
                \int_{(\Z_1)^k\times(\Z_1^c)^{n-k}} 1\wedge\Big( \frac{P_{[+]}}{P_{[-]}} (z_1)\cdots 
                \frac{P_{[+]}}{P_{[-]}}(z_k) \Big) dP_{[-]}^{\otimes n}(z_1,\dots,z_n)\\
        = & \sum_{k=0}^n \binom{n}{k} \mu^{n-k}(\Z_1^c) 
                \int_{(\Z_1)^k} 1\wedge\Big( \frac{P_{[+]}}{P_{[-]}} (z_1)\cdots 
                \frac{P_{[+]}}{P_{[-]}}(z_k) \Big) dP_{[-]}^{\otimes n}(z_1,\dots,z_k)\\
        = & \sum_{k=0}^n \binom{n}{k} \mu^{n-k}(\Z_1^c) \mu^k(\Z_1) 
                \S_{\wedge} \big( P_{\X_1,+}^{\otimes k} , P_{\X_1,-}^{\otimes k} \big)\\
        \ge & \sum_{k=0}^n \binom{n}{k} (1-w)^{n-k} w^k \gamma(k)\\
		= & \E \gamma(V) 
        \endarlab
where $V$ is a Binomial distribution with parameters $n$ and $w$.
By Jensen's inequality, we have $\E \gamma(V) \ge \gamma[\E (V)]=\gamma(nw),$
which ends the proof.
\end{proof}

The interest of the previous lemma is to provide a lower bound on the similarity between
representatives of the hypercube from a lower bound on the similarity between restricted representatives,
restricted representatives being much simpler to study.
The following result lower bounds the $\wedge$-similarity between the restricted representatives of the hypercube.

\begin{lemma} \label{lem:restr}
For any non-negative integer $k$, we have
        \begarlab{eq:swedge1a}
        \S_{\wedge} \big( P_{\X_1,+}^{\otimes k} , P_{\X_1,-}^{\otimes k} \big)
                \ge 1-\sqrt{1 - [1-\db]^k}
                \ge 1-\sqrt{k \db },
        \endarlab
When the hypercube is symmetrical and constant, for $N$ 
a centered gaussian random variable with variance $1$, we have
        \begarlab{eq:swedge2}
        \S_{\wedge} \big( P_{\X_1,+}^{\otimes k} , P_{\X_1,-}^{\otimes k} \big)
                \ge \P\Big(|N| > \sqrt{\frac{k\db}{1-\db}} \Big) - \db^{1/4}.
        \endarlab
\end{lemma}

\begin{proof}
First, we recall that $P_0$ denotes the base of the hypercube (see Definition \ref{def:hyp2}).
The conditional law of $(X,Y)$ knowing $X\in\X_1$, 
when $(X,Y)$ is drawn from $P_0$, is denoted
$\Prestz$. 

For any $r\in\{-,0,+\}$, introduce $\Prr{x}$ the probability distribution on the output space such that
        $\Prr{x}(dY) = \Prestr(dY|X=x)$.
We have
        {\renewcommand{\arraystretch}{1.5}
        \begarlab{eq:lbpr1}
        \S_\wedge\big( \Prestp^{\otimes k} , \Prestm^{\otimes k} \big)
                = \expec{Z_1^k}{\Prestz^{\otimes k} } \Big[\frac{\Prestp^{\otimes k}}{\Prestz^{\otimes k}}(Z_1^k) 
                \wedge \frac{\Prestm^{\otimes k}}{\Prestz^{\otimes k}}(Z_1^k) \Big]\\
        \qquad\quad
                = \expec{X_1^k}{\Prestz^{\otimes k} } \E_{Y_1^k\sim\Prestz^{\otimes k}|X_1^k} 
                \Big[ \prod_{i=1}^k \frac{\Prp{X_i}}{\Prz{X_i}}(Y_i) 
                \wedge \prod_{i=1}^k \frac{\Prm{X_i}}{\Prz{X_i}}(Y_i)  \Big]\\
        \qquad\quad
         = \expec{X_1^k}{\Prestz^{\otimes k} } \S_\wedge\Big( \otimes_{i=1}^k \Prp{X_i} , \otimes_{i=1}^k \Prm{X_i} \Big),
        \endarlab \par}
\noindent where $\otimes_{i=1}^k \Prr{X_i}$, $r\in\{-1;1\}$ denotes the law of the $k$-tuple
$(Y_1,\dots,Y_k)$ when the $Y_i$ are independently drawn from $\Prr{X_i}$.

To study divergences (or equivalently similarities) between $k$-fold product distributions, 
the standard way is to link the divergence (or similarity) of the product with
the ones of base distributions. This lead to tensorization equalities or inequalities.
To obtain a tensorization inequality for $\S_\wedge$, we introduce the similarity
associated with the square root function (which is non-negative and concave):
        \begar
        \S_{\sqrt{\;}}(\P,\Q) \eqdef \int \sqrt{d\P d\Q}
        \endar
and use the following lemmas:
\begin{lemma}
For any probability distributions $\P$ and $\Q$, we have 
        \begar
        \S_{\wedge}(\P,\Q) \ge 1 - \sqrt{1-\S_{\sqrt{\;}}^2(\P,\Q)}.
        \endar
\end{lemma}
\begin{proof}
Introduce the variational distance $V(\P,\Q)$ as the $f$-divergence associated with
the convex function $f: u \mapsto \frac{1}{2} |u-1|$. From Scheff\'e's theorem, we
have $\S_{\wedge}(\P,\Q) = 1 - V(\P,\Q)$ for any distributions $\P$ and $\Q$.
Introduce the Hellinger distance $H$, which is defined as $H(\P,\Q)\ge 0$ and 
        $1 - \frac{H^2( \P , \Q )}{2} = \S_{\sqrt{\;}}(\P,\Q)$
for any probability distributions $\P$ and $\Q$.
The variational and Hellinger distances are known (see e.g. \cite[Lemma 2.2]{Tsy04c})
to be related by 
        \begar
        V( \P,\Q ) \le \sqrt{1-\big(1-\frac{H^2(\P,\Q)}{2}\big)^2},
        \endar
hence the result.
\end{proof}
\begin{lemma}
For any distributions $\P^{(1)},\dots,\P^{(k)},$ $\Q^{(1)},\dots,\Q^{(k)}$, we have
        \begar
        \S_{\sqrt{\;}}( \P^{(1)}\otimes \cdots \otimes \P^{(k)} , \Q^{(1)}\otimes \cdots \otimes \Q^{(k)} )\\
        \qquad\qquad\qquad\qquad\qquad
                = \S_{\sqrt{\;}}( \P^{(1)} , \Q^{(1)} )
                \times \cdots \times \S_{\sqrt{\;}}( \P^{(k)} , \Q^{(k)} )
        \endar
\end{lemma}
\begin{proof}
When it exists, the density of $\P^{(1)}\otimes \cdots \otimes \P^{(k)}$ $\wrt$ $\Q^{(1)}\otimes \cdots \otimes \Q^{(k)}$
is the product of the densities of $\P^{(i)}$ $\wrt$ $\Q^{(i)}$, $i=1,\dots,k$, hence the desired tensorization equality.

\end{proof}

From the last two lemmas, we obtain     
        \begarlab{eq:wedgetensor}
        \S_\wedge\Big( \otimes_{i=1}^k \Prp{X_i} , \otimes_{i=1}^k \Prm{X_i} \Big)
                \ge 1 - \sqrt{1-\prod_{i=1}^k 
                S^2_{\sqrt{\;}}\Big( \Prp{X_i} , \Prm{X_i} \Big)}
        \endarlab
From \eqref{eq:lbpr1}, \eqref{eq:wedgetensor} and Jensen's inequality, we obtain        
        \begar
        \S_\wedge\big( \Prestp^{\otimes k} , \Prestm^{\otimes k} \big)
                \ge 1 - \expec{X_1^k}{\Prestz^{\otimes k} } \sqrt{1-\prod_{i=1}^k 
                S^2_{\sqrt{\;}}\Big( \Prp{X_i} , \Prm{X_i} \Big)}\\
                \ge 1 - \sqrt{1-\expec{X_1^k}{\Prestz^{\otimes k} } \prod_{i=1}^k 
                S^2_{\sqrt{\;}}\Big( \Prp{X_i} , \Prm{X_i} \Big)}\\
                = 1 - \sqrt{1- \Big[\expec{X}{\Prestz} S^2_{\sqrt{\;}}\Big( \Prp{X} , \Prm{X} \Big)\Big]^k}\\
                = 1 - \sqrt{1- \Big[\expec{X}{\mu(\bullet|\X_1)} S^2_{\sqrt{\;}}\Big( \Prp{X} , \Prm{X} \Big)\Big]^k}\\
        \endar
Now we have
        \begar
        \expec{X}{\mu(\bullet|\X_1)} S^2_{\sqrt{\;}}\Big( \Prp{X} , \Prm{X} \Big)\\
        \qquad\qquad
                = \expecd{X}{\mu(\bullet|\X_1)} \big[ \sqrt{\phantom{(}\pa\,\pb\phantom{)}}+\sqrt{(1-\pa)(1-\pb)}\,\big]^2\\
        \qquad\qquad
                = 1 - \expecd{X}{\mu(\bullet|\X_1)} \big[\sqrt{\pa(1-\pb)}-\sqrt{(1-\pa)\pb}\,\big]^2\\
        \qquad\qquad
                = 1 - \db
        \endar
So we get
        \begarlab{eq:reflbswedge}
        \S_\wedge\big( \Prestp^{\otimes k} , \Prestm^{\otimes k} \big)
                \ge 1 - \sqrt{1- (1-\db)^k} \ge 1 - \sqrt{k\db},
        \endarlab
where the second inequality follows 
from the inequality $1-x^k\le k (1-x)$ that holds for any $0\le x\le 1$ and $k\ge 1$.
This ends the proof of \eqref{eq:swedge1a}.

For \eqref{eq:swedge2}, since we assume that the hypercube is symmetrical and constant, we can tighten
\eqref{eq:reflbswedge} for $k\sqrt{\db}\ge 1$.
We have
        \begarlab{eq:s1}
        \S_\wedge\big( \Prestp^{\otimes k} , \Prestm^{\otimes k} \big)
                = \Prestp^{\otimes k}\Big( \frac{\Prestp^{\otimes k}}{\Prestm^{\otimes k}}(Z_1^k)\le 1 \Big)
                + \Prestm^{\otimes k}\Big( \frac{\Prestp^{\otimes k}}{\Prestm^{\otimes k}}(Z_1^k) > 1 \Big).
        \endarlab
Since $\frac{\Prestp}{\Prestm}(z) = \frac{\Prestp(Y=y|X=x)}{\Prestm(Y=y|X=x)} = \frac{\Prp{x}(y)}{\Prm{x}(y)}$ for any $z=(x,y)\in\Z$, we have
        \begarlab{eq:density}
        \frac{\Prestp^{\otimes k}}{\Prestm^{\otimes k}}(Z_1^k)
                & = & \prod_{i=1}^k \frac{\Prp{X_i}(Y_i)}{\Prm{X_i}(Y_i)} \\
        & = & \prod_{i=1}^k \Big( \frac{\pa(X_i)}{\pb(X_i)} \Big)^{\ds1_{Y_i=\ha(X_i)}}
                \Big( \frac{1-\pa(X_i)}{1-\pb(X_i)}\Big)^{\ds1_{Y_i=\hb(X_i)}}.
        \endarlab
Using that the hypercube is symmetrical and constant, \eqref{eq:density} leads to
        \begarlab{eq:s2}
        \frac{\Prestp^{\otimes k}}{\Prestm^{\otimes k}}(Z_1^k)
                & = & \prod_{i=1}^k \Big( \frac{\pa(X_i)}{1-\pa(X_i)} \Big)^{\ds1_{Y_i=\ha(X_i)}-\ds1_{Y_i=\hb(X_i)}} \\
                & = & \big( \frac{\pa}{1-\pa} \big)^{\sum_{i=1}^k [\ds1_{Y_i=\ha(X_i)}-\ds1_{Y_i=\hb(X_i)}]}.
        \endarlab
Without loss of generality, we may assume that $\pa > 1/2$. 
Then we have $\pa=1-\pb=\frac{1+\sqrt{\db}}{2}$.
Introduce $W_i\eqdef\ds1_{Y_i=\ha(X_i)}-\ds1_{Y_i=\hb(X_i)}$.
From \eqref{eq:s1} and \eqref{eq:s2}, we obtain
        \begar
        \S_\wedge\big( \Prestp^{\otimes k} , \Prestm^{\otimes k} \big)
                & = & \Prestp^{\otimes k}\big( \sum_{i=1}^k W_i \le 0 \big)
                + \Prestm^{\otimes k}\big( \sum_{i=1}^k W_i  > 0 \big)\\
                & = & \Prestm^{\otimes k}\big( \sum_{i=1}^k W_i \ge 0 \big)
                + \Prestm^{\otimes k}\big( \sum_{i=1}^k W_i  > 0 \big)
        \endar
The law of $U\eqdef\sum_{i=1}^k W_i$ when $(X_1,Y_1),\dots,(X_n,Y_n)$
are independently drawn from $\Prestm$ is the binomial distribution of parameter $\big(k,\frac{1-\sqrt{\db}}{2}\big)$.
Let $\integ{x}$ still denote the largest integer $k$ such that $k \le x$.
We get 
        \begar
        \S_\wedge\big( \Prestp^{\otimes k} , \Prestm^{\otimes k} \big) 
                & = & \P(U > k/2)+\P(U \ge k/2)\\
                & \ge & 2 \P\big( U \ge \integ{k/2} \big) - 2\P\big( U = \integ{k/2} \big) \\
        \endar
When $k\sqrt{\db}\ge 1$, this last $\rhs$ can be lower bounded by Slud's theorem \cite{Slu77} 
for the first term and by using 
Stirling's formula for the second term (see e.g. \cite[Appendix A.8]{Dev96}).
It gives
        \begar
        \S_\wedge\big( \Prestp^{\otimes k} , \Prestm^{\otimes k} \big) 
                & \ge & 2 \P\Big(N \ge \frac{2 \integ{k/2}-k(1-\sqrt{\db})}{\sqrt{k(1-\db)}} \Big) 
                - \sqrt{\frac{2}{k\pi}} \\
                & \ge & 2 \P\Big(N \ge \sqrt{\frac{k\db}{1-\db}} \Big) 
                - \sqrt{\frac{2}{\pi}} \db^{1/4},
        \endar
where we recall that $N$ is a normalized gaussian random variable. Finally we have
        \begar
        \S_\wedge\big( \Prestp^{\otimes k} , \Prestm^{\otimes k} \big) 
                \ge \left\{ \begin{array}{lll}
                1-\sqrt{k}\db^{1/4} & \text{for any }k\ge 1\\
                \P\Big(|N| > \sqrt{\frac{k\db}{1-\db}} \Big) 
                - \sqrt{\frac{2}{\pi}} \db^{1/4} & \text{for any }k\ge \frac{1}{\sqrt{\db}}\\
                \end{array} \right.
        \endar
which can be weakened into for any non-negative integer $k$     
        \begar
        \S_\wedge\big( \Prestp^{\otimes k} , \Prestm^{\otimes k} \big) 
                \ge \P\Big(|N| > \sqrt{\frac{k\db}{1-\db}} \Big) 
                - \db^{1/4},
        \endar  
that is \eqref{eq:swedge2}.
\end{proof}

%

By computing the second derivative of $u\mapsto \sqrt{1-e^{-u}}$ and 
$u \mapsto \int_0^{\sqrt{u}} e^{- t^2} dt$, we obtain that these functions are concave.
So for any $a\in[0;1]$, the functions
$x\mapsto 1-\sqrt{1-a^x}$, $x\mapsto 1-\sqrt{a x}$ 
and $x\mapsto \P\Big(|N| > \sqrt{\frac{xa}{1-a}} \Big) - a^{1/4}$ are convex.
The convexity of these functions and Lemmas \ref{lem:simconv} and \ref{lem:restr}
imply Lemma \ref{lem:wun}.

\subsecprooflem{wdeux}

Let $\theta:u \mapsto u/(u+1)$ denote the non-negative concave function on which the
similarity $\S_{\bullet}$ is defined.
For any $u>0$, we have
        \begar
        \theta(u) & = & \frac{1}{4}\big( u+1- \frac{(u-1)^2}{u+1} \big)\\
                & \ge & \frac{1}{4}\big( u+1- \frac{(u-1)^2}{2\sqrt{u}} \big)\\
                & = & \frac{1}{4}\big( u+1+\sqrt{u}- \frac{u^{3/2}}{2} - \frac{1}{2\sqrt{u}}\big),
        \endar
hence for any probability distributions $\P$ and $\Q$,
        \begar
        \S_{\bullet}\big( \P , \Q \big)
                & = & \int \theta\big(\frac{\P}{\Q}\big) d\Q\\
                & \ge & \frac{1}{4} \int \big( d\P+d\Q+\sqrt{d\P d\Q}- \frac{d\P^{3/2}}{2d\Q^{1/2}} 
                        - \frac{d\Q^{3/2}}{2d\P^{1/2}} \big)\\
                & = & \demi + \frac{1}{4} \int\sqrt{d\P d\Q} - \frac{1}{8} \int \frac{d\P^{3/2}}{d\Q^{1/2}} 
                        - \frac{1}{8} \int \frac{d\Q^{3/2}}{d\P^{1/2}} 
        \endar
The goal of this bound is to obtain a form for which tensorization equalities hold. Precisely,
let 
        $I_1 \eqdef \int \sqrt{dP_{[+]} dP_{[-]}}$
and
        $I_2 \eqdef \int \frac{dP_{[+]}^{3/2}}{dP_{[-]}^{1/2}} 
        = \int \frac{dP_{[-]}^{3/2}}{dP_{[+]}^{1/2}},$
where the last equality holds since the hypercube is symmetrical. We have
        \begar
        \S_{\bullet}\big( P_{[+]}^{\otimes n} , P_{[-]}^{\otimes n} \big)
                \ge \demi + \frac{1}{4} I_1^n - \frac{1}{4} I_2^n               
        \endar
Since the hypercube is symmetrical and constant, without loss of generality, we may assume 
that $\pa\ge\demi$ on $\X_1$. Then we have $1-\pb=\pa=(1+\sqrt{\db})/2$, hence
        $I_1 = 1-w+w\sqrt{1-\db}$ 
and 
        \begar
        I_2=1-w+\frac{w}{2} \Big( \frac{(1+\sqrt{\db})^{3/2}}{(1-\sqrt{\db})^{1/2}}
                + \frac{(1-\sqrt{\db})^{3/2}}{(1+\sqrt{\db})^{1/2}} \Big)
                = 1-w+w\frac{1+\db}{\sqrt{1-\db}},
        \endar
which gives the desired result.

\subsection{Proof of Theorems \ref{th:lqpt} and \ref{th:qsup1}} \label{sec:proofthlqpg}

We consider a $(\tm,\tw,\tdb)$-hypercube (see Definition \refp{def:mwxi}) with
        \begar
        \tm=\lfloor \log_2 |\G| \rfloor,
        \endar
$\ha\equiv -B$ and $\hb\equiv B$, and
with $\tw$ and $\tdb$ to be taken in order to (almost) maximize the bound.


\paragraph{Case $q=1$ : } 

From \eqref{eq:1power}, we have
        $
        \da = \frac{\sqrt{\db}}{2}|\hb-\ha|= B\sqrt{\db}
        $
so that, choosing $\tw=1/\tm$, \eqref{eq:w1a} gives 
        \begar
        \underset{P \in \H}{\sup} \big\{ \E R( \hg ) 
                - \und{\min}{g} R(g) \big\} & \ge & 
                B \sqrt{\db} \big(1-\sqrt{n\db/\tm}\big).
        \endar
Maximizing the lower bound $\wrt$ $\db$, we choose $\db=\frac{\tm}{4n}\wedge 1$ and obtain the announced result.

\paragraph{Case $1<q \le 1+\sqrt{\frac{\tm}{4n}\wedge 1}$ : } 
From \eqref{eq:dacstsym} and \eqref{eq:qpower}, for any $0<\eps\le1$, we have
        \begar
        \da & \ge & \frac{\db}{2} 
                \int_{\frac{1-\eps}{2}}^{\frac{1+\eps}{2}} [t\wedge (1-t)]
                \Big| \phi_{\ha,\hb}''\Big(\frac{1-\sqrt{\db}}{2}+\sqrt{\db}t\Big)\Big| \, dt\\
        & \ge & \frac{\db}{2} \frac{\eps(2-\eps)}{4} 
                \und{\inf}{u\in\big[\frac{1-\eps \sqrt{\db}}{2};\frac{1+\eps \sqrt{\db}}{2}\big]} 
                \big|\phi''_{\ha,\hb}(u)\big|\\
        & \ge & \frac{\eps(2-\eps)}{8} \db 
                \times \Big| \phi''_{-B,B}\Big(\frac{1-\eps\sqrt{\db}}{2}\Big) \Big|\\
        & \ge & \frac{\eps(2-\eps)}{8} \db
                \times \frac{q}{q-1} \big[\frac{1-\eps^2\db}{4}\big]^{\frac{2-q}{q-1}}
                \frac{(2B)^q}{2^{q+1} [(1+\eps\sqrt{\db})/2]^{\frac{q+1}{q-1}}}\\
        & \ge & \frac{\eps(2-\eps)}{8} \db 
                \times \frac{4 q B^q}{q-1} (1-\eps\sqrt{\db})^{\frac{2-q}{q-1}}
                (1+\eps\sqrt{\db})^{\frac{1-2q}{q-1}}\\
        & = & (1-\eps/2) q B^q \frac{\eps\db}{q-1} 
                (1-\eps\sqrt{\db})^{\frac{2-q}{q-1}}
                (1+\eps\sqrt{\db})^{\frac{1-2q}{q-1}}\\
        \endar
Let $K = (1-\eps\sqrt{\db})^{\frac{2-q}{q-1}} (1+\eps\sqrt{\db})^{\frac{1-2q}{q-1}}$.
From \eqref{eq:w1a}, taking $\tw=1/\tm$, we get
        \begarlab{eq:cur3}
        \underset{P \in \H}{\sup} \big\{ \E R( \hg ) 
                - \und{\min}{g} R(g) \big\} & \ge & 
                (1-\eps/2) KqB^q \frac{\eps\db}{q-1} \big(1-\sqrt{n\db/\tm}\big).
        \endarlab
This leads us to choose $\db=\frac{\tm}{4n} \wedge 1$ and $\eps=(q-1)\sqrt{\frac{n}{\tm}\vee \frac{1}{4}} \le \demi$
and obtain
        \begar
        \E R( \hg ) - \undc{\min}{g} R(g) 
                \ge \frac{3 q B^q}{8} K \big\{ 
                \big( \frac{1}{4}\sqrt{\frac{\tm}{n}} \big) \vee \big( 1 - \sqrt{\frac{n}{\tm}} \big) \big\}.
        \endar 
Since $1<q \le 2$ and $\eps\sqrt{\db}=\frac{q-1}{2},$ we may check that $K\ge 0.29$
(to be compared with $\lim_{q\rightarrow 1} K = e^{-1}\approx 0.37$).

\paragraph{Case $q> 1+\sqrt{\frac{\tm}{4n}}$ : }
%
We take $\tw=\frac{1}{n+1} \wedge \frac{1}{\tm}.$
From \eqref{eq:psidef}, \eqref{eq:defed1a} and \eqref{eq:qphi}, we get
        $
        \da=\psi_{1,0,-B,B}(1/2)
        =\phi_{-B,B}(1/2)=B^{q}.
        $
From \eqref{eq:w1c}, we obtain
        \begarlab{eq:forallq}
        \E R(\hg) - \und{\min}{g\in\G} R(g) 
                & \ge & \big( \frac{\lfloor \log_2 |\G| \rfloor}{n+1} \wedge 1\big) B^q
                                \big(1-\frac{1}{n+1} \wedge \frac{1}{\lfloor \log_2 |\G| \rfloor}\big)^n\\
                & \ge & e^{-1}B^q\big( \frac{\lfloor \log_2 |\G| \rfloor}{n+1} \wedge 1\big),
        \endarlab
where the last inequality uses $[1-1/(n+1)]^{n} \searrow e^{-1}$.

\paragraph{Improvement when $1+\sqrt{\frac{\tm}{4n}\wedge 1}< q < 2$ : }
%
From \eqref{eq:cur3}, by choosing $\eps=1/2$ and introducing $K'\eqdef (1-\sqrt{\db}/2)^{\frac{2-q}{q-1}}
                (1+\sqrt{\db}/2)^{\frac{1-2q}{q-1}}$,
we obtain
        \begar
        \underset{P \in \H}{\sup} \big\{ \E R( \hg ) 
                - \und{\min}{g} R(g) \big\} & \ge & 
                \frac{3qB^q}{8} K'  \frac{\db}{q-1} \big(1-\sqrt{n\db/\tm}\big).
        \endar
This leads us to choose $\db=\frac{4\tm}{9n} \wedge 1$. 
Since $\sqrt{\frac{\tm}{4n}\wedge 1}<q-1$, we have $\sqrt\db \le \frac{4}{3}(q-1)$, hence
$K' \ge \big(1-\frac{2}{3}(q-1)\big)^{\frac{2-q}{q-1}}
                \big(1+\frac{2}{3}(q-1)\big)^{\frac{1-2q}{q-1}}$. For any
$1<q <2$, this last quantity is greater than $0.2$.
So we have proved that for $1+\sqrt{\frac{\tm}{4n}\wedge 1}< q < 2$,
        \begarlab{eq:qinf2}
        \E R(\hg) - \und{\min}{g\in\G} R(g) 
                \ge \frac{q}{90(q-1)} B^q \frac{\lfloor \log_2 |\G| \rfloor}{n} .
        \endarlab
Theorem \ref{th:qsup1} follows from \eqref{eq:forallq} and \eqref{eq:qinf2}.

\secproofth{lqunbounded}

\subsubsection{Proof of the first inequality of Theorem \ref{th:lqunbounded}.}

Let $\tm=\lfloor \log_2 |\G| \rfloor$.
Contrary to other lower bounds obtained in this work, this learning setting requires asymmetrical
hypercubes of distributions.
Here we consider a constant $\tm$-dimensional hypercube of distributions with edge probability $\tw$ such that
$\pa \equiv p$, $\pb \equiv 0$, $\ha\equiv +B$ and $\hb\equiv 0$,
where $\tw$, $p$ and $B$ are positive real parameters to be chosen
according to the strategy described at the beginning of Section \ref{sec:examples}.
To have $\E |Y|^s \le A$, we need that 
$\tm \tw p B^s \le A$.
To ensure that a best prediction function has infinite norm bounded by $b$, from 
the computations at the beginning of Appendix \ref{app:proofphisecond}, we need that
        \begar
        B \le \frac{p^{1/(q-1)}+(1-p)^{1/(q-1)}}{p^{1/(q-1)}} \, b.
        \endar
This inequality is in particular satisfied for $B = C p^{-1/(q-1)}$ for appropriate
small constant $C$ depending on $b$ and $q$.
From the definition of the edge discrepancy of type II, we have $\db=p$.
In order to have the $\rhs$ of \eqref{eq:w1a} of order $mw\da$, we want to have 
$n \tw p \le C < 1$. 
All the previous constraints lead us to take the parameters 
$\tw, p$ and $B$ such that
        \lbegar
        B = C p^{-1/(q-1)}\\
        \tm \tw p B^s= A\\
        n \tw p = 1/4
        \rendarp

Let $Q=\frac{\tm}{n}\wedge 1$.
This leads to
        $p = C Q^{(q-1)/s}$,
        $B = C Q^{-1/s}$
and
        $\tw= C \tm^{-1} Q^{1-(q-1)/s}$
with $C$ small positive constants depending on $b$, $A$, $q$ and $s$.
Now from the definition of the edge discrepancy of type I and 
\eqref{eq:psidev}, we have
        \begar
        \da & = & \frac{p^2}{2} \int_0^1 [t\wedge(1-t)] 
                \big| \phi''_{0,B}(tp) \big| dt \\
        & \ge & \frac{p^2}{2} \int_{1/4}^{3/4} \frac{1}{4}
                \min_{[p/4;3p/4]} \big|\phi''_{0,B}(tp)\big| dt\\
        & \ge & C p^2 p^{\frac{2-q}{q-1}} B^q\\
        & = & C
        \endar
where the last inequality comes from \eqref{eq:qpower}.
From \eqref{eq:w1a}, we get
        \begar
        \underset{P \in \calP}{\sup} \big\{ \E R( \hg ) 
                - \und{\min}{g\in\G} R(g) \big\} 
                \ge C Q^{1-\frac{q-1}{s}}.
        \endar

\subsubsection{Proof of the second inequality of Theorem \ref{th:lqunbounded}.} \label{sec:proofsecond}

We still use $\tm=\lfloor \log_2 |\G| \rfloor$.
We consider a $(\tm,\tw,\tdb)$-hypercube with $\ha\equiv -B$ and $\hb\equiv +B$,
where $\tw, \tdb$ and $B$ are positive real parameters to be chosen
according to the strategy described at the beginning of Section \ref{sec:examples}.
To have $\E |Y|^s \le A$, we need that 
$\tm \tw B^s \le A$.
To ensure that a best prediction function has infinite norm bounded by $b$, from 
the computations at the beginning of Appendix \ref{app:proofphisecond}), we need that
        \begarlab{eq:bandb}
        B \le \frac{[1+(\tdb)^{1/2}]^{1/(q-1)}+[1-(\tdb)^{1/2}]^{1/(q-1)}}
                {[1+(\tdb)^{1/2}]^{1/(q-1)}-[1-(\tdb)^{1/2}]^{1/(q-1)}} \, b.
        \endarlab
For fixed $q$ and $b$, this inequality essentially means that $B \le C \tdb^{-1/2}$ since we intend to take 
$\tdb$ close to $0$. 
In order to have the $\rhs$ of \eqref{eq:w1a} of order $mw\da$, we want to have 
$n \tw \tdb \le 1/4$ where, once more, this last constant is arbitrarily taken. 
The previous constraints lead us to choose
        \lbegar
        B = C \tdb^{-1/2}\\
        \tm \tw B^s= A\\
        n \tw \tdb = 1/4
        \rendarp
We still use $Q=\frac{\tm}{n}\wedge 1$.
This leads to
        $\tdb = C Q^{2/(s+2)}$,
        $B = C Q^{-1/(s+2)}$
and
        $\tw= C \tm^{-1} Q^{s/(s+2)}$
with $C$ small positive constants depending on $b$, $A$, $q$ and $s$.
Now from \eqref{eq:qpower}, the differentiability assumption is satisfied
for $\zeta=C B^q = C Q^{-q/(s+2)}$.
From \eqref{eq:w1a} and \eqref{eq:dadb}, we obtain
        \begar
        \underset{P \in \calP}{\sup} \big\{ \E R( \hg ) 
                - \und{\min}{g\in\G} R(g) \big\} 
                \ge C Q^{1-\frac{q}{s+2}}.
        \endar

\secproofth{l2ub}

The starting point is similar to the one in Section \ref{sec:proofsecond}. Since $q=2$, \eqref{eq:bandb}
simplifies into $B\le b (\tdb)^{-1/2}$. We take $B=b (\tdb)^{-1/2}$ and $\tw=A/(\tm B^2)$
and we optimize the parameter $\tdb$ in order to maximize the lower bound.
From \eqref{eq:ls}, we get $\tm \tw \da = A \tdb$.
Introducing $a\eqdef n \tw \tdb = \frac{nA}{\tm b^2} (\tdb)^2,$ we obtain
$\tm \tw \da = b\sqrt{A \tm/n} \sqrt{a}$. The results then follow from Corollary \ref{cor:appassouad}
and the fact that the differentiability assumption \eqref{eq:assump1}
holds for $\zeta = 8B^2 = \frac{8\da}{\db}.$

\fi


\ifRR
\appendix


\section{Computations of the second derivative of $\phi$ for the $L_q$-loss} \label{app:proofphisecond}

Let $\ya$ and $\yb$ be fixed.
We start with the computation of $\phi_{\ya,\yb}$.
For any $p\in[0;1]$, the quantity 
        $\varphi_{p,y_1,y_2}(y) = p|y-\ya|^q+(1-p)|y-\yb|^q$
is minimized when $y\in[\ya\wedge\yb;\ya\vee\yb]$ and
        $pq (y-\ya)^{q-1} = (1-p)q(\yb-y)^{q-1}$.
Introducing $r=\frac{1}{q-1}$ and $D=p^r+(1-p)^r$, the minimizer can be written as
        $y=\frac{p^r \ya + (1-p)^r \yb}{D}$ and the minimum is
        \begar
        \phi_{\ya,\yb}(p) & = & \Big(p \frac{(1-p)^{rq}}{D^q} + (1-p) \frac{p^{rq}}{D^q} \Big) |\yb-\ya|^q\\
        & = & p(1-p) \frac{|\yb-\ya|^q}{D^{q-1}},
        \endar
where we use the equality $rq=1+r$.
We get
        \begar
        \frac{1}{|\yb-\ya|^q} \phi'_{\ya,\yb}(p) & = & 
                \frac{1-2p}{D^{q-1}}+p(1-p)(1-q)r D^{-q} [p^{r-1}-(1-p)^{r-1}]\\
        & = & D^{-q} \big\{ (1-2p)[p^r+(1-p)^{r}]-(1-p)p^r+p(1-p)^{r}] \big\}\\
        & = & D^{-q} \big\{ (1-p)^{r+1}-p^{r+1} \big\},
        \endar
hence
        \begar
        \frac{1}{|\yb-\ya|^q} \phi''_{\ya,\yb}(p) & = & 
                -qr D^{-q-1}[p^{r-1}-(1-p)^{r-1}][(1-p)^{r+1}-p^{r+1}]\\
        & & 
                -qr D^{-q-1}[p^{r}-(1-p)^{r}]^2\\
        & = & -qr D^{-q-1} p^{r-1}(1-p)^{r-1}\\
        & = & - \frac{q}{q-1} 
                \frac{[p(1-p)]^{\frac{2-q}{q-1}}}{\big[p^{\frac{1}{q-1}}+(1-p)^{\frac{1}{q-1}}\big]^{q+1}}.
        \endar

\section{Expected risk bound from Hoeffding's inequality} \label{app:stdhoeff}

Let $\lam'>0$ and $\rho$ be a probability distribution on $\G$. Let $r(g)$ denote the empirical 
risk of a prediction function $g$, that is
        $
        r(g) = \frac{1}{n} \sum_{i=1}^n L(Z_i,g).
        $
Hoeff\-ding's inequality applied to the random variable $W=\expec{g}{\rho} L(Z,g)-L(Z,g') \in [-(b-a);b-a]$
for a fixed $g'$ gives
        \begar
        \undc{\E}{Z\sim P} e^{\eta[W-\E W]} \le e^{\eta^2 (b-a)^2/2} 
        \endar
for any $\eta>0$. For $\eta=\lam'/n$, this leads to
        \begar
        \Ezun e^{\lam'[R(g') - \expec{g}{\rho} R(g) - r(g') + \expec{g}{\rho} r(g)]} \le e^{(\lam')^2 (b-a)^2/(2n)}
        \endar
Consider the Gibbs distribution $\hrho=\pi_{-\lam' r}$. This distribution 
satisfies
        \begar
        \expec{g'}{\hrho} r(g') + K(\hrho,\pi)/\lam' \le \expec{g}{\rho} r(g) + K(\rho,\pi)/\lam'.
        \endar
We have
        \begar
        & \Ezun \expec{g'}{\hrho} R(g') -\expec{g}{\rho} R(g)\\
        \le & \Ezun \Big\{ \expec{g'}{\hrho}\big[ R(g')-\expec{g}{\rho} R(g)
                - r(g')-\expec{g}{\rho} r(g) \big] + \frac{K(\rho,\pi)-K(\hrho,\pi)}{\lam'} \Big\} \\
        \le & \frac{K(\rho,\pi)}{\lam'}+\Ezun \frac{1}{\lam'}\log \undc{\E}{g'\sim\pi} e^{\lam'
                [ R(g')-\expec{g}{\rho} R(g)
                - r(g')-\expec{g}{\rho} r(g) ]}\\
        \le & \frac{K(\rho,\pi)}{\lam'}+\frac{1}{\lam'}\log \undc{\E}{g'\sim\pi} \Ezun e^{\lam'
                [ R(g')-\expec{g}{\rho} R(g)
                - r(g')-\expec{g}{\rho} r(g) ]}\\
        \le & \frac{K(\rho,\pi)}{\lam'}+\frac{\lam' (b-a)^2}{2n}.
        \endar
This proved that for any $\lam>0$, the 
generalization error of the algorithm which draws its prediction function 
according to the Gibbs distribution $\pi_{-\lam \Sigma_n/2}$ satisfies 
        \begar
        \Ezun \expec{g'}{\pi_{-\lam \Sigma_n/2}} R(g')
                \le \und{\min}{\rho\in\M} \Big\{ 
                \expec{g}{\rho} R(g) 
                + 2 \Big[ \frac{\lam (b-a)^2}{8} 
                + \frac{K(\rho,\pi)}{\lam n} \Big]\Big\},
        \endar
where we use the change of variable $\lam=2\lam'/n$ in order to underline the
difference with \eqref{eq:corhoeff1}.

\else
\fi

\ifRR

\section{Proof of Inequality \eqref{eq:swedgen3}} \label{app:proofasympt}

To prove \eqref{eq:swedgen3}, we need to uniformly control the difference between the tail
of the sum of i.i.d. random variables and the gaussian approximate. This is done by the following result.
\begin{theorem}[Berry\cite{Ber41}-Esseen\cite{Ess42} inequality] \label{th:petrov}
Let $N$ be a centered gaussian variable of variance $1$.
Let $U_1,\dots,U_n$ be real-valued independent identically distributed random variables such that $\E U_1=0$,
$\E U_1^2=1$ and $\E |U_1|^{3} < +\infty$. Then
        \begarlab{eq:petrov}
        \und{\sup}{x\in\R} \Big|\P\big(n^{-1/2} \sum_{i=1}^n U_i > x \big) - \P( N > x ) \Big|
                \le C n^{-1/2} \E |U_1|^{3}
        \endarlab
for some universal positive constant $C$.

\end{theorem}

To shorten the notation, let $P^-=P_{[-]}^{\otimes n}$ and $P^+=P_{[+]}^{\otimes n}$ 
be the $n$-fold product of the representatives of the hypercube.
Since we have $\pa>1/2>\pb$ (by definition of a $(m,w,\db)$-hypercube (p.\pageref{def:mwxi})),
the set of sequences $Z_1^n$ for which 
$\frac{P^+}{P^-}(Z_1^n)<1$ is 
        \begar
        E \eqdef \big\{ \sum_{i=1}^n \ds1_{Y_i=\ha(X_i),X_i\in\X_1} < \sum_{i=1}^n \ds1_{Y_i=\hb(X_i),X_i\in\X_1} \big\}
        \endar
Introduce the quantity
        \begar
        S_n \eqdef \sum_{i=1}^n \big( 2\ds1_{Y_i=\ha(X_i)}-1 \big) \ds1_{X_i\in\X_1}.
        \endar
We have
        \begarlab{eq:interm3}
        \S_{\wedge} ( P^+ , P^- )
                & = & 1-P^-( E ) + P^+( E )\\
                & = & 1-P^-( S_n < 0) + P^+( S_n < 0)\\
        & = & P^-( S_n \ge 0) + P^-( S_n > 0).
        \endarlab
Introduce 
        \begar
        W_i=\big( 2\ds1_{Y_i=\ha(X_i)}-1 \big) \ds1_{X\in\X_1}.
        \endar
From now on, we consider that the pairs $Z_i=(X_i,Y_i)$ are generated by $P^-$, so that
$\E W_i$ and $\Var W_i$ simply denote the expectation and variance of $W_i$ 
when $(X_i,Y_i)$ is drawn according to $P_{-1,1,\dots,1}$.
Define the normalized quantity 
        \begar
        U_i \eqdef (W_i-\E W_i)/\sqrt{\Var W_i}.
        \endar
We have         
        \begar
        P^-( S_n > 0) = P^-\big(n^{-1/2} \sum_{i=1}^n U_i > t_n \big)
        \endar
where $t_n\eqdef-\sqrt{\frac{n}{\Var W_1}}\E W_1$.
By Berry-Esseen's inequality (Theorem \ref{th:petrov}), we get
        \begarlab{eq:pet2}
        \big| P^-( S_n > 0) - \P( N > t_n ) \big| \le C n^{-1/2} \E |U_1|^{3}
        \endarlab
Let us now upper bound $n^{-1/2} \E |U_1|^{3}$.

Since we have $\pa=(1+\xi)/2=1-\pb$ for a $(m,w,\db)$-hypercube (p.\pageref{def:mwxi})),
the law of $W_1$ is described by
        \lbegar
        \P(W_1=1)= w \frac{1-\xi}{2}\\
        \P(W_1=0) = 1-w\\
        \P(W_1=-1) = w \frac{1+\xi}{2}
        \rendar
where $w$ still denotes $\mu(\X_1)$.
We get $\E W_1=-w \xi$, $\Var W_1 = w(1-w\xi^2)$ and
since $0<w\le 1$ and $0<\xi \le 1$
        \begarlab{eq:gammas}
        \E |W_1-\E W_1|^3 & = & (1-w) (w\xi)^3+w \frac{1-\xi}{2} (1+w\xi)^3
                + w \frac{1+\xi}{2} (1-w\xi)^3\\
        & \le & w + w[1+3(w\xi)^2]\\
        & \le & 5 w.
        \endarlab
We obtain
        \begar
        n^{-1/2} \E|U_1|^3 \le 
                n^{-1/2} \frac{5 w}{[w(1-w\xi^2)]^{3/2}} \le 5 \frac{(nw)^{-1/2}}{(1-\xi^2)^{3/2}}
                \und{\longrightarrow}{nw\rightarrow+\infty,\xi\rightarrow 0} 0.
        \endar
From \eqref{eq:pet2}, we get that 
        $\big| P^-( S_n > 0) - \P( N > t_n ) \big|$ 
converges to zero when $nw$ goes to infinity and $\xi$ goes to zero.
Now the previous convergence also holds when `$>$' is replaced with '$\ge$' (since it suffices
to consider the random variables $-U_i$ in Theorem \ref{th:petrov}).
Using both convergence and \eqref{eq:interm3}, we obtain
        \begar
        \big| \S_{\wedge} ( P^+ , P^- ) - \P( |N| > t_n ) \big|
        \und{\longrightarrow}{nw\rightarrow+\infty,\xi\rightarrow 0} 0,
        \endar
which is the desired result since 
        $t_n=-\sqrt{\frac{n}{\Var W_1}}\E W_1=\sqrt{\frac{nw\xi^2}{1-w\xi^2}}=\sqrt{\frac{nw\db}{1-w\db}}$
and
	\begar
	\Big| \P\Big( |N| > \sqrt{\frac{nw\db}{1-w\db}} \Big) - \P( |N| > \sqrt{nw\db} ) \Big| \und{\longrightarrow}{nw\rightarrow+\infty,\db\rightarrow 0} 0.
	\endar

\section{Towards adaptivity for the temperature parameter} \label{sec:ada}

%

Once the distribution $\pi$ is fixed, an appropriate choice for the parameter $\lam$ of the SeqRand algorithm is the minimizer of 
the $\rhs$ of \eqref{eq:1gen}. This minimizer is unknown by the statistician. This section proposes to
modify the $\lam$ during the iterations so that it automatically fits to a value close to this minimizer.
The adaptive SeqRand algorithm is described in Figure~\ref{fig:ada}. The idea of incremental 
updating of the temperature parameter to solve the adaptivity problem has been successfully developed in
\cite[Section 2]{AueCesGen02} and \cite[Lemma 3]{Ces07}. Here we improve the argument by using Lemma \ref{lem:partition}.

\begin{figure} 
\hspace*{0.5cm}\hbox{\raisebox{0.4em}{\vrule depth 0pt height 0.4pt width 12cm}}
\begin{enumerate}
\item[] Input: \begin{enumerate} \item[$\bullet$] $\lam_1\ge \lam_2\ge \cdots\ge \lam_{n+1} >0$ with $\lam_{i+1}$ possibly depending on the values of $Z_1,\dots,Z_i$
\item[$\bullet$] $\pi$ a distribution on the set $\G$
\end{enumerate}
\item Define $\hrho_0 \eqdef \hpi(\pi)$ in the sense of the variance inequality (p.\pageref{varcond})
and draw a function $\hg_{0}$ according to this distribution.
Let $S_0(g)=0$ for any $g\in\G$. 
\item For any $i\in\{1,\dots,n\}$, iteratively define 
        \begarlab{eq:adaSi}
        S_{i}(g) \eqdef S_{i-1}(g) + L(Z_{i},g) + \delta_{\lam_i}(Z_{i},g,\hg_{i-1}) \quad \text{for any } g\in\G.
        \endarlab
and
        \begar
        \hrho_i \eqdef \hpi(\pi_{-\lam_{i+1} S_i}) \qquad\text{in the sense of the variance inequality (p.\pageref{varcond})}
        \endar
and draw a function $\hg_i$ according to the distribution $\hrho_i$.
\item Predict with a function drawn according to the uniform distribution on the finite set $\{\hg_0,\dots,\hg_n\}$.\\
Conditionally to the training set, the distribution of the output prediction function will 
be denoted $\hmua$.
\end{enumerate}
\hspace*{0.5cm}\hbox{\raisebox{0.4em}{\vrule depth 0pt height 0.4pt width 12cm}}
\caption{The adaptive SeqRand algorithm} \label{fig:ada}
\end{figure}

The following theorem upper bounds the generalization error of the adaptive SeqRand algorithm.
\begin{theorem} \label{th:ada}
Let $\Delta_\lam(g,g') \eqdef \undc{\E}{Z\sim P} \delta_{\lam}(Z,g,g')$ for $g \in G$ and $g'\in \barG$,
where we recall that $\delta_{\lam}$ is a function satisfying the variance inequality (see p.\pageref{varcond}). 
The expected risk of the adaptive SeqRand algorithm satisfies
        \begarlab{eq:ada1gen}
        \Ezun \expec{g'}{\hmua} R(g')
                & \le & \und{\min}{\rho\in\M} \bigg\{ 
                \expec{g}{\rho} R(g) 
                + \frac{K(\rho,\pi)}{\lam_{n+1}(n+1)} \\
                & & \qquad + \expec{g}{\rho}  \Ezun \expecc{\hg_0^n}{\Omega_n} \frac{\sum_{i=0}^{n} \Delta_{\lam_{i+1}}(g,\hg_i)}{n+1}  
                \bigg\}
        \endarlab
\end{theorem}

\begin{proof}
Let $\calE$ denote the expected risk of the adaptive SeqRand algorithm:
        \begar 
        \calE \eqdef \Ezun \expec{g}{\hmua} R(g)
                 = \frac{1}{n+1} \sum_{i=0}^{n} \E_{Z_1^i} \expec{\hg_0^i}{\Omega_i} R(\hg_i).
        \endar
We recall that $Z_{n+1}$ is a random variable 
independent of the training set $Z_1^n$ and with the same distribution $P$.
Define $S_{n+1}$ by \eqref{eq:adaSi} for $i=n+1$.
To shorten formulae, let $\hpi_i \eqdef \pi_{-\lam_{i+1} S_i}$ so that by definition we have
$\hrho_i=\hpi(\hpi_i)$. The variance inequality implies that
        \begar
        \expec{g'}{\pirho} R(g') \le - \frac{1}{\lam_{i+1}} 
                \E_Z \expec{g'}{\pirho} \log \expec{g}{\rho} e^{-\lam_{i+1} [L(Z,g)+\delta_{\lam_{i+1}}(Z,g,g')]}.
        \endar
So for any $i\in\{0,\dots,n\},$ for fixed $\hg_0^{i-1}=(\hg_0,\dots,\hg_{i-1})$
and fixed $Z_1^i$, we have
        \begar
        \expe{g'}{\hrho_i} R(g') \le - \frac{1}{\lam_{i+1}} 
                \E_{Z_{i+1}} \expec{g'}{\hrho_i} \log \expec{g}{\hpi_i} 
                e^{-\lam_{i+1} [L(Z_{i+1},g)+\delta_{\lam_{i+1}}(Z_{i+1},g,g')]}
        \endar
Taking the expectations $\wrt$ $(Z_1^i,\hg_0^{i-1})$, we get
        \begar
        \E_{Z_1^i} \expecc{\hg_0^i}{\Omega_i} R(\hg_i)
                & = & \E_{Z_1^i} \expecc{\hg_0^{i-1}}{\Omega_{i-1}} \expec{g'}{\hrho_i} R(g')\\
                & \le & - \E_{Z_1^{i+1}} \expecc{\hg_0^i}{\Omega_i} \frac{1}{\lam_{i+1}} 
                \log \expec{g}{\hpi_i} e^{-\lam_{i+1} [L(Z_{i+1},g)+\delta_{\lam_{i+1}}(Z_{i+1},g,\hg_i)]}.
        \endar
Consequently, by the chain rule (i.e. cancellation in the sum of logarithmic terms; \cite{Bar87}) 
and by intensive use of Fubini's theorem, we get        
        \begarlab{eq:aa}
        \calE = \frac{1}{n+1} \sum_{i=0}^{n} \E_{Z_1^i} \expecc{\hg_0^i}{\Omega_i} R(\hg_i)\\
        \, \le - \frac{1}{n+1} \sum_{i=0}^{n}
                \E_{Z_1^{i+1}} \expecc{\hg_0^i}{\Omega_i} \frac{1}{\lam_{i+1}} \log \expec{g}{\hpi_i} e^{-\lam_{i+1} 
                [ L(Z_{i+1},g)+\delta_{\lam_{i+1}}(Z_{i+1},g,\hg_i) ]} \\
        \, = \frac{1}{n+1} \E_{Z_1^{n+1}} \expecc{\hg_0^n}{\Omega_n} \sum_{i=0}^{n} a_i,
        \endarlab
with $a_i = - \frac{1}{\lam_{i+1}} \log \expe{g}{\hpi_i} e^{-\lam_{i+1} [L(Z_{i+1},g)+\delta_{\lam_{i+1}}(Z_{i+1},g,\hg_i)]}$.
Introduce the function $\phi_i(\lam) = \frac{1}{\lam} \log \expec{g}{\pi} e^{-\lam S_{i}(g)}.$
Let us now concentrate on the last sum.
		\begarlab{eq:ab}
		\sum_{i=0}^{n} a_i & = & - \sum_{i=0}^{n} \frac{1}{\lam_{i+1}} \log \left( \frac{\expec{g}{\pi} e^{-\lam_{i+1} S_{i+1}(g)}}
                {\expec{g}{\pi} e^{-\lam_{i+1} S_{i}(g)}} \right)\\
		& = & - \frac{1}{\lam_{n+1}} \log \expec{g}{\pi} e^{-\lam_{n+1} S_{n+1}(g)}
                + \sum_{i=0}^{n} \left[ \phi_i(\lam_{i+1}) - \phi_i(\lam_{i})  \right]\\
        & \le & - \frac{1}{\lam_{n+1}} \log \expec{g}{\pi} e^{-\lam_{n+1} S_{n+1}(g)},
        \endarlab
where the last inequality uses that $\lam_{i+1} \le \lam_i$ and that the functions $\phi_i$ are nondecreasing according to
the following lemma.

\begin{lemma} \label{lem:partition}
Let $\W$ be a real-valued measurable function defined on the space $\A$ and
let $\mu$ be a probability distribution on $\A$. The mapping
$\lam \mapsto \frac{1}{\lam} \log \expec{a}{\mu} e^{-\lam \W(a)}$
is nondecreasing on the interval on which it is defined.
\end{lemma}
\begin{proof}[Proof of Lemma \ref{lem:partition}]
It is easy to check that the function is well-defined on an interval (possibly empty),
and is smooth on the interior of this interval. The derivative of this function is
	$$
	- \frac{1}{\lam^2} \log \expecc{a}{\mu} e^{-\lam \W(a)} - \frac{1}{\lam} \frac{\expecc{a}{\mu} \W(a) e^{-\lam \W(a)}}{\expecc{a}{\mu} e^{-\lam \W(a)}}
		= \frac{1}{\lam^2} K( \mu_{- \lam \W(a)} , \mu )\ge 0.
	$$
\end{proof}

Plugging the inequality \eqref{eq:ab} into \eqref{eq:aa} and using Lemma \ref{eq:concpart}, we obtain
        \begar
        \calE & \le & - \frac{1}{n+1} \expecc{Z_1^{n+1}}{P^{\otimes (n+1)}} 
                \expecc{\hg_0^n}{\Omega_n} \frac{1}{\lam_{n+1}} \log \expec{g}{\pi} e^{-\lam_{n+1} S_{n+1}(g)}\\
        & = & \frac{1}{n+1} \expecc{Z_1^{n+1}}{P^{\otimes (n+1)}} \expecc{\hg_0^n}{\Omega_n} \und{\min}{\rho\in\M} \left\{ 
                \expec{g}{\rho} \sum_{i=1}^{n+1} \big[ L(Z_i,g) + \delta_{\lam_{i}}(Z_i,g,\hg_{i-1}) \big] 
                + \frac{K(\rho,\pi)}{\lam_{n+1}} \right\}\\
        & \le & \und{\min}{\rho\in\M} \left\{ 
                \expec{g}{\rho} R(g) + \expec{g}{\rho} \expecc{Z_1^{n}}{P^{\otimes (n+1)}} \expecc{\hg_0^n}{\Omega_n} 
                \frac{1}{n+1} \sum_{i=1}^{n+1} \Delta_{\lam_{i}}(g,\hg_{i-1}) 
                + \frac{K(\rho,\pi)}{\lam_{n+1}(n+1)} \right\}\\
        \endar  
\end{proof}

\vspace{0.5cm}
\noindent {\bf Acknowledgement.} 
I would like to thank Nicolas Vayatis, Alexandre Tsybakov, Gilles Stoltz and Olivier Catoni
for their helpful comments.
\fi

\bibliographystyle{plain}
\ifHAL
\bibliography{ref_hal}
\else
\bibliography{../../ref}
\fi
\end{document}